\documentclass[10pt]{article}
\usepackage{setspace}

\usepackage{amsmath} 
\usepackage{amssymb}
\usepackage{latexsym}
\usepackage{xcolor}
\usepackage{fancyhdr}
\allowdisplaybreaks 

 \usepackage{comment}

\def\Xint#1{\mathchoice {\XXint\displaystyle\textstyle{#1}} 
{\XXint\textstyle\scriptstyle{#1}}
{\XXint\scriptstyle\scriptscriptstyle{#1}} 
{\XXint\scriptscriptstyle\scriptscriptstyle{#1}}\!\int} 
\def\XXint#1#2#3{{\setbox0=\hbox{$#1{#2#3}{\int}$} 
\vcenter{\hbox{$#2#3$}}\kern-.5\wd0}}   

 \numberwithin{equation}{section}
\newtheorem{theorem}[equation]{Theorem}
\newtheorem{proposition}[equation]{Proposition}
\newtheorem{definition}[equation]{Definition}
\newtheorem{remark}[equation]{Remark}
\newtheorem{lemma}[equation]{Lemma}
\newtheorem{corollary}[equation]{Corollary}

\title{
A nonvariational form of the Neumann  problem for H\"{o}lder continuous harmonic functions
 } 
 
\author{  
Massimo Lanza de Cristoforis
\\
Dipartimento di Matematica `Tullio Levi-Civita', 
\\
Universit\`a degli Studi di Padova, 
\\
Via Trieste 63, Padova 35121, 
Italy. 
\\
E-mail: mldc@math.unipd.it   }

\date{\ }

\begin{document}
 
 \maketitle

\noindent
{\bf Abstract:}  We present a nonvariational setting for the Neumann problem for harmonic functions that are H\"{o}lder continuous and that may have infinite Dirichlet integral. Then we introduce a space of distributions on the boundary that we characterize (a  space of first order traces for H\"{o}lder continuous harmonic functions), we analyze the properties of the corresponding distributional single layer potential and  we prove a representation theorem for harmonic H\"{o}lder continuous functions in terms of distributional single layer potentials. As an application, we solve  the interior and exterior Neumann problem with distributional data in the space of first order traces that has been introduced.
 \vspace{\baselineskip}

\noindent
{\bf Keywords:} distributional normal derivative, nonvariational Neumann problem,  first order traces, H\"{o}lder continuous functions, integral representation, harmonic functions

\par
\noindent   
{{\bf 2020 Mathematics Subject Classification:}}   31B20, 31B10, 
35J25, 47B92.

\section{Introduction} Our starting point are the classical examples of Prym \cite{Pr1871} and  Hadamard \cite{Ha1906} of harmonic functions in the unit ball of the plane that are continuous up to the boundary and have infinite energy, \textit{i.e.},  whose
 gradient is not square-summable. Such functions solve the classical Dirichlet problem in the unit ball, but not the corresponding weak (variational) formulation. For a discussion on this point we refer to 
Maz’ya and Shaposnikova \cite{MaSh98}, Bottazzini and Gray \cite{BoGr13} and
 Bramati,  Dalla Riva and  Luczak~\cite{BrDaLu23}. 
 
 We plan to define a normal derivative on the boundary  in distributional form  for H\"{o}lder continuous harmonic functions 
 (cf.~Definitions \ref{defn:conoderd}, \ref{defn:conoederd}) and to introduce a nonvariational form  for the corresponding Neumann problem. Then our main results consist in the Definition \ref{defn:njsdoa} and in the characterization Theorems \ref{thm:v-1achar}, \ref{thm:v-1an} of the space of traces for the distributional normal derivatives of H\"{o}lder continuous harmonic functions, in the representation Theorem  \ref{thm:isoslv-1a}  for H\"{o}lder continuous harmonic functions in terms of   (distributional) single layer potentials and in showing that one can develop the layer potential theoretic method to show the  existence Theorems \ref{thm:nvINexist}, \ref{thm:nvENexist3}, \ref{thm:nvENexist2}  for the solution of the (nonvariational) interior and exterior Neumann problems. Instead, we plan to consider the extension of the arguments of the present paper to the solution of the Poisson equation in a forthcoming paper.

 Our approach holds in the (nonseparable)  spaces of H\"{o}lder continuous functions, but could be extended to different function spaces whose distributional gradient has no summability properties and differs from the so-called transposition method   of Lions and Magenes~\cite[Chapt.~II, \S 6]{LiMa68}, R\u{o}itberg and S\v{e}ftel' \cite{RoSe69}, 
 Aziz and Kellog \cite{AzKe95} which exploit  the form of the dual of a Sobolev space of functions  and which accordingly are suitable in a Sobolev space setting. 
 
 We also mention the distributional definition of the (difference) of normal derivatives of harmonic functions that has been exploited by Dalla Riva and Mishuris \cite[Defn.~4.8]{DaMi15} in order to define a transmission condition and that differs from the one adopted here for we replace the test functions by functions that are harmonic in the interior and in the exterior, a choice that serves our purposes of describing the space of first order traces of H\"{o}lder continuous harmonic functions.  
 
 {\color{black} In Section \ref{sec:dsila}, we introduce the  
  (distributional) single layer potential, some of its classical properties and 
we single out two conditions on  	 
 distributional moments $\mu$ 	}
 on the boundary of a bounded open subset $\Omega$ of class $C^{1,\alpha}$ of   the Euclidean space ${\mathbb{R}}^n$  for some $\alpha\in]0,1[$ {\color{black} such that one can consider the boundary value $V_\Omega[\mu]$ of the single layer potential 	}
(see Definition \ref{defn:Vpm}). 
{\color{black} Such conditions on $\mu$ require }   that the restrictions of the distributional single layer to   the sets
 \[
 \Omega^+\equiv \Omega \qquad\text{and}\qquad \Omega^-\equiv {\mathbb{R}}^n\setminus\overline{\Omega}
 \]
  can be extended   to  	 {\color{black}
   continuous	}
    functions on 
  the closures of $\Omega^+$ and $\Omega^-$, respectively and that 
  the two extensions are equal on $\partial\Omega$.
  
  In order to write the third Green identities in nonvariational form we must summarize some known properties of the classical double layer potential  in H\"{o}lder spaces (see Section  \ref{sec:hadolapo}),

As we explain in Section \ref{sec:conoderd}, we cannot exploit the first Green Identity in order to introduce a  distributional  normal derivative on the boundary as done by  
 Lions and Magenes~\cite{LiMa68}, Ne\v{c}as \cite[Chapt.~5]{Ne12},
 Nedelec and Planchard \cite[p.~109]{NePl73},  Costabel \cite{Co88}, McLean~\cite[Chapt.~4]{McL00}, Mikhailov~\cite{Mi11}, Mitrea, Mitrea and Mitrea \cite[\S 4.2]{MitMitMit22}. Instead, in Section \ref{sec:conoderd},  we introduce a distributional  normal derivative on the boundary for   harmonic functions in the interior by taking the transpose of the Steklov-Poincar\'{e} operator, \textit{i.e.}, the Dirichlet to Neumann map,  in a suitable duality pairing in the sense of Wendland \cite{We67}, \cite{We70} (cf.~Kress \cite{Kr14}). Then we do similarly for   harmonic functions in the exterior domain that are harmonic at infinity (see Section \ref{sec:conoederd}).

 Next in Sections \ref{sec:nvneupb} and \ref{sec:nveneupb} we introduce a nonvariational form of the interior and exterior Neumann problems for H\"{o}lder continuous functions. Here we note that even to prove known facts as the uniqueness up to locally constant functions, we cannot resort to energy arguments as in the variational setting.

  In Sections \ref{sec:pokein1a} and  \ref{sec:pokeex1a}, we present some known facts on the
   classical Green function for the Dirichlet problem with H\"{o}lder continuous solutions    both in $\Omega$ and in the exterior $\Omega^-$ in a form  that we need in the sequel.
 
 Next we turn to prove a representation theorem for the double layer potential  in terms of distributional single layer potentials (cf.~Section \ref{sec:thm:dlintesl}), that we exploit to write distributional forms of the third Green Identity for H\"{o}lder continuous functions (cf.~Section \ref{sec:3rdGreen0a}).

  By  the third Green Identity for H\"{o}lder continuous functions  of Theorem \ref{thm:3rdGreen0a}, we prove the  representation Theorem \ref{thm:slreth} for H\"{o}l\-der continuous harmonic functions in terms of single layer potentials of distributions in some auxiliary spaces of distributions $V^{-1,\alpha,\pm}(\partial\Omega)$  (cf.~Definition \ref{defn:njsdopm} of Section \ref{sec:slreth}).
 
 In Section \ref{sec:trdslay},  we show some technical statements that we later exploit to prove that  the auxiliary  spaces $V^{-1,\alpha,\pm}(\partial\Omega)$ are actually contained in the space $V^{-1,\alpha}(\partial\Omega)$ {\color{black} that we introduce below in  Section \ref{sec:ndslv}. }
  
 In Section \ref{sec:ndslv}, we prove a   jump formula for the distributional single layer potentials 
 {\color{black} of a class of distributions on the boundary that we define as 
 $V^{-1,\alpha}(\partial\Omega)$ in Definition \ref{defn:njsdoa}. The space $V^{-1,\alpha}(\partial\Omega)$ consists of distributions on the boundary such that the restriction  of the distributional single layer to   $\Omega$ can be extended   to a $\alpha$-H\"{o}lder continuous function on 
  the closure  of $\Omega$ and the restriction  of the distributional single layer to   $\Omega^-$    can be extended   to a locally $\alpha$-H\"{o}lder continuous function  on 
  the closure    of $\Omega^-$ and   such that 
  the two extensions are equal on $\partial\Omega$ and satisfy a certain symmetry condition of the distributional single layer potential   (see Proposition \ref{prop:ndslv0a}).	In}  Section \ref{sec:dcomvw}
 we prove distributional form of the classical Plemelj Symmetrization Principle.
 
 In Section \ref{sec:kpmi+wtv} we analyze the properties of the transpose of the
 operator  that corresponds to the trace of the double layer potential 
 in $C^{0,\alpha}(\partial\Omega)$ and we prove a representation theorem for 
 H\"{o}lder continuous functions in $\Omega$ and $\Omega^-$ by means of a distributional  single layer potential with density (or moment) in 
 $V^{-1,\alpha,\pm}(\partial\Omega)$ (see Theorem \ref{thm:isoslv-1a} and its Corollary \ref{corol:isoslvc-1a}).
 
 In Section \ref{sec:v-1achar}, we show that $V^{-1,\alpha}(\partial\Omega)=
V^{-1,\alpha,\pm}(\partial\Omega)$, an equality that we exploit to characterize and introduce a norm on $V^{-1,\alpha}(\partial\Omega)$  (see Theorems \ref{thm:v-1achar}, \ref{thm:v-1an}).
 
 In Section \ref{sec:trdlwd}, we consider 
 the duality   pairing in the sense of Wendland \cite{We67}, \cite{We70} of
 \[
 \left(V^{-1,\alpha}(\partial\Omega),C^{1,\alpha} (\partial\Omega)\right)
 \]
and we show that the transpose of the boundary integral operator
that corresponds to the double layer potential is compact in $V^{-1,\alpha}(\partial\Omega)$ (see Corollary \ref{corol:thm:wt-1acp}), a fact that we exploit in Section \ref{sec:eineu} to solve the interior and exterior Neumann problem for those  data in $V^{-1,\alpha}(\partial\Omega)$ on the boundary that satisfy certain compatibility conditions that generalize the classical ones (see Theorems \ref{thm:nvINexist}, \ref{thm:nvENexist3}, \ref{thm:nvENexist2}).  

In the appendix at the end of the paper, we have collected some classical results on the 
 solvability of the Dirichlet problem for  H\"{o}lder continuous harmonic functions in the specific form that we need in the paper (cf.~Appendix  \ref{sec:clsohad}). For most of the classical results in potential theory,  we refer to the joint monograph \cite{DaLaMu21} with Dalla Riva and Musolino.

\section{Preliminaries and notation}\label{sec:prelnot} Unless otherwise specified,  we assume  throughout the paper that
\[
n\in {\mathbb{N}}\setminus\{0,1\}\,,
\]
where ${\mathbb{N}}$ denotes the set of natural numbers including $0$. Let $M_n({\mathbb{R}})$ denote the set of $n\times n$ matrices with real entries.   $\delta_{l,j}$ denotes the Kronecker\index{Kronecker symbol}  symbol. Namely,  $\delta_{l,j}=1$ if $l=j$, $\delta_{l,j}=0$ if $l\neq j$, with $l,j\in {\mathbb{N}}$. $|A|$ denotes the operator norm of a matrix $A$, 
       $A^{t}$ denotes the transpose matrix of $A$.	 
 Let $\Omega$ be an open subset of ${\mathbb{R}}^n$. $C^{1}(\Omega)$ denotes the set of continuously differentiable functions from $\Omega$ to ${\mathbb{R}}$. 
 Let $s\in {\mathbb{N}}\setminus\{0\}$, $f\in \left(C^{1}(\Omega)\right)^{s} $. Then   $Df$ denotes the Jacobian matrix of $f$. For the (classical) definition of open set 
   of class $C^{m}$ or of class $C^{m,\alpha}$
  and of the H\"{o}lder and Schauder spaces $C^{m,\alpha}(\overline{\Omega})$
  on the closure $\overline{\Omega}$ of  an open set $\Omega$ and 
  of the H\"{o}lder and Schauder spaces
   $C^{m,\alpha}(\partial\Omega)$ 
on the boundary $\partial\Omega$ of an open set $\Omega$ for some $m\in{\mathbb{N}}$, $\alpha\in ]0,1]$, we refer for example to
    Dalla Riva, the author and Musolino  \cite[\S 2.3, \S 2.6, \S 2.7, \S 2.11, \S 2.13,   \S 2.20]{DaLaMu21}.    $C^{m}_b(\overline{\Omega})$ denotes the space of $m$-times continuously differentiable functions from $\Omega$ to ${\mathbb{R}}$ such that all 
the partial derivatives up to order $m$ have a bounded continuous extension to    $\overline{\Omega}$ and we set
\[
\|f\|_{   C^{m}_{b}(
\overline{\Omega} )   }\equiv
\sum_{|\beta|\leq m}\, \sup_{x\in \overline{\Omega}}|D^{\beta}f(x)|
\qquad\forall f\in C^{m}_{b}(
\overline{\Omega} )\,.
\]
$C^{m,\alpha}_b(\overline{\Omega})$ denotes the space of functions of $C^{m}_{b}(
\overline{\Omega}) $  such that the  partial derivatives of order $m$ are $\alpha$-H\"{o}lder continuous in $\Omega$. Then we equip $C^{m,\alpha}_{b}(\overline{\Omega})$ with the norm
\[
\|f\|_{  C^{m,\alpha}_{b}(\overline{\Omega})  }\equiv 
\|f\|_{  C^{m }_{b}(\overline{\Omega})  }
+\sum_{|\eta|=m}|D^{\eta}f|_{\alpha}\qquad\forall f\in C^{m,\alpha}_{b}(\overline{\Omega})\,,
\]
where $|D^{\eta}f|_{\alpha}$ denotes the $\alpha$-H\"{o}lder constant of the partial derivative $D^{\eta}f$ of order $\eta$ of $f$ in $\Omega$. 
If $\Omega$ is bounded, we obviously have $C^{m }_{b}(\overline{\Omega})=C^{m } (\overline{\Omega})$ and $C^{m,\alpha}_{b}(\overline{\Omega})=C^{m,\alpha} (\overline{\Omega})$. 
Then $C^{m,\alpha}_{{\mathrm{loc}}}(\overline{\Omega }) $\index{$C^{m,\alpha}_{{\mathrm{loc}}}(\overline{\Omega }) $} denotes 
the space  of those functions $f\in C^{m}(\overline{\Omega} ) $ such that $f_{|\overline{ \Omega\cap{\mathbb{B}}_{n}(0,\rho) }} $ belongs to $
C^{m,\alpha}(   \overline{ \Omega\cap{\mathbb{B}}_{n}(0,\rho) }   )$ for all $\rho\in]0,+\infty[$.
The space of real valued functions of class $C^\infty$ with compact support in an open set $\Omega$ of ${\mathbb{R}}^n$ is denoted ${\mathcal{D}}(\Omega)$. Then its dual ${\mathcal{D}}'(\Omega)$ is known to be the space of distributions in $\Omega$. 

If  $\Omega$ is a bounded open subset of class $C^1$ of ${\mathbb{R}}^n$, then $\Omega$ is known to have a finite number  $\varkappa^+$ connected components and the exterior  $\Omega^-$ of $\Omega$ is known to have a finite number   $\varkappa^-+1$
connected components. Then, the (bounded) connected components of $\Omega$  are denoted by $\Omega_1$, \dots, $\Omega_{\varkappa^+}$, the unbounded connected component of $\Omega^-$ is denoted by $(\Omega^-)_0$, and the bounded connected components of $\Omega^{-}$ are denoted by $(\Omega^-)_1$, \dots, $(\Omega^-)_{\varkappa^-}$. We denote by $\nu_\Omega$ or simply by $\nu$ the outward unit normal of $\Omega$ on $\partial\Omega$. Then $\nu_{\Omega^-}=-\nu_\Omega$ is the outward unit normal of $\Omega^-$ on $\partial\Omega=\partial\Omega^-$ {\color{black}(cf.~\textit{e.g.}, \cite[Lem.~2.38]{DaLaMu21}).		}

Now let $\alpha\in]0,1]$, $m\in {\mathbb{N}}$.  If $\Omega$ is a  bounded open subset of ${\mathbb{R}}^{n}$ of class $C^{\max\{m,1\},\alpha}$, then we find convenient to consider the dual $(C^{m,\alpha}(\partial\Omega))'$  of $C^{m,\alpha}(\partial\Omega)$ with its usual (normable) topology and the corresponding duality pairing $<\cdot,\cdot>$ and we say that the elements of $(C^{m,\alpha}(\partial\Omega))'$ are distributions in 
 $\partial\Omega$. Since $C^{m,\alpha}(\partial\Omega)$ is easily seen to be dense in $C^{m}(\partial\Omega)$, the transpose mapping of the canonical injection of $C^{m,\alpha}(\partial\Omega)$ into $C^{m}(\partial\Omega)$ is a continuous injective operator from  $(C^m(\partial\Omega))'$ into  $(C^{m,\alpha}(\partial\Omega))'$.      
   
Also, if $X$ is a vector subspace of the space $L^{1}(\partial\Omega)$ of Lebesgue integrable functions on $\partial\Omega$, we find convenient to set 
\begin{equation}
\label{eq:x0}
X_{0}\equiv
\left\{
f\in X:\,\int_{\partial\Omega}f\,d\sigma=0 
\right\}\,.
\end{equation}
Similarly, if   $X$ is a vector subspace of  $\left(C^{m,\alpha}(\partial\Omega)\right)'$, we find convenient to set 
\begin{equation}
\label{eq:dx0}
X_{0}\equiv
\left\{
f\in \left(C^{m,\alpha}(\partial\Omega)\right)':\,<f,1>=0 
\right\}\,.
\end{equation}
 {\color{black}Finally, we retain the standard notation for the Lebesgue spaces $L^p$ for $p\in [1,+\infty]$ (cf.~\textit{e.g.}, Folland \cite[Chapt.~6]{Fo95}, \cite[\S 2.1]{DaLaMu21}).	}
 
\section{The distributional  harmonic single layer potential}
  \label{sec:dsila}
  
  Since we are going to exploit the layer potential theoretic method, 
we  introduce the fundamental solution $S_{n}$ of the Laplace operator. Namely, we set
\[
S_{n}(\xi)\equiv
\left\{
\begin{array}{lll}
\frac{1}{s_{n}}\ln  |\xi| \qquad &   \forall \xi\in 
{\mathbb{R}}^{n}\setminus\{0\},\quad & {\mathrm{if}}\ n=2\,,
\\
\frac{1}{(2-n)s_{n}}|\xi|^{2-n}\qquad &   \forall \xi\in 
{\mathbb{R}}^{n}\setminus\{0\},\quad & {\mathrm{if}}\ n>2\,,
\end{array}
\right.
\]
where $s_{n}$ denotes the $(n-1)$ dimensional measure of 
$\partial{\mathbb{B}}_{n}(0,1)$. If $n\geq 2$, then there exists $\varsigma\in]0,+\infty[$ such that
\begin{equation}\label{thm:slayh2a}
\sup_{\xi\in {\mathbb{R}}^{n}\setminus\{0\} }|\xi|^{|\eta|+(n-2)}|D^{\eta}S_{n}(\xi)|\leq \varsigma^{|\eta|}|\eta|!\qquad\forall\eta\in{\mathbb{N}}^{n}\setminus\{0\}\,,
\end{equation}
 (cf. reference \cite[Lem.~A.6]{LaMu18} of the author and Musolino). Let $\alpha\in]0,1]$, $m\in {\mathbb{N}}$.  If $\Omega$ is a  bounded open subset of ${\mathbb{R}}^{n}$ of class $C^{\max\{m,1\},\alpha}$, then we can consider the restriction map $r_{|\partial\Omega}$ from ${\mathcal{D}}({\mathbb{R}}^n)$ to $C^{m,\alpha}(\partial\Omega)$. Then the transpose map $r_{|\partial\Omega}^t$ is linear and continuous from $(C^{m,\alpha}(\partial\Omega))'$ to ${\mathcal{D}}'({\mathbb{R}}^n)$. Moreover, if $\mu\in (C^{m,\alpha}(\partial\Omega))'$, then $r_{|\partial\Omega}^t\mu$ has compact support. Hence, it makes sense to consider the convolution of 
  $r_{|\partial\Omega}^t\mu$ with any distribution in ${\mathbb{R}}^n$ and thus in particular with the fundamental solution of the Laplace operator. Thus we are now ready to introduce the following known definition.
\begin{definition}\label{defn:dslsn}
 Let $\alpha\in]0,1]$, $m\in {\mathbb{N}}$. 
 Let   $\Omega$ be a bounded open subset of ${\mathbb{R}}^{n}$ of class $C^{\max\{m,1\},\alpha}$. If $\mu\in (C^{m,\alpha}(\partial\Omega))'$, then the (distributional) single layer potential relative to $S_{n }$ and $\mu$ is the distribution
\[
v_\Omega[\mu]=(r_{|\partial\Omega}^t\mu)\ast S_{n }\in {\mathcal{D}}'({\mathbb{R}}^n)\,.
\]
\end{definition}
In general, $(r_{|\partial\Omega}^t\mu)\ast S_{n}$ is not a function, \textit{i.e.} $(r_{|\partial\Omega}^t\mu)\ast S_{n}$ is not a distribution that is associated to a locally integrable function in ${\mathbb{R}}^n$. However, this is the case if for example $\mu$ is associated to an integrable function. So if $\mu\in L^1(\partial\Omega)$, then it is known that the (distributional) single layer potential relative to $S_{n}$ and $\mu$ is associated to the function
\begin{equation}\label{prop:dsslh1}
\int_{\partial\Omega}S_{n}(x-y)\mu(y)\,d\sigma_y\qquad{\mathrm{a.a.}}\ x\in {\mathbb{R}}^n\,,
\end{equation}
that is locally integrable in ${\mathbb{R}}^n$ and that with some abuse of notation we still denote by the symbol $v_\Omega[\mu]$.

If $\mu\in(C^{m,\alpha}(\partial\Omega))'$, then it is known that the restriction of the distributional single layer potential $(r_{|\partial\Omega}^t\mu)\ast S_{n}$ to ${\mathbb{R}}^n\setminus \partial\Omega$ is associated to the function $\theta$ defined by 
\begin{equation}\label{prop:dslfun1}
\theta(x)\equiv<(r_{|\partial\Omega}^t\mu)(y),S_{n}(x-y)>\qquad\forall x\in {\mathbb{R}}^n\setminus \partial\Omega\,.
\end{equation}
Namely,
\begin{equation}\label{prop:dslfun2}
<(r_{|\partial\Omega}^t\mu)\ast S_{n},\varphi>=\int_{{\mathbb{R}}^n\setminus \partial\Omega}
 <(r_{|\partial\Omega}^t\mu)(y),S_{n}(x-y)>\varphi(x)\,dx
 \quad\forall \varphi\in {\mathcal{D}}({\mathbb{R}}^n\setminus \partial\Omega)\,.
 \end{equation}
Then known properties of the convolution imply that
 \begin{equation}\label{prop:dslfun3}
 \Delta\left((r_{|\partial\Omega}^t\mu)\ast S_{n}\right)
 =(r_{|\partial\Omega}^t\mu)\ast (\Delta S_n)=
 (r_{|\partial\Omega}^t\mu)\ast\delta_0=(r_{|\partial\Omega}^t\mu) \quad\text{in}\ {\mathcal{D}}'({\mathbb{R}}^n)\,,
 \end{equation} 
 where $\delta_0$ is the Dirac measure with mass at $0$. 
 Since $(r_{|\partial\Omega}^t\mu)$ vanishes in ${\mathbb{R}}^n\setminus \partial\Omega$, the Weyl lemma implies that the function $\theta$ of (\ref{prop:dslfun1})  that represents the  restriction of $(r_{|\partial\Omega}^t\mu)\ast S_{n}$ to ${\mathbb{R}}^n\setminus \partial\Omega$
  is real analytic. Under the assumptions of Definition \ref{defn:dslsn}, we  set
\begin{eqnarray}\label{prop:dslfun6}
v_\Omega^+[\mu](x)&\equiv&<(r_{|\partial\Omega}^t\mu)(y),S_{n}(x-y)>
\qquad\forall x\in\Omega\,,
\\ \nonumber
v_\Omega^-[\mu](x)&\equiv&<(r_{|\partial\Omega}^t\mu)(y),S_{n}(x-y)>
\qquad\forall x\in\Omega^-\,.
\end{eqnarray}
As we have pointed out above,   $v_\Omega^\pm[\mu]$ are real analytic in  $ \Omega^\pm$. Next we introduce the following statement that generalizes the known condition for classical harmonic single layer potentials to be harmonic at infinity to the case of the distributional single layer potential.
\begin{proposition}\label{prop:sldinfty}
Let $\alpha\in]0,1]$, $m\in {\mathbb{N}}$.  Let $\Omega$ be a bounded open subset of ${\mathbb{R}}^{n}$ of class $C^{\max\{m,1\},\alpha}$. 
 Let $\mu\in \left(C^{m,\alpha}(\partial\Omega)\right)'$. 
 \begin{enumerate}
\item[(i)]  $v_\Omega^+[\mu]$ is harmonic in $\Omega$ and $v_\Omega^-[\mu]$ is harmonic in $\Omega^-$.
\item[(ii)] If $n\geq 3$, then $v_\Omega^-[\mu]$ is harmonic at infinity.
In particular, $\lim_{\xi\to\infty}v_\Omega^-[\mu](\xi)$  equals $0$.

 \item[(iii)]  If $n= 2$, then $v_\Omega^-[\mu]$ is harmonic at infinity if and only if $<\mu,1>=0$. If such a condition holds, then $\lim_{\xi\to\infty}v_\Omega^-[\mu](\xi)$ equals $0$.
\end{enumerate}
\end{proposition}
{\bf Proof.} (i) Since $r_{|\partial\Omega}^t\mu$ vanishes on ${\mathbb{R}}^n\setminus\partial\Omega$,  equality (\ref{prop:dslfun3}),  implies that $(r_{|\partial\Omega}^t\mu)\ast S_{ n }$ is harmonic in ${\mathbb{R}}^n\setminus\partial\Omega$ and 
statement (i) holds true.

 (ii) We first note that 
 \begin{eqnarray}\label{prop:sldinfty1}
\lefteqn{
|v_\Omega^-[\mu](x)|=|<\mu(y),S_n(x-\cdot>|
}
\\ \nonumber
&&\qquad\qquad\qquad
\leq\|\mu\|_{\left(C^{m,\alpha}(\partial\Omega)\right)'}
\left\|S_n(x-\cdot)\right\|_{C^{m,\alpha}(\partial\Omega)} \qquad\forall x\in\Omega^-\,.
\end{eqnarray}
Now let $r_0\in ]0,+\infty[$ be such $\overline{\Omega}\subseteq {\mathbb{B}}_n(0,r_0)$. By the definition of $S_n$ and by the inequalities (\ref{thm:slayh2a}), there exists $\varsigma\in]0,+\infty[$ such that
\begin{equation}\label{prop:sldinfty2}
|D^{\eta}S_{n}(x-y)|\leq \varsigma^{|\eta|}|\eta|!|x-y|^{-|\eta|-(n-2)}\qquad \forall x\in \Omega^-\,,\ y\in\overline{{\mathbb{B}}_n(0,r_0)}\setminus\Omega\,,
\end{equation}
for all $\eta\in{\mathbb{N}}^{n}\setminus\{0\}$
and accordingly, 
 we have $\lim_{x\to\infty}\left\|S_n(x-\cdot)\right\|_{C^{m{\color{black}+1	}
 }(\overline{{\mathbb{B}}_n(0,r_0)}\setminus\Omega)}=0$.
 Since  the restriction operator from $C^{m
 {\color{black}+1	}
 }(\overline{{\mathbb{B}}_n(0,r_0)}\setminus\Omega)$ to $C^{m,\alpha}(\partial\Omega)$ is continuous, we have
\[
\lim_{x\to\infty}\left\|S_n(x-\cdot)\right\|_{C^{m,\alpha}(\partial\Omega)}=0
\]
 and accordingly the above inequality (\ref{prop:sldinfty1}) implies that statement (ii) holds true.

(iii) Let $x_0\in \Omega$. Then we have
\begin{eqnarray}\label{prop:sldinfty3}
\lefteqn{
v_\Omega^-[\mu](x)=<\mu(y),S_2(x-y)>
}
\\ \nonumber
&&\qquad
-<\mu(y),1>S_2(x-x_0)+<\mu(y),1>S_2(x-x_0)
\qquad\forall x\in\Omega^-
\end{eqnarray}
 and
\begin{eqnarray}\label{prop:sldinfty4}
\lefteqn{
|<\mu(y),S_2(x-y)>-<\mu(y),1>S_2(x-x_0)|
}
\\ \nonumber
&&\qquad
=|<\mu(y), S_2(x-y)-S_2(x-x_0)>|
\\ \nonumber
&&\qquad
\leq\|\mu\|_{\left(C^{m,\alpha}(\partial\Omega)\right)'}
\left\|S_2(x-\cdot)-S_2(x-x_0)\right\|_{C^{m,\alpha}(\partial\Omega)} \qquad\forall x\in\Omega^-\,.
\end{eqnarray}
Since 
\[
S_2(x-y)-S_2(x-x_0)=
\frac{1}{2\pi}\log \left(
1+\left(\frac{|x-y|}{|x-x_0|}-1\right)
\right)     
\]
for all $x\in\Omega^-$, $y\in\overline{{\mathbb{B}}_n(0,r_0)}\setminus\Omega$ and
\[
\left|\frac{|x-y|}{|x-x_0|}-1\right|
=\left|\frac{|x-y|-|x-x_0|}{|x-x_0|}\right|
\leq \frac{|y-x_0|}{|x-x_0|}\qquad\forall x\in\Omega^-,   y\in\overline{{\mathbb{B}}_n(0,r_0)}\setminus\Omega\,,
\]
inequalities (\ref{prop:sldinfty2}) imply that $\lim_{x\to\infty}\left\|S_2(x-\cdot)-S_2(x-x_0)\right\|_{C^{m
{\color{black}+1	}
}(\overline{{\mathbb{B}}_n(0,r_0)}\setminus\Omega)}=0$. Since  the restriction operator from $C^{m
{\color{black}+1	}
}(\overline{{\mathbb{B}}_n(0,r_0)}\setminus\Omega)$ to $C^{m,\alpha}(\partial\Omega)$ is continuous,  we have $\lim_{x\to\infty}\left\|S_2(x-\cdot)-S_2(x-x_0)\right\|_{C^{m,\alpha}(\partial\Omega)}=0$ and accordingly inequality (\ref{prop:sldinfty4}) implies that
the harmonic function 
\[
<\mu(y),S_2(x-y)>-<\mu(y),1>S_2(x-x_0)
\]
 of the variable $x\in\Omega^-$ is   harmonic at infinity. Hence, equality (\ref{prop:sldinfty3}) implies that the function 
$v_\Omega^-[\mu]$ is harmonic at infinity if and only if the harmonic  function $<\mu(y),1>S_2(x-x_0)$ is harmonic at infinity in the variable $x\in\Omega^-$. 

Since $<\mu(y),1>S_2(x-x_0)$ is harmonic at infinity in the variable $x\in\Omega^-$ if and only if 
$<\mu(y),1>=0$, the proof of (iii) is complete.\hfill  $\Box$ 

\vspace{\baselineskip}

Then we introduce the following definition. 
\begin{definition}\label{defn:Vpm}
 Let $\alpha\in]0,1]$, $m\in {\mathbb{N}}$.  Let $\Omega$ be a bounded open subset of ${\mathbb{R}}^{n}$ of class $C^{\max\{m,1\},\alpha}$.  If $\mu\in (C^{m,\alpha}(\partial\Omega))'$ and if
$v_\Omega^+[\mu]$ admits a continuous extension to $\overline{\Omega}$ (that we denote by the same symbol) and $v_\Omega^-[\mu]$ admits a continuous extension to $\overline{\Omega^-}$ (that we denote by the same symbol), and if
\begin{equation}\label{defn:Vpm0}
v_\Omega^+[\mu](x)=v_\Omega^-[\mu](x)\qquad\forall x\in\partial\Omega\,,
\end{equation}
then we set
\begin{equation}\label{defn:Vpm1}
V_{\Omega}[\mu](x)\equiv v_\Omega^+[\mu](x)=v_\Omega^-[\mu](x)\qquad\forall x\in\partial\Omega\,.
\end{equation}
\end{definition}
We say that condition (\ref{defn:Vpm0}) is the no jump condition on the boundary for the single layer potential.
We collect in the following statement  some known properties of the single layer potential which we exploit in the sequel and that we formulate as in \cite[Thm.~7.1]{DoLa17} (cf. Miranda~\cite{Mi65}, Wiegner~\cite{Wi93}, Dalla Riva \cite{Da13}, Dalla Riva, Morais and Musolino \cite{DaMoMu13} and references therein.) 
\begin{theorem}
\label{slay}
Let  $\alpha\in]0,1[$, {\color{black}$m\in {\mathbb{N}}\setminus\{0\}$. 	}
Let $\Omega$ be a bounded open subset of ${\mathbb{R}}^{n}$ of class $C^{m,\alpha}$. Then the following statements hold. 
\begin{enumerate}
\item[(i)] If $\mu\in C^{m-1,\alpha}(\partial\Omega)$, then the function 
$v^{+}_\Omega[\mu]$ belongs to $C^{m,\alpha}(\overline{\Omega})$ and  the function 
$v^{-}_\Omega[\mu]$ belongs to $C^{m,\alpha}_{{\mathrm{loc}}}(\overline{\Omega^{-})}$.
 Moreover the map which takes $\mu$ to the function  $v^{+}_\Omega[\mu]$ is  continuous from   $C^{m-1,\alpha}(\partial\Omega)$ to $C^{m,\alpha}(\overline{\Omega})$ and the map from  the space  $C^{m-1,\alpha}(\partial\Omega)$ to $C^{m,\alpha}(\overline{{\mathbb{B}}_{n}(0,r)}\setminus \Omega) $ which takes $\mu$ to $v^{-}_\Omega[\mu]_{|\overline{{\mathbb{B}}_{n}(0,r)}\setminus \Omega}$ is  continuous for all 
$r\in]0,+\infty[$ such that $\overline{\Omega}\subseteq {\mathbb{B}}_{n}(0,r)$.
\item[(ii)] Let $l\in\{1,\dots,n\}$. If $\mu\in C^{0,\alpha}(\partial\Omega)$, then 
we have the following jump relation
\begin{eqnarray*}
\lefteqn{
\frac{\partial}{\partial x_{l}}v^{\pm}_\Omega[\mu](x)
}
\\
&&
=
\mp\frac{\nu_{l}(x)}{
2 
}\mu (x)
+
\int_{\partial\Omega} \partial_{x_{l}}S_{n}(x-y)\mu(y)\,d\sigma_{y}\qquad\forall x\in \partial\Omega
\,,
\end{eqnarray*}
where the integral in the right hand side exists in the sense of the principal value. 
If $\mu\in C^{0,\alpha}(\partial\Omega)$, then 
we have the following jump relation
\[
\frac{\partial}{\partial \nu_\Omega}v^{\pm}_\Omega[\mu](x)
=\mp\frac{1}{2}\mu(x)+\int_{\partial\Omega}\frac{\partial}{\partial \nu_{\Omega,x}}
S_{n}(x-y)\mu(y)\,d\sigma_{y}\qquad\forall x\in \partial\Omega
\,,
\]
where the integral in the right hand side exists in the sense of Lebesgue. Moreover, the no jump condition (\ref{defn:Vpm0}) holds true.
 \end{enumerate}
\end{theorem}
Next we introduce the following two lemmas that state a variant of two facts that have been shown in the proof of Corollary~4.26 of  reference
\cite{DaLaMu21} of Dalla Riva, the author and Musolino for $m=1$.
\begin{lemma}\label{lem:embour}
 Let $m\in{\mathbb{N}}$, $n\in{\mathbb{N}}\setminus\{0\}$,  $\alpha\in]0,1]$, $r\in]0,+\infty[$. Then  the space $C^{m+1}_b( {\mathbb{R}}^{n}\setminus {\mathbb{B}}_{n}(0,r) )$ is continuously imbedded into $C^{m,\alpha}_b( {\mathbb{R}}^{n}\setminus {\mathbb{B}}_{n}(0,r) )$. 
\end{lemma}
{\bf Proof.} By Lemma 2.48 of \cite{DaLaMu21},  $ {\mathbb{R}}^{n}\setminus \overline{{\mathbb{B}}_{n}(0,r)}$ is regular in the sense of Whitney (cf.~\cite[Defn.~2.27]{DaLaMu21}) and thus Lemma 2.48 and  Proposition 2.54 (i), (ii) of \cite{DaLaMu21} imply that the space $C^{m+1}_b( {\mathbb{R}}^{n}\setminus {\mathbb{B}}_{n}(0,r) )$ is continuously imbedded into $C^{m,\alpha}_b( {\mathbb{R}}^{n}\setminus {\mathbb{B}}_{n}(0,r) )$.\hfill  $\Box$ 

\vspace{\baselineskip}

\begin{lemma}\label{lem:cmaspest}
Let $m\in{\mathbb{N}}$, $n\in{\mathbb{N}}\setminus\{0\}$,  $\alpha\in]0,1]$
 Let $\Omega$ be a bounded open subset of ${\mathbb{R}}^n$. Let $r\in]0,+\infty[$ be such that $\overline{\Omega}\subseteq {\mathbb{B}}_{n}(0,r)$. If 
 $u\in C^{m}_b(\overline{\Omega^-})$ and 
\[
u_{|\overline{{\mathbb{B}}_{n}(0,r)}\setminus\Omega)}
 \in C^{m,\alpha}(\overline{{\mathbb{B}}_{n}(0,r)}\setminus\Omega)\,,\qquad u_{|{\mathbb{R}}^{n}\setminus {\mathbb{B}}_{n}(0,r) }\in C^{m,\alpha}_b( {\mathbb{R}}^{n}\setminus {\mathbb{B}}_{n}(0,r) )\,,
 \]
  then $u\in C^{m,\alpha}_b(\overline{\Omega^-})$ and
 \[
\|u\|_{C^{m,\alpha}_b(\overline{\Omega^-})}\leq 2\max\biggl\{\|u\|_{C^{m,\alpha}(\overline{{\mathbb{B}}_{n}(0,r)}\setminus\Omega)},\|u\|_{C^{m,\alpha}_b( {\mathbb{R}}^{n}\setminus {\mathbb{B}}_{n}(0,r) )}\biggr\}\,.
\]
\end{lemma}
{\bf Proof.} It suffices to check that the H\"{o}lder constant of a continuous function $g$ in $\overline{\Omega^-}$ is less or equal to $2$ times the maximum of  the  H\"{o}lder constants the restrictions of $g$ to $\overline{{\mathbb{B}}_{n}(0,r)}\setminus\Omega$ and to ${\mathbb{R}}^{n}\setminus {\mathbb{B}}_{n}(0,r) $ as done at the end of the   proof of
\cite[Cor.~4.26]{DaLaMu21}.\hfill  $\Box$ 

\vspace{\baselineskip}

Then  we can prove the following statement by following the argument of  
\cite[Cor.~4.26]{DaLaMu21}.
\begin{theorem}\label{thm:slayh}
Let  $m\in{\mathbb{N}}\setminus\{0\}$, $\alpha\in]0,1[$.
 Let $\Omega$ be a bounded open subset of ${\mathbb{R}}^{n}$ of class $C^{m,\alpha}$. Then the following statements hold.
 \begin{enumerate}
\item[(i)] Let $n=2$.  If 
 $\mu\in C^{m-1,\alpha}(\partial\Omega)_0$, then   the function 
$v^{-}_\Omega[\mu]$ belongs to $C^{m,\alpha}_{b}(\overline{\Omega^{-}})$.
 Moreover the map from $C^{m-1,\alpha}(\partial\Omega)_0$ to $C^{m,\alpha}_{b}(\overline{\Omega^{-}})$
 which takes $\mu$ to the function    $v^{-}_\Omega[\mu]_{|\overline{\Omega^-}}$ is  continuous.
\item[(ii)] Let $n\geq 3$.  If 
 $\mu\in C^{m-1,\alpha}(\partial\Omega)$, then   the function 
$v^{-}_\Omega[\mu]$ belongs to $C^{m,\alpha}_{b}(\overline{\Omega^{-}})$.
 Moreover the map from $C^{m-1,\alpha}(\partial\Omega)$ to $C^{m,\alpha}_{b}(\overline{\Omega^{-}})$
 which takes $\mu$ to the function    $v^{-}_\Omega[\mu]_{|\overline{\Omega^-}}$ is  continuous.
\end{enumerate}
\end{theorem}
{\bf Proof.} Let $r\in]1,+\infty[$ 
{\color{black}be	}
such that $\overline{\Omega}\subseteq {\mathbb{B}}_{n}(0,r-1)$. By Lemmas
\ref {lem:embour}, \ref{lem:cmaspest}  
 and Theorem \ref{slay}, it suffices to show that  ${v}^-_\Omega[\cdot]_{| {\mathbb{R}}^{n}\setminus {\mathbb{B}}_{n}(0,r) }$  is linear and continuous from  $C^{m-1,\alpha}(\partial\Omega)_0$ to  $C^{m+1}_b( {\mathbb{R}}^{n}\setminus {\mathbb{B}}_{n}(0,r) )$ in case $n=2$ and 
 from  $C^{m-1,\alpha}(\partial\Omega)$ to  $C^{m+1}_b( {\mathbb{R}}^{n}\setminus {\mathbb{B}}_{n}(0,r) )$ in case $n\geq 3$.  
 
 If $n=2$, there exists $r_1\in ]r,+\infty[$ such that
 \begin{equation}\label{thm:slayh1}
 |2\pi {v}^-_\Omega[\mu](x)|=
\left|
\int_{\partial\Omega}\mu(y)\log\left(
1+\left(
\frac{|x-y|}{|x|}-1
\right)
\right)\,d\sigma_{y}
\right|
 \leq\frac{3}{2}\frac{r}{|x|}\int_{\partial\Omega }|\mu(y)|
\,d\sigma_{y}
\end{equation}
for all $x\in {\mathbb{R}}^{2}\setminus{\mathbb{B}}_{2}(0,r_{1})$ 
and $\mu\in C^{0,\alpha}(\partial\Omega)_0$ (cf.  \cite[(4.68), proof of 
{\color{black} Thm.~4.26	}
]{DaLaMu21}).

If $n\geq 3$, then
 \begin{equation}\label{thm:slayh2}
 |v^{-}_{\Omega}[\mu](x)|\leq
 \frac{1}{(n-2)s_{n}}\frac{1}{||x|-r|^{n-2}}
 \int_{\partial\Omega}
|\mu(y)|\,d\sigma_{y}\qquad\forall x\in {\mathbb{R}}^{n}\setminus
\overline{{\mathbb{B}}_{n}(0,r)}\,,
 \end{equation}
for all  $\mu\in C^{0,\alpha}(\partial\Omega)$  (cf.  \cite[(4.67)]{DaLaMu21}).

 If $n\geq 2$, inequality (\ref{thm:slayh2a}) and the classical differentiability theorem for integrals depending on a parameter and  the elementary inequality
\begin{equation}\label{thm:slayh2b}
|x-y|\geq |x|-(r-1)\geq r-(r-1)=1\qquad\forall (x,y)\in ({\mathbb{R}}^{n}\setminus {\mathbb{B}}_{n}(0,r)) \times (\partial\Omega)\,,
\end{equation}
imply that
\begin{eqnarray}\label{thm:slayh3}
\lefteqn{|D^\eta{v}^-_\Omega[\mu](x)|
=\left|\int_{\partial\Omega}D^\eta S_n(x-y)\mu(y)\,d\sigma_y\right|
}
\\ \nonumber
&&\qquad\qquad\qquad\qquad\qquad 
\leq \int_{\partial\Omega} \frac{ \varsigma^{|\eta|}|\eta|!d\sigma_y}{|x-y|^{|\eta|+(n-2)}}
\|\mu\|_{C^{0,\alpha}(\partial\Omega)}
\\ \nonumber
&&\qquad\qquad\qquad\qquad\qquad 
\leq \varsigma^{|\eta|}|\eta|!
m_{n-1}(\partial\Omega)\|\mu\|_{C^{0,\alpha}(\partial\Omega)}
\ \forall x\in {\mathbb{R}}^{n}\setminus {\mathbb{B}}_{n}(0,r)
\end{eqnarray}
for all $\eta\in {\mathbb{N}}^n\setminus\{0\}$ such that $|\eta|\leq m+1$ and $\mu\in C^{0,\alpha}(\partial\Omega)$. Then by combining inequalities (\ref{thm:slayh1}), (\ref{thm:slayh3}) we deduce the validity of statement (i) and by combining inequalities (\ref{thm:slayh2}), (\ref{thm:slayh3}) we deduce the validity of statement (ii). Indeed, $C^{m-1,\alpha}(\partial\Omega)$ is continuously embedded into $C^{0,\alpha}(\partial\Omega)$.\hfill  $\Box$ 

\vspace{\baselineskip}

 \section{Classical properties of the harmonic double layer potential}
  \label{sec:hadolapo}
  
  We collect in the present section some variants of classical known facts that we need in the sequel. For the convenience of the reader, we also present some of the proofs by following the scheme of \cite{DaLaMu21}.
 \begin{definition}\label{bvps.wpsi} 
 Let   $\alpha\in]0,1]$. 
 Let $\Omega$ be a bounded open   subset of ${\mathbb{R}}^{n}$ of class $C^{1,\alpha}$.  If   $\psi\in   C^{0}  (\partial\Omega)$, then we  denote by ${w}_\Omega[\psi]$ the double layer potential with moment (or density) $\psi$, i.e., the function from ${\mathbb{R}}^{n}$ to ${\mathbb{R}}$ defined by\index{${w}_\Omega[\psi]$} 
 \[
{w}_\Omega[\psi](x)\equiv \int_{\partial\Omega}\psi(y)\frac{\partial}{\partial \nu_{\Omega_y} }
\left(S_n(x-y)\right)\,d\sigma_y\qquad\forall x\in\mathbb{R}^n\,,
\]
where
\[
\frac{\partial}{\partial \nu_{\Omega_y} }
\left(S_n(x-y)\right)\equiv
 - \nu_{\Omega}(y)\cdot D S_n(x-y) \qquad\forall (x,y)\in\mathbb{R}^n\times\partial\Omega\,, x\neq y\,.
\]
\end{definition} 
Next we state the following classical result on the double layer potential.
\begin{theorem}\label{thm:wpsi+-}
Let   $m\in{\mathbb{N}}$, $\alpha\in]0,1[$.
 Let $\Omega$ be a bounded open subset of ${\mathbb{R}}^{n}$ of class $C^{\max\{m,1\},\alpha}$. 
 \begin{enumerate}
 \item[(i)]  If $\psi\in C^{m,\alpha}(\partial\Omega)$, then the  restriction ${w}_\Omega[\psi]_{|\Omega}$ extends to a function ${w}^+_\Omega[\psi]$\index{${w}^+_\Omega[\psi]$} of class $C^{m,\alpha}(\overline\Omega)$.
 Moreover, the map  from $C^{m,\alpha}(\partial\Omega)$ to 
 $C^{m,\alpha}(\overline\Omega)$ that takes $\psi$ to ${w}^+_\Omega[\psi]$ is linear and continuous.
 \item[(ii)] If $\psi\in C^{m,\alpha}(\partial\Omega)$, then the restriction ${w}_\Omega[\psi]_{| \Omega^-}$ extends to a function ${w}^-_\Omega[\psi]$  of $C^{m,\alpha}_{b}(\overline{\Omega^-})$.  Moreover, then the map
 from $C^{m,\alpha}(\partial\Omega)$ to $C^{m,\alpha}_b(\overline{\Omega^-})$ that takes $\psi$ to ${w}^-_\Omega[\psi]_{| \overline{\Omega^-}  }$ is linear and continuous. 
\item[(iii)] If $\psi\in C^{0,\alpha}(\partial\Omega)$, then  the following jump formula holds true
\[
w^{\pm}_\Omega[\psi](x)
=\pm\frac{1}{2}\psi(x)+w_\Omega[\psi](x)
\qquad\forall x\in\partial\Omega\,.
\]
\end{enumerate}
\end{theorem}
{\bf Proof.}  (i) If $m\geq 1$, we refer  to   Miranda~\cite{Mi65}, Wiegner~\cite{Wi93}, Dalla Riva \cite{Da13}, Dalla Riva, Morais and Musolino \cite{DaMoMu13} and references therein. See also reference \cite[Thm.~7.3]{DoLa17} with Dondi.
 If instead $m=0$, we note that
 the definition of double layer potential implies that
\begin{equation}\label{thm:wpsi+-2}
w_{\Omega}[\psi](x)  =-\sum_{j=1}^{n}\frac{\partial}{\partial x_{j}}v_{\Omega}[\psi  \cdot  \nu_{j}] (x)\qquad\forall x\in  {\mathbb{R}}^n\setminus\partial\Omega\,,
\end{equation}
for all $\psi\in C^{0,\alpha}(\partial\Omega)$. 
Since the components of $\nu$ belong to $C^{0,\alpha}(\partial\Omega)$, the continuity of the pointwise product 
from 
 $C^{0,\alpha}(\partial\Omega)\times C^{0,\alpha}(\partial\Omega)$ to $C^{0,\alpha}(\partial\Omega)$
 and Theorem \ref{slay} on the regularity of the single layer potential imply the validity of statement (i). 

(ii) Let $r\in]1,+\infty[$ be such that $\overline{\Omega}\subseteq{\mathbb{B}}_{n}(0,r-1) $. By Lemmas
\ref {lem:embour}, \ref{lem:cmaspest} 
   it suffices to show that the map  ${w}^-_\Omega[\cdot]_{| \overline{{\mathbb{B}}_{n}(0,r)}\setminus\Omega}$  is linear and continuous from  $C^{m,\alpha}(\partial\Omega)$ to  $C^{m,\alpha}( \overline{{\mathbb{B}}_{n}(0,r) }\setminus\Omega)$
   and that  the map ${w}^-_\Omega[\cdot]_{| {\mathbb{R}}^{n}\setminus {\mathbb{B}}_{n}(0,r) }$  is linear and continuous from  $C^{m,\alpha}(\partial\Omega)$ to  $C^{m+1}_b( {\mathbb{R}}^{n}\setminus {\mathbb{B}}_{n}(0,r) )$.  We first consider ${w}^-_\Omega[\cdot]_{| \overline{{\mathbb{B}}_{n}(0,r)}\setminus\Omega}$. For case $m>1$, we refer to     Miranda~\cite{Mi65}, Wiegner~\cite{Wi93}, Dalla Riva \cite{Da13}, Dalla Riva, Morais and Musolino \cite{DaMoMu13} and references therein. See also reference \cite[Thm.~7.3]{DoLa17} with Dondi.  If instead $m=0$, we note that the membership of 
  the components of $\nu$  to $C^{0,\alpha}(\partial\Omega)$,  the continuity of the pointwise product 
from 
 $C^{0,\alpha}(\partial\Omega)\times C^{0,\alpha}(\partial\Omega)$ to $C^{0,\alpha}(\partial\Omega)$, Theorem \ref{slay} on the regularity of the single layer potential and equality equality (\ref{thm:wpsi+-2}) imply that ${w}^-_\Omega[\cdot]_{| \overline{{\mathbb{B}}_{n}(0,r)}\setminus\Omega}$  is linear and continuous from  $C^{0,\alpha}(\partial\Omega)$ to  $C^{0,\alpha}( \overline{{\mathbb{B}}_{n}(0,r)}\setminus\Omega)$. 
 
 Next we consider ${w}^-_\Omega[\cdot]_{| {\mathbb{R}}^{n}\setminus {\mathbb{B}}_{n}(0,r) }$ and we note that 
 {\color{black}equality (\ref{thm:wpsi+-2}),	}
 the classical differentiability theorem for integrals depending on a parameter and  the  inequalities (\ref{thm:slayh2a}), 
 imply that
\begin{eqnarray}\label{thm:wpsi+-3}
\lefteqn{|D^\eta{w}^-_\Omega[\mu](x)|
=\left|-\int_{\partial\Omega}\sum_{j=1}^nD^\eta_x \frac{\partial}{\partial x_j}S_n(x-y) \mu(y)\nu_j(y)\,d\sigma_y\right|
}
\\ \nonumber
&&\qquad\qquad
\leq  \int_{\partial\Omega} n\frac{ \varsigma^{|\eta|+1|}(|\eta|+1)!d\sigma_y}{|x-y|^{|\eta|+1+(n-2)}}
\|\mu\|_{C^{0,\alpha}(\partial\Omega)}
\\ \nonumber
&&\qquad\qquad  
\leq n\varsigma^{|\eta|+1|}(|\eta|+1)!
m_{n-1}(\partial\Omega)\|\mu\|_{C^{0,\alpha}(\partial\Omega)}
\ \forall x\in {\mathbb{R}}^{n}\setminus {\mathbb{B}}_{n}(0,r)
\end{eqnarray}
for all $\eta\in {\mathbb{N}}^n$ such that $|\eta|\leq m+1$ and $\mu\in C^{0,\alpha}(\partial\Omega)$. Then inequality (\ref{thm:wpsi+-3}) implies continuity of ${w}^-_\Omega[\cdot]_{| {\mathbb{R}}^{n}\setminus {\mathbb{B}}_{n}(0,r) }$. Indeed, $C^{m,\alpha}(\partial\Omega)$ is continuously embedded into $C^{0,\alpha}(\partial\Omega)$.
Hence, statement (ii) holds true.  For the jump formula of statement (iii), we refer for example to \cite[Thm.~4.30]{DaLaMu21}\hfill  $\Box$ 

\vspace{\baselineskip}

 Then the following result of  Schauder \cite[Hilfsatz VII, p.~112]{Sc31} is well-known to hold (cf.~\textit{e.g.}, \cite[Thm.~4.33]{DaLaMu21}).
 \begin{theorem}
\label{thm:schauin}
Let $\alpha\in ]0,1[$. Let $\Omega$ be a bounded open  subset of ${\mathbb{R}}^{n}$ of class $C^{1,\alpha}$. Then the following statements hold. 
\begin{enumerate}
\item[(i)] Let $\psi\in L^{\infty}(\partial\Omega)$. Then the function $W_{\Omega}[\psi]$ from $\partial\Omega$ to 
${\mathbb{R}}$ defined by \index{$W_{\Omega}[\psi]$}
\begin{equation}\label{bvps.W}
 W_{\Omega}[\psi](x) 
 \equiv 
-\int_{\partial\Omega}\psi(y)\;\nu_{\Omega}(y)\cdot\nabla S_n(x-y)\, d\sigma_y
\qquad\forall x\in \partial\Omega
\end{equation}
belongs to $C^{0,\alpha}(\partial\Omega)$.
\item[(ii)] The operator $W_{\Omega}$ from  $L^{\infty}(\partial\Omega)$ to $C^{0,\alpha}(\partial\Omega)$ that takes $\psi$ to  $W_{\Omega}[\psi]$ is linear and continuous. 
\end{enumerate}
\end{theorem}
Since $C^{0,\alpha}(\partial\Omega)$ is well known to be compactly embedded into $L^{\infty}(\partial\Omega)$, the operator $W_{\Omega}$ is compact from $C^{0,\alpha}(\partial\Omega)$ to itself. As a consequence, its transpose $W_{\Omega}^t$ is linear and continuous from $(C^{0,\alpha}(\partial\Omega))'$ to itself. Next we point out the validity of the following extension of a known classical lemma (cf.~\textit{e.g.}, \cite[Thm.~6.11]{DaLaMu21}).
\begin{lemma}\label{lem:Sphi0a}
Let   $m\in{\mathbb{N}}$, $\alpha\in]0,1[$.
 Let $\Omega$ be a bounded open subset of ${\mathbb{R}}^{n}$ of class $C^{\max\{m,1\},\alpha}$. 
  Then $W_{\Omega}[1]=\frac{1}{2}$ on $\partial\Omega$ and
 \[
 <\left(\frac{1}{2}I+W_{\Omega}^t\right)[\mu],1>=<\mu,1>\qquad\forall \mu\in (C^{m,\alpha}(\partial\Omega))'\,.
 \]
\end{lemma}
 {\bf Proof.} Equality $W_{\Omega}[1]=\frac{1}{2}$ on $\partial\Omega$ is well-known
 (cf.~\textit{e.g.}, \cite[Thm.~6.9]{DaLaMu21}). 
  If $\mu\in (C^{m,\alpha}(\partial\Omega))'$, then
 \begin{eqnarray*}
\lefteqn{
<\left(\frac{1}{2}I+W_{\Omega}^t\right)[\mu],1>=< \mu,\left(\frac{1}{2}I+W_{\Omega}^t\right)[1]>
}
\\ \nonumber
&&\qquad
=< \mu,\left(-\frac{1}{2}I+W_{\Omega}\right)[1]>+<\mu,1>=0+<\mu,1>\,.
\end{eqnarray*}
\hfill  $\Box$ 

\vspace{\baselineskip}

Next we introduce the following classical result on the Green operator for the interior Dirichlet problem. 
 \begin{theorem}\label{thm:idwp}
Let $m\in {\mathbb{N}}$, $\alpha\in ]0,1[$. Let $\Omega$ be a bounded open  subset of ${\mathbb{R}}^{n}$ of class $C^{\max\{m,1\},\alpha}$.
Then the map ${\mathcal{G}}_{d,+}$  from $C^{m,\alpha}(\partial\Omega)$ to the closed subspace 
 \begin{equation}\label{thm:idwp1}
C^{m,\alpha}_h(\overline{\Omega}) \equiv \{
u\in C^{m,\alpha}(\overline{\Omega}), u\ \text{is\ harmonic\ in}\ \Omega\}
\end{equation}
of $ C^{m,\alpha}(\overline{\Omega})$ that takes $v$ to the only solution $v^\sharp$ of the Dirichlet problem
\begin{equation}\label{defn:cinspo3}
\left\{
\begin{array}{ll}
 \Delta v^\sharp=0 & \text{in}\ \Omega\,,
 \\
v^\sharp_{|\partial\Omega} =v& \text{on}\ \partial\Omega 
\end{array}
\right.
\end{equation}
is a linear homeomorphism.
\end{theorem}
{\bf Proof.}  Since the uniform limit of harmonic functions is harmonic, 
$C^{m,\alpha}_h(\overline{\Omega})$ is a closed subspace of the Banach space $C^{m,\alpha}(\overline{\Omega})$. The restriction map from the Banach subspace 
$C^{m,\alpha}_h(\overline{\Omega})$
of $C^{m,\alpha}(\overline{\Omega})$ to $C^{m,\alpha}(\partial\Omega)$ is linear and continuous. By 
{\color{black}the classical }
existence and uniqueness Theorem \ref{thm:existenceO0} 
{\color{black}of the Appendix	}
 for the Dirichlet problem, such a map is also a linear isomorphism. Then the Open Mapping Theorem implies that its inverse is continuous. Since such inverse coincides with the map ${\mathcal{G}}_{d,+}$, the proof is complete.\hfill  $\Box$ 

\vspace{\baselineskip}

Similarly, we introduce the following classical result on the Green operator for the exterior Dirichlet problem. 
\begin{theorem}\label{thm:edwp}
 Let $m\in {\mathbb{N}}$, $\alpha\in ]0,1[$. Let $\Omega$ be a bounded open  subset of ${\mathbb{R}}^{n}$ of class $C^{\max\{m,1\},\alpha}$. The map ${\mathcal{G}}_{d,-}$ from $C^{m,\alpha}(\partial\Omega)$ to 
the closed subspace
\begin{equation}\label{thm:edwp1}
C^{m,\alpha}_{bh}(\overline{\Omega^-}) \equiv \{
u\in C^{m,\alpha}_b(\overline{\Omega^-}), u\ \text{is\ harmonic\ in}\ \Omega^- 
 \text{and\ is\ harmonic\ at}\ \infty\}
\end{equation}
of $C^{m,\alpha}_{b}(\overline{\Omega^-})$ that takes $v$ to the only solution $v^\sharp$ of   in $ C^{m,\alpha}_{b}(\overline{\Omega^-})$ of the exterior Dirichlet problem
\begin{equation}\label{defn:censpo3}
\left\{
\begin{array}{ll}
 \Delta v^\sharp=0 & \text{in}\ \Omega^-\,,
 \\
v^\sharp_{|\partial\Omega} =v& \text{on}\ \partial\Omega\,,
\\
v^\sharp\  \text{is\ harmonic\ at}\ \infty   &
\end{array}
\right.
\end{equation}
is a linear homeomorphism. Moreover, 
\begin{equation}\label{cor:existenceO-01a}
\lim_{\xi\to\infty}{\mathcal{G}}_{d,-}[g](\xi)=\left\{
\begin{array}{ll}
 0 &\quad\text{if}\ n\geq 3\,,
 \\
 \Xint{-}_{ \partial{\mathbb{B}}_{n}(0,r) }{\mathcal{G}}_{d,-}[g]\,d\sigma &\quad\text{if}\ n=2 
\end{array}
\right.
\end{equation}
for all $r\in]0,+\infty[$ such that $\overline{\Omega}\subseteq {\mathbb{B}}_{n}(0,r)$
and the map from $C^{0,\alpha}(\partial\Omega)$ to ${\mathbb{R}}$ that takes $g$ to
$\lim_{\xi\to\infty}{\mathcal{G}}_{d,-}[g]$ is continuous.
\end{theorem} 
{\bf Proof.} 
Since the uniform limit of harmonic functions is harmonic and the uniform  limit of harmonic functions at infinity is harmonic at infinity, $C^{m,\alpha}_{bh}(\overline{\Omega^-})$ is a closed subspace of the Banach space $C^{m,\alpha}_{b}(\overline{\Omega^-})$. The restriction map from the Banach subspace $C^{m,\alpha}_{bh}(\overline{\Omega^-})$ 
of $C^{m,\alpha}_b(\overline{\Omega^-})$ to $C^{m,\alpha}(\partial\Omega)$ is linear and continuous. By the {\color{black}classical	}
 existence Theorem \ref{thm:existenceO-0} {\color{black}of the Appendix,	}
 such a map is also a linear isomorphism. Then the Open Mapping Theorem implies that its inverse is continuous. Since such inverse coincides with the map ${\mathcal{G}}_{d,-}$, the proof of the first part of the statement is complete. For a proof of the  classical formula (\ref{cor:existenceO-01a}), we refer for example to \cite[Cor.~3.24]{DaLaMu21}. Then the continuity statement for the limit follows by the  classical formula (\ref{cor:existenceO-01a}), by the continuity of ${\mathcal{G}}_{d,-}[\cdot]$ and by the continuity  the restriction map from $C^{0,\alpha}_{bh}(\overline{\Omega^-})$ to the space of integrable functions in $\partial{\mathbb{B}}_n(0,r)$.\hfill  $\Box$ 

\vspace{\baselineskip}

\section{A distributional  normal derivative for harmonic H\"{o}lder continuous functions in a bounded open set}\label{sec:conoderd}
 
Let $\alpha\in]0,1{\color{black}[	}
$.  Let  $\Omega$ be a  bounded open subset of ${\mathbb{R}}^{n}$ of class $C^{1,\alpha}$.
We plan to define a  normal derivative for a function $u\in C^{0,\alpha}(\overline{\Omega})$ that is harmonic in $\Omega$.

If $u$ were to belong  to the Sobolev space $H^{1}(\Omega)$ of  functions in $L^2(\Omega)$ which have first order distributional derivatives in $L^2(\Omega)$ and $u$ were to be harmonic in $\Omega$, then one could classically  define the distributional  normal derivative $\frac{\partial u}{\partial\nu}$ of $u$ on $\partial\Omega$ to be the only element of the dual $H^{-1/2}(\partial\Omega)$ of the space $H^{1/2}(\partial\Omega)$ of traces on $\partial\Omega$ of $H^{1}(\Omega)$  that is defined by the equality
\begin{equation}\label{eq:conoder}
<\frac{\partial u}{\partial\nu},v>\equiv\int_{\Omega}DuD(Ev)\,dx\qquad\forall v\in H^{1/2}(\partial\Omega)\,,
\end{equation}
 where $E$ is any bounded extension operator from $H^{1/2}(\partial\Omega)$ to $H^{1}(\Omega)$ (cf.~\textit{e.g.}, Lions and Magenes~\cite{LiMa68}, Ne\v{c}as \cite[Chapt.~5]{Ne12},
 Nedelec and Planchard \cite[p.~109]{NePl73},  Costabel \cite{Co88}, McLean~\cite[Chapt.~4]{McL00}, Mikhailov~\cite{Mi11}, Mitrea, Mitrea and Mitrea \cite[\S 4.2]{MitMitMit22}). 
 We also note that if $u$ were to be continuously differentiable, then the first Green identity would
 imply the validity of such equality with $\int_{\partial\Omega}\frac{\partial u}{\partial\nu} v\,d\sigma$ instead of 
 $<\frac{\partial u}{\partial\nu},v>$. Thus definition (\ref{eq:conoder}) implies that $\frac{\partial u}{\partial\nu}$ is required to satisfy a generalized form of the classical first Green Identity.
 
 However functions in $C^{0,\alpha}(\overline{\Omega})$ are not necessarily in $H^{1}(\Omega)$ (for a discussion on this point we refer to  Bramati,  Dalla Riva and  Luczak~\cite{BrDaLu23}). Thus we now introduce a different notion of distributional  normal derivative $\frac{\partial u}{\partial\nu}$ that requires that $\frac{\partial u}{\partial\nu}$ satisfies a generalized form of the classical second Green Identity. To do so, we need some preliminaries. 
We find convenient to introduce the (classical) interior Steklov-Poincar\'{e} operator.
\begin{definition}\label{defn:cinspo}
 Let $\alpha\in]0,1[$.  Let  $\Omega$ be a  bounded open subset of ${\mathbb{R}}^{n}$ of class $C^{1,\alpha}$. The classical interior Steklov-Poincar\'{e} operator is defined to be the operator $S_+$ from
 \begin{equation}\label{defn:cinspo1}
C^{1,\alpha}(\partial\Omega)\quad\text{to}\quad C^{0,\alpha}(\partial\Omega)
\end{equation}
 takes $v\in C^{1,\alpha}(\partial\Omega)$ to the function 
 \begin{equation}\label{defn:cinspo2}
S_+[v](x)\equiv \frac{\partial  }{\partial\nu}{\mathcal{G}}_{d,+}[v](x)\qquad\forall x\in\partial\Omega\,.
\end{equation}
 \end{definition}
   Since   the classical normal derivative is continuous from $C^{1,\alpha}(\overline{\Omega})$ to $C^{0,\alpha}(\partial\Omega)$, the continuity of ${\mathcal{G}}_{d,+}$ implies  that $S_+[\cdot]$ is linear and continuous from 
  $C^{1,\alpha}(\partial\Omega)$ to $C^{0,\alpha}(\partial\Omega)$. We are now ready to introduce the following definition.
\begin{definition}\label{defn:conoderd}
 Let $\alpha\in]0,1[$.  Let  $\Omega$ be a  bounded open subset of ${\mathbb{R}}^{n}$ of class $C^{1,\alpha}$. If  $u\in C^{0}(\overline{\Omega})$ and $u$ is harmonic in $\Omega$, then we define the distributional  normal derivative $\frac{\partial u}{\partial\nu}$ to be the only element of the dual $(C^{1,\alpha}(\partial\Omega))'$ that satisfies the following equality
 \begin{equation}\label{defn:conoderd1}
<\frac{\partial u}{\partial\nu},v>\equiv\int_{\partial\Omega}uS_+[v]\,d\sigma
\qquad\forall v\in C^{1,\alpha}(\partial\Omega)\,.
\end{equation}
\end{definition}
Here we have introduced the Definition  \ref{defn:conoderd} for functions of $C^{0}(\overline{\Omega})$, but one could do the same also for harmonic functions in other function spaces that have a  trace operator on $\partial\Omega$ and  where we do not have information on the integrability of the first order partial derivatives of $u$ in $\Omega$.\par 

The following Lemma is well known and is an immediate consequence of the H\"{o}lder inequality.
\begin{lemma}\label{lem:caincl}
 Let $m\in {\mathbb{N}}$, $\alpha\in ]0,1[$. Let $\Omega$ be a bounded open  subset of ${\mathbb{R}}^{n}$ of class $C^{\max\{m,1\},\alpha}$.  Then the canonical inclusion  ${\mathcal{J}}$ from the Lebesgue space $L^1(\partial\Omega)$ of integrable functions in $\partial\Omega$ to $(C^{m,\alpha}(\partial\Omega))'$ that takes $\mu$ to the functional ${\mathcal{J}}[\mu]$ defined by 
 \begin{equation}\label{lem:caincl1}
<{\mathcal{J}}[\mu],v>\equiv \int_{\partial\Omega}\mu v\,d\sigma\qquad\forall v\in C^{
{\color{black}m	},
\alpha}(\partial\Omega)\,,
\end{equation}
is linear continuous and injective.
\end{lemma}
As customary, we say  that ${\mathcal{J}}[\mu]$ is the `distribution' that is canonically associated to $\mu$ and we omit the indication of the inclusion map ${\mathcal{J}}$. By Lemma \ref{lem:caincl}, the space $C^{0,\alpha}(\partial\Omega)$ is continuously embedded into $(C^{m,\alpha}(\partial\Omega))'$.  
 
 By applying  the classical second Green Identity for harmonic functions in $\Omega$ to the function $u$ and to the function ${\mathcal{G}}_{d,+}[v]$, one can show
   that if $ u\in C^{1}(\overline{\Omega})$ is harmonic in $\Omega$, then 
 \[
 <\frac{\partial u}{\partial\nu},v>=\int_{\partial\Omega}\frac{\partial u}{\partial\nu} v\,d\sigma\qquad\forall v\in C^{1,\alpha}(\partial\Omega)\,,
 \]
  where $\frac{\partial u}{\partial\nu}$ in the right hand side coincides with the classical normal derivative of $u$, \textit{i.e.}, the distributional  normal derivative $\frac{\partial u}{\partial\nu}$ of definition (\ref{defn:conoderd1}) coincides with the distribution of $(C^{1,\alpha}(\partial\Omega))'$ that is associated to the classical normal derivative of $u$. Next we point out the validity of the following statement on the continuity of the (distributional) normal derivative.\par
\begin{theorem}\label{thm:contidnu}
 Let $\alpha\in ]0,1[$. Let $\Omega$ be a bounded open  subset of ${\mathbb{R}}^{n}$ of class $C^{1,\alpha}$. Then the following statements hold.
 \begin{enumerate}
\item[(i)] The transpose operator $S_+^t$ of the map $S_+$ is linear and continuous from 
\begin{equation}\label{eq:trispo}
\left(C^{0,\alpha}(\partial\Omega)\right)'\quad\text{to}\quad \left(C^{1,\alpha}(\partial\Omega)\right)'\,.
\end{equation}
\item[(ii)] The map $\frac{\partial}{\partial\nu}$ from the closed subspace $C^{0,\alpha}_h(\overline{\Omega}) $ of $ C^{0,\alpha}(\overline{\Omega})$ to $\left(C^{1,\alpha}(\partial\Omega)\right)'$ is linear and continuous. Moreover,
\begin{equation}\label{eq:trispo1}
\frac{\partial u}{\partial\nu} =S_+^t[u_{|\partial\Omega}]
  \qquad\forall u\in C^{0,\alpha}_h(\overline{\Omega})\,.
  \end{equation}
\end{enumerate}
\end{theorem}
{\bf Proof.} The continuity of $S_+^t$ follows by the continuity of  the linear operator  $S_+$ from $C^{1,\alpha}(\partial\Omega)$ to $C^{0,\alpha}(\partial\Omega)$. Since
 \[
 <S_+^t[u_{|\partial\Omega}],v>=<u_{|\partial\Omega},S_+[v]>=\int_{\partial\Omega}uS_+[v]\,d\sigma
 =<\frac{\partial u}{\partial\nu},v>\quad\forall v\in C^{1,\alpha}(\partial\Omega)\,,
\]
for all $u\in C^{0,\alpha}_h(\overline{\Omega})$, equality (\ref{eq:trispo1}) holds true. Since the restriction map is continuous from $C^{0,\alpha}_h(\overline{\Omega})$ to $C^{0,\alpha}(\partial\Omega)$  and the canonical embedding of $C^{0,\alpha}(\partial\Omega)$ into the dual $\left(C^{
{\color{black}0	}
,\alpha}(\partial\Omega)\right)'$ is continuous, the continuity of $S_+^t$ as in statement (i) implies the validity of (ii).\hfill  $\Box$ 

\vspace{\baselineskip}  
  
    In other words, statement (ii) says that the
distributional  normal derivative of $u$ coincides with the transpose of the classical interior Steklov operator applied to the distribution that is canonically associated to the restriction of $u$ to $\partial\Omega$. Once more, one 
  could do the same also for harmonic functions in $\Omega$ in function spaces that have a  trace operator on $\partial\Omega$, even in cases where we do not have information on the integrability of the first order partial derivatives of $u$ in $\Omega$.\par 

  \section{A distributional  normal derivative for harmonic H\"{o}lder continuous functions in the  exterior of a bounded open set}\label{sec:conoederd}
  
  Let $\alpha\in]0,1
  {\color{black}[	}
  $.  Let  $\Omega$ be a  bounded open subset of ${\mathbb{R}}^{n}$ of class $C^{1,\alpha}$.
  We plan to define a normal derivative for a function $u\in C^{0,\alpha}_{{\mathrm{loc}}}(\overline{\Omega^-})$ that is both harmonic in $\Omega^-$ and harmonic at infinity and we follow the same pattern for harmonic functions in $\Omega$. To do so, we need some preliminaries.  
 We find convenient to introduce the  (classical) exterior Steklov-Poincar\'{e} operator (for the definition of the exterior Steklov-Poincar\'{e} operator in a variational context, we refer to Arendt and ter Elst \cite{Arte15}, where a different approach in a Sobolev space setting has been introduced).
\begin{definition}\label{defn:censpo}
 Let $\alpha\in]0,1[$.  Let  $\Omega$ be a  bounded open subset of ${\mathbb{R}}^{n}$ of class $C^{1,\alpha}$. The classical interior Steklov-Poincar\'{e} operator is defined to be the operator $S_-$ from
 \begin{equation}\label{defn:censpo1}
C^{1,\alpha}(\partial\Omega)\quad\text{to}\quad C^{0,\alpha}(\partial\Omega)
\end{equation}
 takes $v\in C^{1,\alpha}(\partial\Omega)$ to the function 
 \begin{equation}\label{defn:censpo2}
S_-[v](x)\equiv \frac{\partial}{\partial\nu_{\Omega^-}}{\mathcal{G}}_{d,-}[v](x)\qquad\forall x\in\partial\Omega\,,
\end{equation}
where $\nu_{\Omega^-}=-\nu_\Omega$ is the outward unit normal to $\Omega^-$ on $\partial\Omega=\partial\Omega^-$.
 \end{definition}
 Since   the classical normal derivative is continuous from $C^{1,\alpha}_b(\overline{\Omega^-})$ to $C^{0,\alpha}(\partial\Omega)$, the continuity of ${\mathcal{G}}_{d,-}$ implies  that  $S_-[\cdot]$ is linear and continuous from 
  $C^{1,\alpha}(\partial\Omega)$ to $C^{0,\alpha}(\partial\Omega)$.  We are now ready to introduce the following definition.
\begin{definition}\label{defn:conoederd}
 Let $\alpha\in]0,1[$.  Let  $\Omega$ be a  bounded open subset of ${\mathbb{R}}^{n}$ of class $C^{1,\alpha}$. If  $u\in C^{0}(\overline{\Omega^-})$ and $u$ is both harmonic in $\Omega^-$ and harmonic at infinity, then we define the distributional  normal derivative $\frac{\partial u}{\partial\nu_{\Omega^-}} $ to be the only element of the dual $(C^{1,\alpha}(\partial\Omega))'$ that satisfies the following equality
 \begin{equation}\label{defn:conoederd1}
<\frac{\partial u}{\partial\nu_{\Omega^-}},v>\equiv\int_{\partial\Omega}uS_-[v]\,d\sigma
\qquad\forall v\in C^{1,\alpha}(\partial\Omega)\,.
\end{equation}
\end{definition}
Here we have introduced the Definition  \ref{defn:conoederd} for harmonic functions of $C^{0}(\overline{\Omega^-})$ that are harmonic at infinity, but one could do the same also for harmonic functions in $\Omega^-$ that are harmonic at infinity in other function spaces that have a  trace operator on $\partial\Omega$ and where we do not have information on the integrability of the first order partial derivatives of $u$ in $\Omega$.\par 

   By applying  the classical second Green Identity for harmonic functions in $\Omega^-$ that are harmonic at infinity to the function $u$ and to the function ${\mathcal{G}}_{d,-}[v]$   (cf.~\textit{e.g.}, \cite[Cor.~4.8]{DaLaMu21}), one can show that if $ u\in C^{1}(\overline{\Omega^-})$ is harmonic in $\Omega^-$ and harmonic at infinity, then 
 \[
 <\frac{\partial u}{\partial\nu_{\Omega^-}},v>=\int_{\partial\Omega}\frac{\partial u}{\partial\nu_{\Omega^-}}  \,v\,d\sigma\qquad\forall v\in C^{1,\alpha}(\partial\Omega)\,,
 \]
  where $\frac{\partial u}{\partial\nu_{\Omega^-}} $ in the right hand side coincides with the classical normal derivative of $u$, \textit{i.e.}, the distributional  normal derivative $\frac{\partial u}{\partial\nu_{\Omega^-}} $ of definition (\ref{defn:conoederd1}) coincides with the distribution of $(C^{1,\alpha}(\partial\Omega))'$ that is associated to the classical normal derivative of $u$.
  Next we point out the validity of the following statement on the continuity of the (distributional) normal derivative.

 \begin{theorem}\label{thm:contednu}
 Let $\alpha\in ]0,1[$. Let $\Omega$ be a bounded open  subset of ${\mathbb{R}}^{n}$ of class $C^{1,\alpha}$. Then the following statements hold.
 \begin{enumerate}
\item[(i)] The transpose operator $S_-^t$ of the map $S_-$ is linear and continuous from 
\begin{equation}\label{eq:trespo}
\left(C^{0,\alpha}(\partial\Omega)\right)'\quad\text{to}\quad \left(C^{1,\alpha}(\partial\Omega)\right)'\,.
\end{equation}
\item[(ii)] The map $\frac{\partial}{\partial\nu_{\Omega^-}}$ from the closed subspace
$
C^{0,\alpha}_{bh}(\overline{\Omega^-}) $ of $ C^{0,\alpha}_b(\overline{\Omega^-})$ to $\left(C^{1,\alpha}(\partial\Omega)\right)'$ is linear and continuous. Moreover,
\begin{equation}\label{eq:trespo2}
\frac{\partial u}{\partial\nu_{\Omega^-}}=S_-^t[u_{|\partial\Omega}]
\qquad\forall u\in C^{0,\alpha}_{bh}(\overline{\Omega^-})\,.
\end{equation}
\end{enumerate}
\end{theorem}
{\bf Proof.} The continuity of the   operator $S_-^t$ follows by the continuity of  the linear operator  $S_-$ from $C^{1,\alpha}(\partial\Omega)$ to $C^{0,\alpha}(\partial\Omega)$. Since
 \[
 <S_-^t[u_{|\partial\Omega}],v>=<u_{|\partial\Omega},S_-[v]>=\int_{\partial\Omega}uS_-[v]\,d\sigma
 =<\frac{\partial u}{\partial\nu_{\Omega^-}},v>\ \forall v\in C^{1,\alpha}(\partial\Omega)\,,
\]
for all $u\in C^{0,\alpha}_{bh}(\overline{\Omega^-})$, equality (\ref{eq:trespo2}) holds true. Since the restriction map is continuous from $C^{0,\alpha}_{bh}(\overline{\Omega^-})$ to $C^{0,\alpha}(\partial\Omega)$  and the canonical embedding of $C^{0,\alpha}(\partial\Omega)$ into the dual $\left(C^{
{\color{black}0	}
,\alpha}(\partial\Omega)\right)'$ is continuous, the continuity of $S_-^t$ as in statement (i) implies the validity of (ii).\hfill  $\Box$ 

\vspace{\baselineskip}  
    
 In other words, statement (ii) says that the
distributional  normal derivative of $u$ coincides with the transpose of the classical exterior Steklov operator applied to the distribution that is canonically associated to the restriction of $u$ to $\partial\Omega$. Once more, one 
  could do the same also for harmonic functions in $\Omega^-$ that are harmonic at infinity  in function spaces that have a  trace operator on $\partial\Omega$, even in cases where we do not have information on the integrability of the first order partial derivatives of $u$ in $\Omega^-$.\par 
  
  We remark that the definition of distributional normal derivative that we have introduced is rather rigid. Indeed, we cannot take the normal derivative of a nonzero constant function in $\Omega^-$ for $n\geq 3$. Indeed, nonzero constant functions in $\Omega^-$ are not harmonic at infinity. Nonetheless such a definition turns out to be effective in the treatment of exterior boundary value problems in a H\"{o}lder space setting.

   \section{A nonvariational form of the interior Neumann problem for harmonic functions}
   \label{sec:nvneupb}
 
 Let $\alpha\in]0,1[$.  Let  $\Omega$ be a  bounded open subset of ${\mathbb{R}}^{n}$ of class $C^{1,\alpha}$. By exploiting the Definition   \ref{defn:conoderd}  of distributional  normal derivative, we can state the following Neumann problem. Given $g\in (C^{1,\alpha}(\partial\Omega))'$, find all $u\in C^{0}(\overline{\Omega})$ such that the following interior Neumann problem is satisfied.
 \begin{equation}\label{eq:nvneupb}
\left\{
\begin{array}{ll}
 \Delta u=0 & \text{in}\ \Omega\,,
 \\
 \frac{\partial u}{\partial\nu} =g& \text{in}\ (C^{1,\alpha}(\partial\Omega))'\,,
\end{array}
\right.
\end{equation}
where $ \frac{\partial u}{\partial\nu} $ is as in  Definition  \ref{defn:conoderd}. 
Since the solutions of the Neumann problem (\ref{eq:nvneupb}) may well have infinite Dirichlet integral, we address to problem (\ref{eq:nvneupb}) as `nonvariational interior Neumann problem'.\par

For the interior nonvariational Neumann problem  to have solutions, the boundary datum $g$ has to satisfy certain compatibility conditions that are akin to the corresponding compatibility conditions for the variational Neumann problem, as we show in the following Lemma (see Section \ref{sec:prelnot} for the notation on the connected components $\Omega_j$ of $\Omega$).
\begin{lemma}\label{lem:nvINcc}
Let $\alpha\in]0,1[$.  Let  $\Omega$ be a  bounded open subset of ${\mathbb{R}}^{n}$ of class $C^{1,\alpha}$. 
Let $g\in (C^{1,\alpha}(\partial\Omega))'$. If the interior Neumann problem (\ref{eq:nvneupb})
 has a solution $u\in C^{
 {\color{black}	0	}
 }(\overline\Omega)$, then 
\begin{equation}\label{lem:nvINcc2}
<g,\chi_{\partial\Omega_j}>=0\qquad\forall j\in \{1,\dots,\kappa^+\}\,.
\end{equation}
\end{lemma}
{\bf Proof.} First we note that $\chi_{\partial\Omega_j}$ is locally constant on $\partial\Omega$ and that accordingly $\chi_{\partial\Omega_j}\in C^{1,\alpha}(\partial\Omega)$  for all $j\in \{1,\dots,\kappa^+\}$. Next we note that $\chi_{\overline{\Omega_j}}$ solves  the Dirichlet problem (\ref{defn:cinspo3}) with $v=\chi_{\partial\Omega_j}$ and that accordingly ${\mathcal{G}}_{d,+}[\chi_{\partial\Omega_j}]=\chi_{\overline{\Omega_j}}$. Hence, the validity of the  interior Neumann problem (\ref{eq:nvneupb}) implies that
\begin{eqnarray*}
\lefteqn{
<g,\chi_{\partial\Omega_j}>=<\frac{\partial u}{\partial\nu},\chi_{\partial\Omega_j}>\equiv\int_{\partial\Omega}u\frac{\partial{\mathcal{G}}_{d,+}[\chi_{\partial\Omega_j}]}{\partial\nu}\,d\sigma
}
\\ \nonumber
&&\qquad\qquad
=\int_{\partial\Omega}u\frac{\partial \chi_{\overline{\Omega_j}}}{\partial\nu}\,d\sigma
=\int_{\partial\Omega}u0\,d\sigma=0 
\qquad\forall j\in \{1,\dots,\kappa^+\}\,.
\end{eqnarray*}
 \hfill  $\Box$ 

\vspace{\baselineskip}

Next we show that the possible solutions of the nonvariational  interior Neumann problem (\ref{eq:nvneupb}) are unique up to locally constant functions, exactly as in the classical case. To do so, we prove the following statement, that we prove by an argument which cannot be the classical energy argument.
\begin{theorem}\label{thm:ivneu.uni}
 Let $\alpha\in]0,1[$.  Let  $\Omega$ be a  bounded open subset of ${\mathbb{R}}^{n}$ of class $C^{1,\alpha}$. Let $g\in (C^{1,\alpha}(\partial\Omega))'$.
 
 If $u_{1}$, $u_{2}\in C^{0}(\overline{\Omega})$ solve the interior Neumann problem (\ref{eq:nvneupb}), then $u_{1}-u_{2}$ is constant  in each connected component of $\Omega$. In particular, all solutions of the interior Neumann problem in $C^{0}(\overline{\Omega})$ can be obtained by adding to $u_{1}$ an arbitrary function which is constant on the closure of each connected component of $\Omega$.
\end{theorem}
{\bf Proof.} Let $u\equiv u_{1}-u_{2}$. Since $\frac{\partial u}{\partial\nu} =0$ in $ (C^{1,\alpha}(\partial\Omega))'$, we have
\begin{equation}\label{thm:ivneu.uni1}
\int_{\partial\Omega}u\frac{\partial}{\partial\nu}
{\mathcal{G}}_{d,+}[v]
\,d\sigma=0\qquad\forall v\in C^{1,\alpha}(\partial\Omega)\,,
\end{equation}
{\color{black}(cf.~ Definition  \ref{defn:conoderd}). 	}
Now let
\begin{equation}\label{thm:ivneu.uni2}
C^{0,\alpha}(\partial\Omega)_{+}\equiv\left\{
h\in C^{0,\alpha}(\partial\Omega):\,\int_{\partial\Omega_j}h\,d\sigma=0\ \forall j\in\{1,\dots,\kappa^+\}\right\}\,.
\end{equation}
If $h\in C^{0,\alpha}(\partial\Omega)_{+}$, 
 then there exists ${\mathcal{N}}_+[h] 
\in C^{1,\alpha}(\overline{\Omega})$ such that
\begin{equation}\label{thm:ivneu.uni3}
\left\{
\begin{array}{ll}
 \Delta {\mathcal{N}}_+[h]=0 & \text{in}\ \Omega\,,
 \\
 \frac{\partial  }{\partial\nu}{\mathcal{N}}_+[h] =h& \text{on}\  \partial\Omega\,,
\end{array}
\right.
\end{equation}
(cf.~\textit{e.g.}, \cite[Thm.~6.42]{DaLaMu21}). Clearly, ${\mathcal{G}}_{d,+}\left[{\mathcal{N}}_+[h]_{|\partial\Omega}\right]={\mathcal{N}}_+[h]$  
and  condition (\ref{thm:ivneu.uni1}) with $v= {\mathcal{N}}_+[h]_{|\partial\Omega}$ implies that
\begin{equation}\label{thm:ivneu.uni4}
\int_{\partial\Omega}uh\,d\sigma=\int_{\partial\Omega}u\frac{\partial }{\partial\nu}{\mathcal{N}}_+[h]\,d\sigma=
\int_{\partial\Omega}u\frac{\partial {\mathcal{G}}_{d,+}\left[{\mathcal{N}}_+[h]_{|\partial\Omega}\right]}{\partial\nu}\,d\sigma=0 \,.
\end{equation}
Now let $j\in\{1,\dots,\kappa^+\}$. By condition (\ref{thm:ivneu.uni4}) for all $h\in C^{0,\alpha}(\partial\Omega)_{+}$, we have
\begin{equation}\label{thm:ivneu.uni5}
\int_{\partial\Omega_j}uh\,d\sigma=0\qquad\forall h\in C^{0,\alpha}(\partial\Omega_j)_0\,, 
\end{equation}
{\color{black} (cf.~(\ref{eq:x0})).	}	
Now let $\phi_j\in  C^{0,\alpha}(\partial\Omega_j)$ be such that $\int_{\partial\Omega_j}\phi_j\,d\sigma=1$. Then we have
\begin{eqnarray*}
\lefteqn{
\int_{\partial\Omega_j}uf\,d\sigma=\int_{\partial\Omega_j}u\left(f-\int_{\partial\Omega_j}f\,d\sigma\phi_j\right)\,d\sigma
+\int_{\partial\Omega_j}u\int_{\partial\Omega_j}f\,d\sigma\phi_j\,d\sigma
}
\\ \nonumber
&&\qquad\qquad\qquad\qquad
=0+\int_{\partial\Omega_j} \left(\int_{\partial\Omega_j}u\phi_j\,d\sigma\right) f\,d\sigma\qquad\forall f\in C^{0,\alpha}(\partial\Omega_j)\,,
\end{eqnarray*}
and accordingly
\[
u=\int_{\partial\Omega_j}u\phi_j\,d\sigma\qquad\text{on}\ \partial\Omega_j\,.
\]
Hence, $u$ is constant on $\partial\Omega_j$ for all $j\in\{1,\dots,\kappa^+\}$. Since $u$ is harmonic in $\Omega_j$ and continuous on $\overline{\Omega_j}$, we conclude that $u$ is constant on $\overline{\Omega_j}$ for all $j\in\{1,\dots,\kappa^+\}$.

\hfill  $\Box$ 

\vspace{\baselineskip}

 \section{A nonvariational form of the exterior Neumann problem for harmonic functions}
 \label{sec:nveneupb}
 Let $\alpha\in]0,1[$.  Let  $\Omega$ be a  bounded open subset of ${\mathbb{R}}^{n}$ of class $C^{1,\alpha}$. By exploiting the Definition   \ref{defn:conoederd}  of distributional  normal derivative, we can state the following Neumann problem in the exterior domain $\Omega^-$. Given $g\in (C^{1,\alpha}(\partial\Omega))'$, find all $u\in C^{0}_b(\overline{\Omega^-})$ such that the following exterior Neumann problem is satisfied.
  \begin{equation}\label{eq:nveneupb}
\left\{
\begin{array}{ll}
 \Delta u=0 & \text{in}\ \Omega\,,
 \\
 -\frac{\partial u}{\partial\nu_{\Omega^-}}  =g& \text{in}\ (C^{1,\alpha}(\partial\Omega))'\,,
 \\
 u\ \text{is\ harmonic\ at}\ \infty\,, &
\end{array}
\right.
\end{equation}
where $ \frac{\partial u}{\partial\nu_{\Omega^-}} $ is as in Definition  \ref{defn:conoederd}. 
Since the solutions of the Neumann problem (\ref{eq:nveneupb}) may well have infinite Dirichlet integral, we address to problem (\ref{eq:nveneupb}) as `nonvariational exterior Neumann problem'.\par

 For the   nonvariational exterior Neumann problem  to have solutions, the boundary datum $g$ has to satisfy certain compatibility conditions that are akin to the corresponding compatibility conditions for the exterior variational Neumann problem, as we show in the following Lemma (see Section \ref{sec:prelnot} for the notation on the connected components $(\Omega^-)_j$ of $\Omega^-$).
\begin{lemma}\label{lem:nvENcc}
Let $\alpha\in]0,1[$.  Let  $\Omega$ be a  bounded open subset of ${\mathbb{R}}^{n}$ of class $C^{1,\alpha}$. 
Let $g\in (C^{1,\alpha}(\partial\Omega))'$. If the exterior Neumann problem (\ref{eq:nveneupb})
 has a solution $u\in C^{0}_{b}(\overline{\Omega^-})$, then 
\begin{equation}\label{lem:nvENcc.eq2}
<g,\chi_{\partial(\Omega^-)_j}>=0\qquad\forall j\in \{1,\dots,\kappa^-\}\,.
\end{equation}
In addition, if $n=2$, then we also have 
\begin{equation}\label{lem:nvENcc.eq3}
<g,\chi_{\partial(\Omega^-)_0}>=0\,.
\end{equation}
\end{lemma}
{\bf Proof.} First we note that $\chi_{\partial(\Omega^-)_j}$ is locally constant on $\partial\Omega^-=\partial\Omega$ and that accordingly $\chi_{\partial(\Omega^-)_j}\in C^{1,\alpha}(\partial\Omega)$  for all $j\in \{0,\dots,\kappa^-\}$. 
 
 Next we note that $\chi_{\overline{(\Omega^-)_j}}$ solves  the Dirichlet problem (\ref{defn:censpo3}) with $v=\chi_{\partial(\Omega^-)_j}$ and that accordingly ${\mathcal{G}}_{d,-}[\chi_{\partial(\Omega^-)_j}]=\chi_{\overline{(\Omega^-)_j}}$
 for all $j\in \{1,\dots,\kappa^-\}$ in case $n\geq 2$ and also for $j=0$ in case $n=2$. Indeed, nonzero constant functions are harmonic at infinity only in case $n=2$. Hence, the validity of the  exterior Neumann problem (\ref{eq:nveneupb}) implies that
\begin{eqnarray*}
\lefteqn{
<g,\chi_{\partial(\Omega^-)_j}>=<-\frac{\partial u}{\partial
\nu_{\Omega^-}} ,\chi_{\partial(\Omega^-)_j}>
=-\int_{\partial\Omega}u\frac{\partial {\mathcal{G}}_{d,-}[\chi_{\partial(\Omega^-)_j}]}{\partial \nu_{\Omega^-}
}\,d\sigma
}
\\ \nonumber
&&\qquad\qquad\qquad\qquad\qquad\qquad
=-\int_{\partial\Omega}u\frac{\partial \chi_{\overline{(\Omega^-)_j}}}{\partial\nu_{\Omega^-}}\,d\sigma
=-\int_{\partial\Omega}u0\,d\sigma=0 
\end{eqnarray*}
for all $j\in \{1,\dots,\kappa^-\}$ in case $n\geq 2$ and also for $j=0$ in case $n=2$.\hfill  $\Box$ 

\vspace{\baselineskip}

 Next we show that the possible solutions of the nonvariational exterior Neumann problem (\ref{eq:nveneupb}) are unique up to locally constant functions, exactly as in the classical case. To do so, we prove the following statement, that we prove by an argument which cannot be the classical energy argument.
 \begin{theorem}\label{thm:evneu.uni}
 Let $\alpha\in]0,1[$.  Let  $\Omega$ be a  bounded open subset of ${\mathbb{R}}^{n}$ of class $C^{1,\alpha}$. Let $g\in (C^{1,\alpha}(\partial\Omega))'$.
 
 If $u_{1}$, $u_{2}\in C^{0}_{b}(\overline{\Omega^{-}})$ solve the exterior Neumann problem (\ref{eq:nveneupb}), {\color{black} then	}  $u_{1}-u_{2}$ is constant  in each connected component of $\Omega^-$. Moreover, if $n\ge 3$, then $u_{1}-u_{2}=0$ in the unbounded connected component $(\Omega^{-})_{0}$. In particular, all the solutions of the exterior Neumann problem  (\ref{eq:nveneupb}) in $C^{0}_{b}(\overline{\Omega^{-}})$ {\color{black} can }
 be obtained by adding to $u_{1}$ an arbitrary function  which is constant on the closure of each connected component of $\Omega^{-}$ and equals zero in $\overline{(\Omega^{-})_{0}}$ when $n\geq 3$.
\end{theorem}
{\bf Proof.} Let $u\equiv u_{1}-u_{2}$. Since $\frac{\partial u}{\partial\nu_{\Omega^-}} =0$ in $ (C^{1,\alpha}(\partial\Omega))'$, we have
\begin{equation}\label{thm:evneu.uni1}
\int_{\partial\Omega}u\frac{\partial }{\partial\nu_{\Omega^-}}{\mathcal{G}}_{d,-}[v]
\,d\sigma=0\qquad\forall v\in C^{1,\alpha}(\partial\Omega)\,,
\end{equation}
{\color{black}(cf.~Definition  \ref{defn:conoederd}).	}
Now let
\begin{eqnarray}\label{thm:evneu.uni2}
\lefteqn{C^{0,\alpha}(\partial\Omega)_{-}
}
\\ \nonumber
&&\qquad
\equiv\left\{
h\in C^{0,\alpha}(\partial\Omega):\,\int_{\partial(\Omega^-)_j}h\,d\sigma=0\ \forall j\in\{1,\dots,\kappa^-\}\right\}
\quad\text{if}\ n\geq 3\,,
\\ \nonumber
\lefteqn{C^{0,\alpha}(\partial\Omega)_{-}
}
\\ \nonumber
&&\qquad
\equiv\left\{
h\in C^{0,\alpha}(\partial\Omega):\,\int_{\partial(\Omega^-)_j}h\,d\sigma=0\ \forall j\in\{0,\dots,\kappa^-\}\right\}
\quad\text{if}\ n= 2\,.
\end{eqnarray}
If $h\in C^{0,\alpha}(\partial\Omega)_{-}$, 
 then there exists ${\mathcal{N}}_-[h]\in C^{1,\alpha}_{b}(\overline{\Omega^-})$ such that
\begin{equation}\label{thm:evneu.uni3}
\left\{
\begin{array}{ll}
 \Delta {\mathcal{N}}_-[h]=0 & \text{in}\ \Omega^-\,,
 \\
 -\frac{\partial {\mathcal{N}}_-[h]}{\partial\nu_{\Omega^-}} =h& \text{on}\  \partial\Omega\,,
 \\
 {\mathcal{N}}_-[h]\ \text{is\ harmonic\ at}\ \infty
\end{array}
\right.
\end{equation}
(cf.~\textit{e.g.}, \cite[Thms.~6.43, 6.44]{DaLaMu21} together with the classical  Theorem  \ref{thm:slayh},  on the membership in $C^{1,\alpha}_{b}(\overline{\Omega^-})$
of the single  layer potential). Clearly, ${\mathcal{G}}_{d,-}\left[{\mathcal{N}}_-[h]_{|\partial\Omega}\right]={\mathcal{N}}_-[h]$  
and  condition (\ref{thm:evneu.uni1}) with $v= {\mathcal{N}}_-[h]_{|\partial\Omega}$ implies that
\begin{equation}\label{thm:evneu.uni4}
\int_{\partial\Omega}uh\,d\sigma=
-\int_{\partial\Omega}u\frac{\partial {\mathcal{N}}_-[h]}{\partial\nu_{\Omega^-}}\,d\sigma=
-\int_{\partial\Omega}u\frac{\partial {\mathcal{G}}_{d,-}\left[{\mathcal{N}}_-[h]_{|\partial\Omega}\right]}{\partial\nu}\,d\sigma=0 \,.
\end{equation}
Now let $j\in\{1,\dots,\kappa^-\}$ if $n\geq 3$ and $j\in\{0,\dots,\kappa^-\}$ if $n=2$. By condition (\ref{thm:evneu.uni4}) for all $h\in C^{0,\alpha}(\partial\Omega)_{-}$, we have
\begin{equation}\label{thm:evneu.uni5}
\int_{\partial(\Omega^-)_j}uh\,d\sigma=0\qquad\forall h\in C^{0,\alpha}( 
{\color{black}\partial(\Omega^-)_j	}
)_0\,,
\end{equation}
{\color{black} (cf.~(\ref{eq:x0})).	}
Now let $\phi_j\in  C^{0,\alpha}(\partial(\Omega^-)_j)$ be such that $\int_{\partial(\Omega^-)_j}\phi_j\,d\sigma=1$. Then we have
\begin{eqnarray*}
\lefteqn{
\int_{\partial(\Omega^-)_j}uf\,d\sigma
}
\\ \nonumber
&&\qquad 
=\int_{\partial(\Omega^-)_j}u\left(f-\int_{\partial(\Omega^-)_j}f\,d\sigma\phi_j\right)\,d\sigma
+\int_{\partial(\Omega^-)_j}u\int_{\partial(\Omega^-)_j}f\,d\sigma\phi_j\,d\sigma
\\ \nonumber
&&\qquad 
=0+\int_{\partial(\Omega^-)_j} \left(\int_{\partial(\Omega^-)_j}u\phi_j\,d\sigma\right) f\,d\sigma\qquad\forall f\in C^{0,\alpha}(\partial(\Omega^-)_j)\,,
\end{eqnarray*}
and accordingly
\[
u=\int_{\partial(\Omega^-)_j}u\phi_j\,d\sigma\qquad\text{on}\ \partial(\Omega^-)_j\,.
\]
Hence, $u$ is constant on $\partial(\Omega^-)_j$ for all $j\in\{1,\dots,\kappa^-\}$ if $n\geq 3$ and for all $j\in\{0,\dots,\kappa^-\}$ if $n=2$.  {\color{black}
Next we note that in case $n\geq3$ condition (\ref{thm:evneu.uni4}) for all $h\in C^{0,\alpha}(\partial\Omega)_{-}$ implies that
\[
\int_{\partial(\Omega^-)_0}uh\,d\sigma=0\qquad\forall h\in C^{0,\alpha}( 
{\color{black}\partial(\Omega^-)_0	}
)\,,
\]
and accordingly that
\[
u=0 \qquad\text{on}\ \partial(\Omega^-)_0\,.
\]
}
Since $u$ is harmonic in $(\Omega^-)_j$ and continuous on $\overline{(\Omega^-)_j}$
{\color{black}for all $j\in\{1,\dots,\kappa^-\}$ and harmonic, harmonic at infinity and continuous on $\overline{(\Omega^-)_0}$, the classical uniqueness Theorem for the exterior Dirichlet problem of Theorem \ref{thm:existenceO-0} of the Appendix implies that 
$u$ is constant on $\overline{(\Omega^-)_j}$ for all $j\in\{1,\dots,\kappa^-\}$ and vanishes on $\overline{(\Omega^-)_0}$ if $n\geq 3$ and 
is constant on $\overline{(\Omega^-)_j}$ for all $j\in\{0,\dots,\kappa^-\}$ if $n=2$. 
}
\hfill  $\Box$ 

\vspace{\baselineskip}

  \section{Preliminaries on the  classical Green function for the Dirichlet problem in bounded open sets of class $C^{1,\alpha}$}
\label{sec:pokein1a}
The present section presents some classical known material on the Green 
{\color{black}function	}  
for the Dirichlet problem in    bounded open sets of class $C^{1,\alpha}$ that we need in the sequel.  By the unique solvability of the Dirichlet problem in a  bounded open set  of class $C^{1,\alpha}$, the following classical Lemma is known to hold (cf.~\textit{e.g.}, Theorem \ref{thm:idwp}).
\begin{lemma}\label{lem:ihx}
 Let $\alpha\in ]0,1[$. Let $\Omega$ be a bounded open  subset of ${\mathbb{R}}^{n}$ of class $C^{1,\alpha}$.
 Then for each $x\in {\mathbb{R}}^n\setminus\partial\Omega$, there exists a unique $h_x\in C^{1,\alpha}(\overline{\Omega})$ that solves the following interior Dirichlet problem
\begin{equation}\label{lem:ihx1}
\left\{
\begin{array}{ll}
\Delta_{y} h_{x}(y)=0 & \forall y\in \Omega\,,
\\
h_{x}(y)=S_{n}(x-y) &  \forall y\in  \partial\Omega\,.
\end{array}
\right.
\end{equation}
\end{lemma}
Then if $u\in C^{1,\alpha}(\overline{\Omega})$ is harmonic in $\Omega$, one can 
 apply  the second Green identity to the pair of functions $u$, $v(\cdot)\equiv h_{x}(\cdot)$, and obtain
 \[
\int_{\partial\Omega}u(y)\frac{\partial}{\partial\nu_{\Omega,y}}h_{x}(y)-
\frac{\partial u}{\partial\nu_{\Omega} }(y) h_{x}(y)\,d\sigma_{y}=0 
\]
(cf.~\textit{e.g.}, \cite[Thm.~4.3]{DaLaMu21}).  Then the third Green Identity for $u$ implies the validity of
   the following statement (cf.~\textit{e.g.}, \cite[Cor.~4.6]{DaLaMu21}).
\begin{theorem}\label{thm:igrfu}
 Let $\alpha\in ]0,1[$. Let $\Omega$ be a bounded open  subset of ${\mathbb{R}}^{n}$ of class $C^{1,\alpha}$.
 If $u\in C^{1,\alpha}(\overline{\Omega})$ is harmonic in $\Omega$, then
 \begin{eqnarray}\label{thm:igrfu1}
u(x)=\int_{\partial\Omega}
u(y)\frac{\partial}{\partial\nu_{\Omega,y} }(S_{n} (x-y)-h_{x}(y))\,d\sigma_{y}
\qquad\forall x\in\Omega\,,
\\ \nonumber
0=\int_{\partial\Omega}
u(y)\frac{\partial}{\partial\nu_{\Omega,y} }(S_{n} (x-y)-h_{x}(y))\,d\sigma_{y}
\qquad\forall x\in\Omega^-\,.
\end{eqnarray}
\end{theorem}
The function $G_+(\cdot,\cdot)$  from $({\mathbb{R}}^n\setminus\partial\Omega)\times\overline{\Omega}\setminus\{(x,x):\,x\in{\mathbb{R}}^n\setminus\partial\Omega\}$ to ${\mathbb{R}}$ defined by 
\[
G_+(x,y)\equiv  S_{n} (x-y)-h_{x}(y)
\quad\forall (x,y)\in \left(
({\mathbb{R}}^n\setminus\partial\Omega)
\times\overline{\Omega}
\right)
\setminus\{(x,x):\,x\in{\mathbb{R}}^n\setminus\partial\Omega\} \,,
\]
is well known to be the Green function for the interior Dirichlet problem in the open set $\Omega$. By the classical existence and uniqueness Theorem for the Dirichlet problem with datum of class $C^{1,\alpha}$ (cf.~\textit{e.g.}, Theorem \ref{thm:idwp}), one can deduce the validity of the following statement.
\begin{theorem}[of the Green function]\label{thm:inpokerco}
 Let $\alpha\in ]0,1[$. Let $\Omega$ be a bounded open  subset of ${\mathbb{R}}^{n}$ of class $C^{1,\alpha}$.
 If $g\in C^{1,\alpha}(\partial\Omega)$, then the function $u_g$ defined by 
 \begin{equation}\label{thm:inpokerco1}
u_g(x)\equiv
\left\{
\begin{array}{ll}
 \int_{\partial\Omega}
g(y)\frac{\partial}{\partial\nu_{\Omega,y} }(S_{n} (x-y)-h_{x}(y))\,d\sigma_{y}
&  \forall x\in\Omega\,,
\\
g(x) &  \forall x\in\partial\Omega\,,
\end{array}
\right.
\end{equation}
is of class $C^{1,\alpha}(\overline{\Omega})$ and is harmonic in $\Omega$. Moreover,
\begin{equation}\label{thm:inpokerco1a}
\int_{\partial\Omega}
g(y)\frac{\partial}{\partial\nu_{\Omega,y} }(S_{n} (x-y)-h_{x}(y))\,d\sigma_{y}=0\qquad\forall x\in\Omega^-\,.
\end{equation}
\end{theorem}
{\bf Proof.} By the classical existence and uniqueness Theorem for the Dirichlet problem with datum of class $C^{1,\alpha}$ (cf.~\textit{e.g.}, Theorem \ref{thm:idwp}), there exists a unique $u\in C^{1,\alpha}(\overline{\Omega})$ that is harmonic in $\Omega$ and that satisfies equality $u_{|\partial\Omega}=g$. Then Theorem \ref{thm:igrfu} implies that $u$ equals $u_g$ and  that equality (\ref{thm:inpokerco1a}) holds true.  Accordingly, $u_g\in C^{1,\alpha}(\overline{\Omega})$, $u_g$ is harmonic in $\Omega$ and $u_g(x)=g(x)$ for all $x\in\partial\Omega$.\hfill  $\Box$ 

\vspace{\baselineskip}

Then one can invoke a   classical  approximation argument (cf.~\textit{e.g.}, Folland~\cite[Prop.~7.38]{Fo95}) and prove the following. For the convenience of the reader we include a proof.
\begin{theorem}[of the Green function]\label{thm:inpokerc0}
 Let $\alpha\in ]0,1[$. Let $\Omega$ be a bounded open  subset of ${\mathbb{R}}^{n}$ of class $C^{1,\alpha}$.
 If $g\in C^{0,\alpha}(\partial\Omega)$, then the function $u_g$ defined by 
 \begin{equation}\label{thm:inpokerc01}
u_g(x)\equiv
\left\{
\begin{array}{ll}
 \int_{\partial\Omega}
g(y)\frac{\partial}{\partial\nu_{\Omega,y} }(S_{n} (x-y)-h_{x}(y))\,d\sigma_{y}
&  \forall x\in\Omega\,,
\\
g(x) &  \forall x\in\partial\Omega\,,
\end{array}
\right.
\end{equation}
is of class $C^{0,\alpha}(\overline{\Omega})$ and is harmonic in $\Omega$.  
Moreover, equality (\ref{thm:inpokerco1a}) holds true and $u_g={\mathcal{G}}_{d,+}[g]$.
 \end{theorem}
{\bf Proof.} By the Weierstrass approximation Theorem, there exists a sequence $\{g_l\}_{l\in {\mathbb{N}}}$ of polynomials in $n$ variables such that the sequence $\{g_{l|\partial\Omega}\}_{l\in {\mathbb{N}}}$
  converges uniformly to $g$ on $\partial\Omega$. Since $g_l\in C^{1,\alpha}(\partial\Omega)$,  Theorem \ref{thm:inpokerco} implies that 
$u_{g_l}\in C^{1,\alpha}(\overline{\Omega})$ is harmonic in $\Omega$ and  $u_{g_l}(x)=g_l(x)$ for all $x\in\partial\Omega$ and $l\in {\mathbb{N}}$ and  that formulas   {\color{black}(\ref{thm:inpokerco1}), 	}
(\ref{thm:inpokerco1a})  hold for $g_l$ for all
 $l\in {\mathbb{N}}$.

By the Maximum Principle, we have
\[
\sup_{x\in \overline{\Omega}}|u_{g_{l_1}}(x)-u_{g_{l_2}}(x)|\leq \sup_{x\in\partial\Omega}|g_{l_1}(x)-g_{l_2} (x)|
\]
for all $l_1$, $l_2\in {\mathbb{N}}$ 
 and accordingly the sequence $\{u_{g_{l}}\}_{l\in {\mathbb{N}}}$ has a uniform limit $u\in C^0(\overline{\Omega})$ in $\overline{\Omega}$. 
Then $u$ is harmonic in $\Omega$ and $u_{|\partial\Omega}=g$. Moreover, if $x\in \Omega$, then we have
\begin{eqnarray*}
\lefteqn{
u(x)=\lim_{l\to\infty}u_{g_{l}}(x)= \int_{\partial\Omega}
g_{l}(y)\frac{\partial}{\partial\nu_{\Omega,y} }(S_{n} (x-y)-h_{x}(y))\,d\sigma_{y}
}
\\ \nonumber
&&\qquad
=\int_{\partial\Omega}
g(y)\frac{\partial}{\partial\nu_{\Omega,y} }(S_{n} (x-y)-h_{x}(y))\,d\sigma_{y}\equiv u_g(x)\,.
\end{eqnarray*}
 Hence, $u_g$ is of class $C^{0}(\overline{\Omega})$,  is harmonic in $\Omega$ and $u_g(x)=g(x)$ for all $x\in\partial\Omega$. Then  the existence  Theorem \ref{thm:existenceO0} of the Appendix for the Dirichlet problem implies that there exists   $u^\sharp \equiv {\mathcal{G}}_{d,+}[g]
 \in C^{0,\alpha}(\overline{\Omega})$ that  is harmonic in $\Omega$ and satisfies the equality $u^\sharp(x)=g(x)$ for all $x\in\partial\Omega$. Then the Maximum Principle implies that $u_g=u^\sharp\in C^{0,\alpha}(\overline{\Omega})$ and thus the proof is complete.\hfill  $\Box$ 

\vspace{\baselineskip}

\section{Preliminaries on the classical Green function for the Dirichlet problem in  the exterior of a bounded open set of class $C^{1,\alpha}$}
\label{sec:pokeex1a}
The present section presents some classical known material on the  Green 
{\color{black}function	}
 for the Dirichlet problem in the   exterior of bounded open sets of class $C^{1,\alpha}$  that we need in the sequel.  By the unique solvability of the  Dirichlet problem in the exterior of a   bounded open set  of class $C^{1,\alpha}$, the following classical Lemma is known to hold (cf.~\textit{e.g.}, Theorem \ref{thm:edwp}).
\begin{lemma}\label{lem:ehx}
 Let $\alpha\in ]0,1[$. Let $\Omega$ be a bounded open  subset of ${\mathbb{R}}^{n}$ of class $C^{1,\alpha}$.
 Then for each $x\in {\mathbb{R}}^n\setminus\partial\Omega$, there exists a unique $h_x^-\in C^{1,\alpha}_{b}(\overline{\Omega^-})$ that solves the following exterior Dirichlet problem
\begin{equation}\label{lem:ehx1}
\left\{
\begin{array}{ll}
\Delta_{y} h_{x-}(y)=0 & \forall y\in \Omega^-\,,
\\
h_{x-}(y)=S_{n}(x-y) &  \forall y\in  \partial\Omega\,,
\\
h_{x-} \text{ is harmonic at infinity.} &
\end{array}
\right.
\end{equation}
\end{lemma}
Then if $u\in C^{1,\alpha}_{b}(\overline{\Omega^-})$ is harmonic in $\Omega^-$ and harmonic at infinity, one can 
 apply  the second Green identity for exterior domains to the pair of functions $u$, $v(\cdot)\equiv h_{x}(\cdot)$, and obtain
 \[
\int_{\partial\Omega}u(y)\frac{\partial}{\partial\nu_{\Omega,y}}h_{x-} (y)-
\frac{\partial u}{\partial\nu_{\Omega} }(y) h_{x-} (y)\,d\sigma_{y}=0 
\]
(cf.~\textit{e.g.}, \cite[Cor.~4.8]{DaLaMu21}). Then the third Green Identity in exterior domains  for $u$ implies the validity of
   the following statement (cf.~\textit{e.g.}, \cite[Cor.~3.24, Thm.~4.10]{DaLaMu21}).
\begin{theorem}\label{thm:egrfu}
 Let $\alpha\in ]0,1[$. Let $\Omega$ be a bounded open  subset of ${\mathbb{R}}^{n}$ of class $C^{1,\alpha}$.
 Let $u\in C^{1,\alpha}_{b}(\overline{\Omega^-})$ be harmonic in $\Omega^-$ and harmonic at infinity. Let $u_{\infty}$ be the limiting value of $u$ at infinity.  Then
 \begin{eqnarray}\label{thm:egrfu1}
u(x)&=&-\int_{\partial\Omega}
u(y)\frac{\partial}{\partial\nu_{\Omega,y} }(S_{n} (x-y)-h_{x-}(y))\,d\sigma_{y}+u_{\infty}
\quad\forall x\in\Omega^{-},
\\ \nonumber
0&=&-\int_{\partial\Omega}
u(y)\frac{\partial}{\partial\nu_{\Omega,y} }(S_{n} (x-y)-h_{x-}(y))\,d\sigma_{y}+u_{\infty}
\quad\forall x\in\Omega 
\end{eqnarray}
and
\begin{equation}\label{thm:egrfu1b}
u_\infty=\lim_{\xi\to\infty}u(\xi)=\left\{
\begin{array}{ll}
 0 &\quad\text{if}\ n\geq 3\,,
 \\
 \Xint{-}_{ \partial{\mathbb{B}}_{n}(0,r) }u \,d\sigma &\quad\text{if}\ n=2\,,
\end{array}
\right.
\end{equation}
for all $r\in]0,+\infty[$ such that $\overline{\Omega}\subseteq {\mathbb{B}}_{n}(0,r)$.
\end{theorem}
The function $G_-(\cdot,\cdot)$  from $\left(({\mathbb{R}}^n\setminus\partial\Omega)\times\overline{\Omega^-}\right)\setminus\{(x,x):\,x\in {\mathbb{R}}^n\setminus\partial\Omega\}$ to ${\mathbb{R}}$ defined by 
\[
G_-(x,y)\equiv  S_{n} (x-y)-h_{x-}(y)
\quad\forall (x,y)\in \left(({\mathbb{R}}^n\setminus\partial\Omega)\times\overline{\Omega^-}\right)\setminus\{(x,x):\,x\in {\mathbb{R}}^n\setminus\partial\Omega\}
\]
is well known to be the Green function for the exterior Dirichlet problem in the open set $\Omega^-$.
 By the existence and uniqueness Theorem   for the exterior Dirichlet problem (cf.~\textit{e.g.}, Theorem \ref{thm:edwp}), one can deduce the validity of the following statement.
\begin{theorem}[of the Green function]\label{thm:enpokerco}
 Let $\alpha\in ]0,1[$. Let $\Omega$ be a bounded open  subset of ${\mathbb{R}}^{n}$ of class $C^{1,\alpha}$.
 If $g\in C^{1,\alpha}(\partial\Omega)$, then there exists a unique $c_g\in {\mathbb{R}}$ such that  the function $u_{g-}$ defined by 
 \begin{equation}\label{thm:enpokerco1}
u_{g-}(x)\equiv
\left\{
\begin{array}{ll}
 -\int_{\partial\Omega}
g(y)\frac{\partial}{\partial\nu_{\Omega,y} }(S_{n} (x-y)-h_{x-}(y))\,d\sigma_{y}+c_g
&  \forall x\in\Omega^-\,,
\\
g(x) &  \forall x\in\partial\Omega\,,
\end{array}
\right.
\end{equation}
is of class $C^{1,\alpha}_{b}(\overline{\Omega^-})$, is harmonic in $\Omega^-$ and harmonic at infinity. Moreover,  
\begin{equation}\label{thm:enpokerco1a}
-\int_{\partial\Omega}
g(y)\frac{\partial}{\partial\nu_{\Omega,y} }(S_{n} (x-y)-h_{x-}(y))\,d\sigma_{y}+c_g
=0  \qquad\forall x\in\Omega\,,
\end{equation}
and
\begin{equation}\label{thm:enpokerco1b}
c_g=\lim_{\xi\to\infty}u_{g-}(\xi)=\left\{
\begin{array}{ll}
 0 &\quad\text{if}\ n\geq 3\,,
 \\
 \Xint{-}_{ \partial{\mathbb{B}}_{n}(0,r) }u_{g-}\,d\sigma &\quad\text{if}\ n=2\,,
\end{array}
\right.
\end{equation}
for all $r\in]0,+\infty[$ such that $\overline{\Omega}\subseteq {\mathbb{B}}_{n}(0,r)$.
\end{theorem}
{\bf Proof.} By Theorem  \ref{thm:edwp}, there exists a unique $u\in C^{1,\alpha}_{b}(\overline{\Omega^-})$ that is harmonic in $\Omega^-$, harmonic at infinity and that satisfies equality 
$u_{|\partial\Omega}=g$. Then Theorem \ref{thm:egrfu} implies that $u$ equals $u_{g-}$ with $c_g$ as in (\ref{thm:enpokerco1b}) and that (\ref{thm:enpokerco1a}) and holds  true. Hence,
  $u_{g-}\in C^{1,\alpha}_{b}(\overline{\Omega^-})$, $u_{g-}$ is harmonic in $\Omega^-$, harmonic at infinity and $u_{g-}(x)=g(x)$ for all $x\in\partial\Omega$.  \hfill  $\Box$ 

\vspace{\baselineskip}
 
 Then one can invoke a classical  approximation argument (cf.~\textit{e.g.}, Folland~\cite[Prop.~7.38]{Fo95}) and prove the following. For the convenience of the reader we include a proof.
\begin{theorem}[of the Green function]\label{thm:enpokerc0}
 Let $\alpha\in ]0,1[$. Let $\Omega$ be a bounded open  subset of ${\mathbb{R}}^{n}$ of class $C^{1,\alpha}$.
 If $g\in C^{0,\alpha}(\partial\Omega)$, then there exists a unique $c_g\in {\mathbb{R}}$ such that  the function $u_{g-}$ defined by  
 \begin{equation}\label{thm:enpokerc01}
u_{g-}(x)\equiv
\left\{
\begin{array}{ll}
-\int_{\partial\Omega}
g(y)\frac{\partial}{\partial\nu_{\Omega,y} }(S_{n} (x-y)-h_{x-}(y))\,d\sigma_{y}+c_g
&  \forall x\in\Omega^-\,,
\\
g(x) &  \forall x\in\partial\Omega\,,
\end{array}
\right.
\end{equation}
is of class $C^{0,\alpha}_{\mathrm{loc}}(\overline{\Omega^-})$,  is harmonic in $\Omega^-$ and harmonic at infinity. Moreover,
\begin{equation}\label{thm:enpokerc01a}
-\int_{\partial\Omega}
g(y)\frac{\partial}{\partial\nu_{\Omega,y} }(S_{n} (x-y)-h_{x-}(y))\,d\sigma_{y}+c_g
=0  \qquad\forall x\in\Omega\,,
\end{equation}
and
\begin{equation}\label{thm:enpokerc01b}
c_g=\lim_{\xi\to\infty}u_{g-}(\xi)=\left\{
\begin{array}{ll}
 0 &\quad\text{if}\ n\geq 3\,,
 \\
 \Xint{-}_{ \partial{\mathbb{B}}_{n}(0,r) }u_{g-}\,d\sigma &\quad\text{if}\ n=2\,,
\end{array}
\right.
\end{equation}
for all $r\in]0,+\infty[$ such that $\overline{\Omega}\subseteq {\mathbb{B}}_{n}(0,r)$   and $u_{g-}={\mathcal{G}}_{d,-}[g]$.
\end{theorem}
{\bf Proof.} By the Weierstrass approximation Theorem, there exists a sequence $\{g_l\}_{l\in {\mathbb{N}}}$ of polynomials in $n$ variables such that the sequence $\{g_{l|\partial\Omega}\}_{l\in {\mathbb{N}}}$
  converges uniformly to $g$ on $\partial\Omega$. Since $g_l\in C^{1,\alpha}(\partial\Omega)$,  Theorem \ref{thm:enpokerco} implies that 
$u_{g_l-}\in C^{1,\alpha}_{b}(\overline{\Omega^-})$ is harmonic in $\Omega^-$, harmonic at infinity, $u_{g_l-}(x)=g_l(x)$ for all $x\in\partial\Omega$ and  that formulas (\ref{thm:enpokerco1}), (\ref{thm:enpokerco1a}) hold for $g_l$ for all
 $l\in {\mathbb{N}}$. By the Maximum Principle for harmonic functions that are harmonic at infinity in an exterior domain, we have
\[
\sup_{x\in \overline{\Omega^-}}|u_{g_{l_1}-}(x)-u_{g_{l_2}-}(x)|\leq \sup_{x\in\partial\Omega}|g_{l_1-}(x)-g_{l_2-} (x)|
\]
for all $l_1$, $l_2\in {\mathbb{N}}$ 
 and accordingly the sequence $\{u_{g_{l}-}\}_{l\in {\mathbb{N}}}$ has a uniform limit $u_-\in C^0_b(\overline{\Omega^-})$ in $\overline{\Omega^-}$. In particular, $u_{-|\partial\Omega}=g$. 
 
 Since the uniform limit of harmonic functions is harmonic and the uniform  limit of harmonic functions at infinity is harmonic at infinity, $u_-$ is harmonic in $\Omega^-$ and harmonic at infinity.
 Moreover,  if $n\geq 3$, we have $c_{g_l}=0$ for all $l\in {\mathbb{N}}$ and the limit at infinity of $u_{-}$ equals $0$. If instead $n=2$, we have
\begin{equation}\label{thm:enpokerc02}
 \lim_{l\to\infty}c_{g_{l}}=\lim_{l\to\infty}\Xint{-}_{ \partial{\mathbb{B}}_{n}(0,r) }u_{g_l-}\,d\sigma
 =\Xint{-}_{ \partial{\mathbb{B}}_{n}(0,r) }u_{-}\,d\sigma=(u_-)_\infty\,.
\end{equation}
Moreover, if $x\in \Omega^-$, then we have
 \begin{eqnarray}\label{thm:enpokerc03}
\lefteqn{
u_-(x)=\lim_{l\to\infty}u_{g_{l}-}(x)
}
\\ \nonumber
&&\qquad
=\lim_{l\to\infty} -\int_{\partial\Omega}
g_{l}(y)\frac{\partial}{\partial\nu_{\Omega,y} }(S_{n} (x-y)-h_{x-}(y))\,d\sigma_{y}+\lim_{l\to\infty} c_{g_l}
\\ \nonumber
&&\qquad
=-\int_{\partial\Omega}
g(y)\frac{\partial}{\partial\nu_{\Omega,y} }(S_{n} (x-y)-h_{x-}(y))\,d\sigma_{y}+(u_-)_\infty\,.
\end{eqnarray}
In particular, the function in the right hand side of (\ref{thm:enpokerc03}) is harmonic in $\Omega^-$ and harmonic at infinity. 
We also note that the function in the right hand side of (\ref{thm:enpokerc03}) coincides with 
the function $u_{g-}$ of the statement with the choice $c_g=(u_-)_\infty$. On the other hand 
  by adding a nonzero constant to the right hand side of (\ref{thm:enpokerc03}) we would no longer obtain a solution of the exterior Dirichlet problem.  Hence, the right hand side of (\ref{thm:enpokerc03}) coincides with the function $u_{g-}$ of the statement and the value $c_g=(u_-)_\infty$ of $c_g$   in {\color{black} (\ref{thm:enpokerc01b})	}
   is uniquely determined.  Then  the existence  Theorem \ref{thm:existenceO-0} of the Appendix for the exterior Dirichlet problem implies that there exists   $u^\sharp_-\equiv {\mathcal{G}}_{d,-}[g]
  \in C^{0,\alpha}_{b}(\overline{\Omega^-})$ that  is harmonic in $\Omega^-$, harmonic at infinity and satisfies the equality $u^\sharp_-(x)=g(x)$ for all $x\in\partial\Omega$. Then the Maximum Principle for harmonic functions that are harmonic at infinity in an exterior domain implies that $u_{g-}=u^\sharp_-\in C^{0,\alpha}_{b}(\overline{\Omega^-})$. Moreover, if $x\in \Omega$, then the validity of  formula (\ref{thm:enpokerco1a}) for {\color{black} $g_l$ for 	}
  all $l\in {\mathbb{N}}$ and the uniform convergence of $\{g_l\}_{l\in {\mathbb{N}}}$ to $g$ and the limiting relation
(\ref{thm:enpokerc02}) imply the validity of formula (\ref{thm:enpokerc01a}).\hfill  $\Box$ 

\vspace{\baselineskip}

 \section{A representation theorem for the double layer potential}
  \label{sec:thm:dlintesl}
We now turn to prove a formula for the double layer potential in terms of (distributional) single layer potentials.
\begin{theorem}\label{thm:dlintesl}
 Let $\alpha\in ]0,1[$. Let $\Omega$ be a bounded open  subset of ${\mathbb{R}}^{n}$ of class $C^{1,\alpha}$.
 Let $\mu\in C^{0,\alpha}(\partial\Omega)$. Then the following statements hold.
 \begin{enumerate}
\item[(i)] Let $ {\mathcal{G}}_{d,+}[\mu] \in C^{0,\alpha}(\overline{\Omega})$ be as in
Theorem \ref{thm:idwp}.  Then the following formulas hold true
\begin{eqnarray}\label{thm:dlintesl1}
 w^+_\Omega[\mu](x)&=&{\mathcal{G}}_{d,+}[\mu](x)+v_\Omega^+[S_+^t[\mu]](x)
 \qquad\forall x\in \Omega\,,
 \\  \nonumber
w^-_\Omega[\mu](x)&=&v_\Omega^-[S_+^t[\mu]](x)\qquad\forall x\in \Omega^-\,.
\end{eqnarray}
\item[(ii)] Let ${\mathcal{G}}_{d,-}[\mu]\in C^{0,\alpha}_{ b}(\overline{\Omega^-})$ be as in Theorem \ref{thm:edwp}. Then the following formulas hold true
\begin{eqnarray}\label{thm:dlintesl2}
w^+_\Omega[\mu](x)&=&-v_\Omega^+[S_-^t[\mu]](x) +\lim_{\xi\to\infty}{\mathcal{G}}_{d,-}[\mu](\xi) \qquad\forall x\in \Omega\,,
 \\ \nonumber
w^-_\Omega[\mu](x)
&=&-{\mathcal{G}}_{d,-}[\mu](x)-
v_\Omega^-[S_-^t[\mu]](x)
+\lim_{\xi\to\infty}{\mathcal{G}}_{d,-}[\mu](\xi) \qquad\forall x\in \Omega^-\,,
\end{eqnarray}
 (cf.~(\ref{cor:existenceO-01a})).
\end{enumerate}
\end{theorem}
{\bf Proof.} (i) By Theorem \ref{thm:inpokerc0}, we can write 
the solution ${\mathcal{G}}_{d,+}[\mu]$ of the interior Dirichlet problem  with datum $\mu$ as ${\mathcal{G}}_{d,+}[\mu]=u_\mu$, where $u_\mu$ is the integral representation in terms of the Green function for the interior Dirichlet problem.  If $x\in \Omega$, then we have
\begin{eqnarray*}
\lefteqn{
w^+_\Omega[\mu](x)=\int_{\partial\Omega}\mu(y)\frac{\partial}{\partial \nu_{\Omega_y} }
\left(S_n(x-y)\right)\,d\sigma_y
}
\\ \nonumber
&&\qquad
=\int_{\partial\Omega}\mu(y)\frac{\partial}{\partial \nu_{\Omega_y} }
\left(S_n(x-y)-h_x(y)\right)\,d\sigma_y+\int_{\partial\Omega}\mu(y)\frac{\partial}{\partial \nu_{\Omega_y} }
\left(h_x(y)\right)\,d\sigma_y
\\ \nonumber
&&\qquad
=u_\mu(x)+ 
\int_{\partial\Omega}\mu(y)S_+[S_n(x-\cdot)](y)
 \,d\sigma_y
 \\ \nonumber
&&\qquad
 =u_\mu(x)+<r_{|\partial\Omega}^tS_+^t[\mu](y),S_n(x-y)>
 =
u_\mu(x)+ v_\Omega^+[S_+^t[\mu]](x)\,,
\end{eqnarray*}
(cf.~(\ref{prop:dslfun6}), (\ref{lem:ihx1})). 
If instead $x\in \Omega^-$, then $S_n(x-\cdot)$ belongs to $C^{1,\alpha}(\overline{\Omega})$, 
is harmonic in $\Omega$ and  we have
\begin{eqnarray*}
\lefteqn{
{\color{black}w^-_\Omega	}	[\mu](x)=\int_{\partial\Omega}\mu(y)\frac{\partial}{\partial \nu_{\Omega_y} }
\left(S_n(x-y)\right)\,d\sigma_y
=
\int_{\partial\Omega}\mu(y)S_+[S_n(x-\cdot)](y)
 \,d\sigma_y
}
\\ \nonumber
&&\qquad\qquad\qquad\qquad\qquad
 =<r_{|\partial\Omega}^tS_+^t[\mu](y),S_n(x-y)>=v_\Omega^-[S_+^t[\mu]](x)\,,
\end{eqnarray*}
(cf.~(\ref{prop:dslfun6})). 

(ii) By Theorem \ref{thm:enpokerc0}, we can write 
the solution ${\mathcal{G}}_{d,-}[\mu]$ of the exterior Dirichlet problem  with datum $\mu$ as
 ${\mathcal{G}}_{d,-}[\mu]=u_{\mu-}$ and we have 
$c_\mu=\lim_{\xi\to\infty}{\mathcal{G}}_{d,-}[\mu](\xi)$.  Let $x\in \Omega$. Then  $S_n(x-\cdot)$ belongs to $C^{1,\alpha}_{{\mathrm{loc}} }(\overline{\Omega^-})$, 
is harmonic in $\Omega^-$, but is harmonic at infinity only for $n\geq 3$. Thus we must consider case  $n\geq 3$ and case $n=2$ separately.  In case $n\geq 3$, we have
\begin{eqnarray*}
\lefteqn{
w^+_\Omega[\mu](x)=\int_{\partial\Omega}\mu(y)\frac{\partial}{\partial \nu_{\Omega_y} }
\left(S_n(x-y)\right)\,d\sigma_y
=-\int_{\partial\Omega}\mu(y)\frac{\partial}{\partial \nu_{\Omega^-_y} }
\left(S_n(x-y)\right)\,d\sigma_y
}
\\ \nonumber
&&\qquad\qquad
=-\int_{\partial\Omega}\mu(y)S_-[S_n(x-\cdot)](y)
 \,d\sigma_y
 =
 -<r_{|\partial\Omega}^tS_-^t[\mu](y),S_n(x-y)>
 \\ \nonumber
&&\qquad\qquad
 =-v_\Omega^+[S_-^t[\mu]](x)\,,
\end{eqnarray*}
(cf.~(\ref{prop:dslfun6})). In case $n=2$, we note that if we fix $x_0\in \Omega$, then
\[
S_2(x-y)-S_2(x_0-y)=\frac{1}{2\pi}\log \frac{|x-y|}{|x_0-y|} \qquad \forall y\in {\mathbb{R}}^n\setminus\{x,x_0\}
\]
and that $\frac{1}{2\pi}\log \frac{|x-\cdot|}{|x_0-\cdot|}$ is harmonic in $\Omega^-$ and harmonic at infinity. Then we have
\begin{eqnarray*}
\lefteqn{
w^+_\Omega[\mu](x)=\int_{\partial\Omega}\mu(y)\frac{\partial}{\partial \nu_{\Omega_y} }
\left(S_2(x-y)\right)\,d\sigma_y
=-\int_{\partial\Omega}\mu(y)\frac{\partial}{\partial \nu_{\Omega^-_y} }
\left(S_2(x-y)\right)\,d\sigma_y
}
\\ \nonumber
&&\qquad\qquad 
=-\int_{\partial\Omega}\mu(y)\frac{\partial}{\partial \nu_{\Omega^-_y} }
\left(S_2(x-y)-S_2(x_0-y)\right)\,d\sigma_y
\\ \nonumber
&&\qquad\qquad\quad
-\int_{\partial\Omega}\mu(y)\frac{\partial}{\partial \nu_{\Omega^-_y} }
\left(S_2(x_0-y)\right)\,d\sigma_y
\\ \nonumber
&&\qquad\qquad 
=-\int_{\partial\Omega}\mu(y)S_-[S_2(x-\cdot)-S_2(x_0-\cdot)](y)
 \,d\sigma_y
 \\ \nonumber
&&\qquad\qquad\quad
 -\int_{\partial\Omega}\mu(y)\frac{\partial}{\partial \nu_{\Omega^-_y} }
\left(S_2(x_0-y)\right)\,d\sigma_y
\\ \nonumber
&&\qquad\qquad 
=-\int_{\partial\Omega}\mu(y)S_-[S_2(x-\cdot)]
 \,d\sigma_y
\\ \nonumber
&&\qquad\qquad \quad
+\int_{\partial\Omega}\mu(y)S_-[h_{x_0-}](y)
 \,d\sigma_y
 -\int_{\partial\Omega}\mu(y)\frac{\partial}{\partial \nu_{\Omega^-_y} }
\left(S_2(x_0-y)\right)\,d\sigma_y
\\ \nonumber
&&\qquad\qquad 
 =
 -<r_{|\partial\Omega}^tS_-^t[\mu](y),S_2(x-\cdot)>
 +c_\mu
\\ \nonumber
&&\qquad\qquad 
 =-v_\Omega^+[S_-^t[\mu]](x)
 +c_\mu
\,,
\end{eqnarray*} (cf.~equality (\ref{thm:enpokerc01a})).
Now let $x\in\Omega^-$. Then
\begin{eqnarray*}
\lefteqn{
w^-_\Omega[\mu](x)=\int_{\partial\Omega}\mu(y)\frac{\partial}{\partial \nu_{\Omega_y} }
\left(S_n(x-y)\right)\,d\sigma_y
}
\\ \nonumber
&&\qquad
=\int_{\partial\Omega}\mu(y)\frac{\partial}{\partial \nu_{\Omega_y} }
\left(S_n(x-y)-h_{x-}(y)\right)\,d\sigma_y-c_\mu
\\ \nonumber
&&\qquad\quad
+\int_{\partial\Omega}\mu(y)\frac{\partial}{\partial \nu_{\Omega_y} }
\left(h_{x-}(y)\right)\,d\sigma_y+c_\mu
\\ \nonumber
&&\qquad
=-u_{\mu-}(x)-\int_{\partial\Omega}\mu(y)\frac{\partial}{\partial \nu_{\Omega^-_y} }
\left(h_{x-}(y)\right)\,d\sigma_y+c_\mu
\\ \nonumber
&&\qquad
=-u_{\mu-}(x)-\int_{\partial\Omega}\mu(y)S_-[h_{x-}]\,d\sigma_y+c_\mu
\\ \nonumber
&&\qquad
=-u_{\mu-}(x)-
<r_{|\partial\Omega}^tS_-^t[\mu](y),S_n(x-\cdot))>
+c_\mu
\\ \nonumber
&&\qquad
=-u_{\mu-}(x)-
v_\Omega^-[S_-^t[\mu]](x)
+c_\mu
\end{eqnarray*}
(cf.~equality (\ref{thm:enpokerc01})).\hfill  $\Box$ 

\vspace{\baselineskip}

\section{A third Green Identity for  H\"{o}l\-der continuous harmonic functions}
\label{sec:3rdGreen0a}

We now exploit the representation Theorem \ref{thm:dlintesl} in order to deduce a third Green identity for H\"{o}l\-der continuous harmonic functions in the interior and in the exterior, which extend the classical ones.
\begin{theorem}\label{thm:3rdGreen0a}
 Let $\alpha\in ]0,1[$. Let $\Omega$ be a bounded open  subset of ${\mathbb{R}}^{n}$ of class $C^{1,\alpha}$.
 \begin{enumerate}
\item[(i)] If $u\in C^{0,\alpha}(\overline{\Omega})$ is harmonic in $\Omega$, then 
\begin{eqnarray}\label{thm:3rdGreen0a1}
u(x)&=&w^+_\Omega[u_{|\partial\Omega}]
-v_\Omega^+[\frac{\partial u}{\partial\nu_{\Omega}}](x) \qquad\forall x\in\Omega\,,
\\ \nonumber
0&=&w^-_\Omega[u_{|\partial\Omega}]
-v_\Omega^-[\frac{\partial u}{\partial\nu_{\Omega}}](x) \qquad\forall x\in\Omega^-\,.
\end{eqnarray}
(where $\frac{\partial u}{\partial\nu_{\Omega}}$ is the distributional  normal derivative of $u$ in the sense of Definition \ref{defn:conoderd}).  
\item[(ii)] If $u\in C^{0,\alpha}_{ {\mathrm{loc}}}(\overline{\Omega^-})$ is harmonic in $\Omega^-$ and harmonic at infinity, then
\begin{eqnarray}\label{thm:3rdGreen0a2}
0&=&-w^+_\Omega[u_{|\partial\Omega}]
-v_\Omega^+[\frac{\partial u}{\partial\nu_{\Omega^-}}](x) +u_\infty\qquad\forall x\in\Omega\,,
\\ \nonumber
u(x)&=&-w^-_\Omega[u_{|\partial\Omega}]
-v_\Omega^-[\frac{\partial u}{\partial\nu_{\Omega^-}}](x)+u_\infty \qquad\forall x\in\Omega^-\,.
\end{eqnarray}
(where $\frac{\partial u}{\partial\nu_{\Omega^-}}$ is the distributional  normal derivative of $u$ in the sense of Definition \ref{defn:conoederd}).  
\end{enumerate}
\end{theorem}
{\bf Proof.} (i) It suffices to set $\mu=u_{|\partial\Omega}$, to apply Theorem \ref{thm:dlintesl} (i) and to observe that $u= {\mathcal{G}}_{d,+}[\mu]$ and  $S_+^t[u_{|\partial\Omega}]=\frac{\partial u}{\partial\nu_{\Omega}}$. Similarly, to  prove statement (ii) it 
suffices to set $\mu=u_{|\partial\Omega}$, to apply Theorem \ref{thm:dlintesl} (ii) and to observe that $u={\mathcal{G}}_{d,-}[\mu]$,  $S_-^t[u_{|\partial\Omega}]=\frac{\partial u}{\partial\nu_{\Omega^-}}$, $\lim_{\xi\to\infty}{\mathcal{G}}_{d,-}[\mu](\xi) =u_\infty$.\hfill  $\Box$ 

\vspace{\baselineskip}

  \section{A representation theorem for H\"{o}l\-der continuous harmonic functions in terms of single layer potentials}\label{sec:slreth}
 By the representation Theorem \ref{thm:dlintesl} for the double layer potential and by the existence theorem for the interior and exterior Dirichlet problems for harmonic functions, we can prove the following representation theorem for harmonic functions in terms of single layer potentials.
\begin{theorem}\label{thm:slreth}
 Let $\alpha\in ]0,1[$. Let $\Omega$ be a bounded open  subset of ${\mathbb{R}}^{n}$ of class $C^{1,\alpha}$.
Then the following statements hold.
\begin{enumerate}
\item[(i)] Let $n\geq 3$. If $u\in C^{0,\alpha}(\overline{\Omega})$ is harmonic in $\Omega$, then there exist
$\mu_0$ and $\mu_1\in C^{0,\alpha}(\partial\Omega)$ such that $u=v^+_\Omega[\mu_0+S_-^t[\mu_1]]$.
\item[(ii)] Let $n\geq 2$. If $u\in C^{0,\alpha}(\overline{\Omega})$ is harmonic in $\Omega$, then there exist
$\mu_0\in C^{0,\alpha}(\partial\Omega)_0$, $\mu_1\in C^{0,\alpha}(\partial\Omega)$ and $c\in{\mathbb{R}}$ such that $u=v^+_\Omega[\mu_0+S_-^t[\mu_1]]+c$.
\item[(iii)] Let $n\geq 3$. If $u\in C^{0,\alpha}_{b}(\overline{\Omega^-})$ is harmonic in $\Omega^-$ and harmonic at infinity, then there exist
$\mu_0$ and $\mu_1\in C^{0,\alpha}(\partial\Omega)$ such that $u=v^-_\Omega[\mu_0+S_+^t[\mu_1]]$.
\item[(iv)]  Let $n=2$. If $u\in C^{0,\alpha}_{b}(\overline{\Omega^-})$ is harmonic in $\Omega^-$ and harmonic at infinity, then there exist
$\mu_0\in C^{0,\alpha}(\partial\Omega)_0$, $\mu_1\in C^{0,\alpha}(\partial\Omega)$ and $c\in{\mathbb{R}}$ such that $u=v^-_\Omega[\mu_0+S_+^t[\mu_1]]+c$.
\end{enumerate}
\end{theorem}
{\bf Proof.} We first prove statements (i) and (ii). By the Existence Theorem \ref{thm:existenceO0} of the Appendix for the interior Dirichlet problem, there exist $\psi_1,\psi_2\in C^{0,\alpha}(\partial\Omega)$ such that $u=w_\Omega^+[\psi_1]+v_\Omega^+[\psi_2]$. Then the representation Theorem \ref{thm:dlintesl} (ii) for the double layer potential implies that
\[
w_\Omega^+[\psi_1]=-v_\Omega^+[S_-^t[\psi_1]] +\lim_{\xi\to\infty}{\mathcal{G}}_{d,-}[{\color{black}\psi_1	}
](\xi)\qquad\text{in}\ \Omega\,,
\]
and $c_{\psi_1}\equiv \lim_{\xi\to\infty}{\mathcal{G}}_{d,-}[
{\color{black}\psi_1	}
](\xi)=0$ if $n\geq 3$. Hence,
\[
u=-v_\Omega^+[S_-^t[\psi_1]] +c_{\psi_1}+v_\Omega^+[\psi_2]\qquad\text{in}\ \Omega\,.
\]
Since   $c_{\psi_1}=0$ if $n\geq 3$, then statement (i) holds true with
$\mu_0=\psi_2$, $\mu_1=-\psi_1$. To prove (ii), we note that Theorem \ref{slay} implies that
$c_{\psi_1}+v_\Omega^+[\psi_2]$ belongs to $C^{1,\alpha}(\overline{\Omega})$ and is harmonic in $\Omega$. Then a classical  representation Theorem implies that 
there exists $(\psi_3,c)\in C^{0,\alpha}(\partial\Omega)_0\times {\mathbb{R}}$ such that
\[
c_{\psi_1}+v_\Omega^+[\psi_2]=v_\Omega^+[\psi_3]+c\,,
\]
(cf.~\cite[Thm.~6.48 (i)]{DaLaMu21}). 
Hence, statement (ii) holds true with   $\mu_0=\psi_3$, $\mu_1=-\psi_1$.  Next we prove statements (iii), (iv).
By the Existence Theorem \ref{thm:existenceO-0} of the Appendix for the exterior Dirichlet problem, there exist $\psi_1 $,  $\psi_2\in C^{0,\alpha}(\partial\Omega)$    in case $n\geq 3$ and 
$\psi_1 \in C^{0,\alpha}(\partial\Omega)$, $\psi_2\in C^{0,\alpha}(\partial\Omega)_0$, $\rho\in {\mathbb{R}}$ in case $n=2$
 such that $u=w_\Omega^-[\psi_1]+v_\Omega^-[\psi_2]+\rho$ in $\Omega^-$ with $\rho=0$ if $n\geq 3$. Then the representation Theorem \ref{thm:dlintesl} (i) for the double layer potential implies that
\[
w_\Omega^-[\psi_1]=v_\Omega^-[S_+^t[\psi_1]]\qquad\text{in}\ \Omega^-\,.
\]
Hence,
\[
u=v_\Omega^-[S_+^t[\psi_1]]+v_\Omega^-[\psi_2]+\rho \qquad\text{in}\ \Omega^-
\]
and we can take $\mu_0=\psi_2$, $\mu_1=\psi_1$, $c=\rho$, with $c=0$ if $n\geq 3$.\hfill  $\Box$ 

\vspace{\baselineskip}

 Theorem \ref{thm:slreth} shows that all H\"{o}lder continuous  maps that  are either harmonic in  $\Omega$ or harmonic in  $\Omega^-$ and harmonic at infinity, can be written (up to an additive constant) in terms of  (distributional) single layer potentials of 
  distributions of the form  $\mu_0+S_+^t[\mu_1]$ and $\mu_0+S_-^t[\mu_1]$ for suitable choices of  $\alpha$-H\"{o}lder continuous functions  $\mu_0$, $\mu_1$ in $\partial\Omega$.  
Hence, we find convenient to introduce the following definition.
 \begin{definition}\label{defn:njsdopm}
 Let   $\alpha\in ]0,1[$. Let $\Omega$ be a bounded open  subset of ${\mathbb{R}}^{n}$ of class $C^{1,\alpha}$. Let 
 \begin{eqnarray}\label{defn:njsdopm1}
 V^{-1,\alpha,+}(\partial\Omega)&\equiv&\biggl\{\mu_0+S_+^t[\mu_1]:\,\mu_0, \mu_1\in C^{0,\alpha}(\partial\Omega)
\biggr\}\,,
\\ \nonumber
V^{-1,\alpha,-}(\partial\Omega)&\equiv&\biggl\{\mu_0+S_-^t[\mu_1]:\,\mu_0, \mu_1\in C^{0,\alpha}(\partial\Omega)
\biggr\}\,.
\end{eqnarray}
 \end{definition}
 According to the above definition, the space $V^{-1,\alpha,\pm}(\partial\Omega)$ is the image of the  map 
\[
A_{\partial\Omega\pm}:\,(C^{0,\alpha}(\partial\Omega))^{2}\to (C^{1,\alpha}(\partial\Omega))'
\]
that takes a pair $(\mu_0,\mu_1)$ to $\mu_0+S_\pm^t[\mu_1]$. By Lemma \ref{lem:caincl} and Theorems \ref{thm:contidnu}, \ref{thm:contednu}, the maps $A_{\partial\Omega\pm}$ are continuous.
Let $\pi_{\partial\Omega\pm}$ denote the canonical projection
\[
\pi_{\partial\Omega\pm}:\,(C^{0,\alpha}(\partial\Omega))^{2}\to (C^{0,\alpha}(\partial\Omega))^{2}/{\mathrm{Ker}}\, A_{\partial\Omega\pm}
\]
of $(C^{0,\alpha}(\partial\Omega))^{2}$ onto the quotient space
$(C^{0,\alpha}(\partial\Omega))^{2}/{\mathrm{Ker}}\, A_{\partial\Omega\pm}$. Let $\tilde{A}_{\partial\Omega\pm}$ be the unique linear injection from $(C^{0,\alpha}(\partial\Omega))^{2}/{\mathrm{Ker}}\, A_{\partial\Omega\pm}$ onto the image $V^{-1,\alpha,\pm}(\partial\Omega)$ of $A_{\partial\Omega\pm}$ such that 
\[
A_{\partial\Omega\pm}=\tilde{A}_{\partial\Omega\pm}\circ\pi_{\partial\Omega\pm}\,. 
\]
Then $\tilde{A}_{\partial\Omega\pm}$ is a linear bijection from $(C^{0,\alpha}(\partial\Omega))^{2}/{\mathrm{Ker}}\, A_{\partial\Omega\pm}$ onto  
 $V^{-1,\alpha,\pm}(\partial\Omega)$.\par
 
 Since the map  $A_{\partial\Omega\pm}$ is continuous,  $
{\mathrm{Ker}}\, A_{\partial\Omega\pm}$ is a closed subspace of the Banach space  $(C^{0,\alpha}(\partial\Omega))^{2}$. Hence,  $(C^{0,\alpha}(\partial\Omega))^{2}/{\mathrm{Ker}}\, A_{\partial\Omega\pm}$ is a Banach space (cf. \textit{e.g.}, \cite[Thm.~2.1]{DaLaMu21}).
 We endow 
$V^{-1,\alpha,\pm}(\partial\Omega)$ with the norm induced by  $\tilde{A}_{\partial\Omega\pm}$, i.e., we set
 \begin{equation}\label{defn:njsdopm2}
\|\tau\|_{  V^{-1,\alpha,\pm}(\partial\Omega) }
\equiv\inf\biggl\{\biggr.
 \|\mu_0\|_{ C^{0,\alpha}(\partial\Omega)  }+\|\mu_1\|_{ C^{0,\alpha}(\partial\Omega)  }
:\,
 \tau=\mu_0+S_\pm^t[\mu_1]\biggl.\biggr\}\,,
\end{equation}
for all $\tau\in  V^{-1,\alpha,\pm}(\partial\Omega)$. 
By definition of the norm $\|\cdot\|_{  V^{-1,\alpha,\pm}(\partial\Omega)  }$, the  linear bijection $\tilde{A}_{\partial\Omega\pm}$ is an isometry of the space $(C^{0,\alpha}(\partial\Omega))^{2}/{\mathrm{Ker}}\, A_{\partial\Omega\pm}$  onto the space $(V^{-1,\alpha,\pm}(\partial\Omega), \|\cdot\|_{  V^{-1,\alpha,\pm}(\partial\Omega)  })$.

Since the quotient  
$(C^{0,\alpha}(\partial\Omega))^{2}/{\mathrm{Ker}}\, A_{\partial\Omega\pm}$ is a Banach space, it follows that $(V^{-1,\alpha,\pm}(\partial\Omega), \|\cdot\|_{  V^{-1,\alpha,\pm}(\partial\Omega)  })$ is also a Banach space.\par

Since $A_{\partial\Omega\pm}$ is continuous from $(C^{0,\alpha}(\partial\Omega))^{2}$ to $(C^{1,\alpha}(\partial\Omega))'$, 
 a fundamental property of the quotient topology implies that the map 
$\tilde{A}_{\partial\Omega\pm}$ is  continuous from the quotient space $(C^{0,\alpha}(\partial\Omega))^{2}/{\mathrm{Ker}}\, A_{\partial\Omega\pm}$ to $(C^{1,\alpha}(\partial\Omega))'$   (cf. \textit{e.g.}, \cite[Prop.~A.5]{DaLaMu21}). 

Hence, the Banach space $(V^{-1,\alpha,\pm}(\partial\Omega), \|\cdot\|_{V^{-1,\alpha,\pm}(\partial\Omega)  })$ is continuously embedded into $(C^{1,\alpha}(\partial\Omega))'$.
Also, the definition of the norm $\|\cdot\|_{V^{-1,\alpha,\pm}(\partial\Omega)  }$ implies that $C^{0,\alpha}(\partial\Omega)$ is continuously embedded into $V^{-1,\alpha,\pm}(\partial\Omega)$ and that $S_\pm^t$ is continuous from 
$C^{0,\alpha}(\partial\Omega)$ to $V^{-1,\alpha,\pm}(\partial\Omega)$. The elements of $V^{-1,\alpha,\pm}(\partial\Omega)$ are not integrable functions, but distributions in $\partial\Omega$. We also point out the validity of the following elementary but useful lemma.
 
\begin{lemma}\label{lem:cov-1apm}
 Let  $\alpha\in]0,1
 {\color{black}[	}	
 $. Let $\Omega$ be a bounded open  subset of ${\mathbb{R}}^{n}$ of class $C^{1
 {\color{black},\alpha	}	
 }$. Let $X$ be a normed space. Let $L$ be a linear map from $V^{-1,\alpha,\pm}(\partial\Omega)$ to $X$. Then $L$ is continuous if and only if the map
 \[
 L \circ  A_{\partial\Omega\pm}
 \]
 is continuous on $C^{0,\alpha}(\partial\Omega)^{2}$.
\end{lemma}
 {\bf Proof.} If $L$ is continuous, then so is the composite map $L \circ  A_{\partial\Omega\pm}$. Conversely, if $L \circ  A_{\partial\Omega\pm}$ 
 {\color{black}is continuous	}	
 we note that
\[
L \circ A_{\partial\Omega}=L \circ  \tilde{A}_{\partial\Omega\pm}\circ\pi_{\partial\Omega\pm}\,.
\]
Then a  fundamental property of the quotient topology implies that the map $L \circ  \tilde{A}_{\partial\Omega\pm}$ is continuous on the quotient $(C^{0,\alpha}(\partial\Omega))^{2}/{\mathrm{Ker}}\, A_{\partial\Omega\pm}$. Since $\tilde{A}_{\partial\Omega}$  is an isometry from
$(C^{0,\alpha}(\partial\Omega))^{2}/{\mathrm{Ker}}\, A_{\partial\Omega\pm}$ onto $V^{-1,\alpha,\pm}(\partial\Omega)$, its inverse map is continuous and accordingly
\[
L=L \circ  \tilde{A}_{\partial\Omega\pm}\circ \left(\tilde{A}_{\partial\Omega\pm}\right)^{(-1)}
\]
is continuous.\hfill  $\Box$ 

\vspace{\baselineskip}
 
 Finally, we point out the validity of the following lemma.
 
\begin{lemma}\label{lem:v-1at1}
 Let   $\alpha\in ]0,1[$. Let $\Omega$ be a bounded open  subset of ${\mathbb{R}}^{n}$ of class $C^{1,\alpha}$. Let $\mu_0$, $\mu_1\in C^{0,\alpha}(\partial\Omega)$. Then
\begin{eqnarray}\label{lem:v-1at11}
 <S_+^t[\mu_1],1>&=&0\qquad \text{if}\  n\geq 2\,,
 \\ \nonumber
 <S_-^t[\mu_1],1>&=&0\qquad \text{if}\  n=2\,.
 \end{eqnarray}
 In particular, if $\tau_\pm\in  V^{-1,\alpha,\pm}(\partial\Omega)$ and if $\tau_\pm=\mu_0+S_\pm^t[\mu_1]$, then 
\begin{eqnarray}\label{lem:v-1at12}
 <\tau_+,1>&=&<\mu_0,1>  \qquad \text{if}\  n\geq 2\,,
  \\ \nonumber
 <\tau_-,1>&=&<\mu_0,1> \qquad \text{if}\  n=2\,.
 \end{eqnarray}
\end{lemma}
{\bf Proof.} It suffices to note that $S_+[1]=0$ for $n\geq 2$, that $S_-[1]=0$ for $n=2$.\hfill  $\Box$ 

\vspace{\baselineskip}

\section{Trace formulas for the single layer potential with distributional densities}
\label{sec:trdslay}
We plan to prove some properties of the single layer potential by exploiting the representation formulas of
Theorem \ref{thm:dlintesl} that involve  the double layer potential. 
\begin{theorem}\label{thm:vpm}
 Let $\alpha\in ]0,1[$. Let $\Omega$ be a bounded open  subset of ${\mathbb{R}}^{n}$ of class $C^{1,\alpha}$. Then the following statements hold.
\begin{enumerate}
\item[(i)] If $\mu\in C^{0,\alpha}(\partial\Omega)$, then $v_\Omega^+[S_+^t[\mu]]$ has an unique extension to $\overline{\Omega}$ that belongs to $C^{0,\alpha}(\overline{\Omega})$ (and that we still denote by the same symbol).
Moreover, the map from $C^{0,\alpha}(\partial\Omega)$ to $C^{0,\alpha}(\overline{\Omega})$ that takes 
$\mu$ to $v_\Omega^+[S_+^t[\mu]]$ is linear and continuous.
\item[(ii)] If $\mu\in C^{0,\alpha}(\partial\Omega)$, then $v_\Omega^-[S_+^t[\mu]]$ has an unique extension to $\overline{\Omega^-}$ that belongs to $C^{0,\alpha}_b(\overline{\Omega^-})$ (and that we still denote by the same symbol).
Moreover, the map from $C^{0,\alpha}(\partial\Omega)$ to $C^{0,\alpha}_b(\overline{\Omega^-})$ that takes 
$\mu$ to $v_\Omega^-[S_+^t[\mu]]$ is linear and continuous.
\item[(iii)] If $\mu\in C^{0,\alpha}(\partial\Omega)$, then 
\begin{equation}\label{thm:vpm1}
v_\Omega^+[S_+^t[\mu]](x)=v_\Omega^-[S_+^t[\mu]](x)\qquad\forall x\in\partial\Omega\,.
\end{equation}
\item[(iv)] If $\mu\in C^{0,\alpha}(\partial\Omega)$, then $v_\Omega^+[S_-^t[\mu]]$ has an unique extension to $\overline{\Omega}$ that belongs to $C^{0,\alpha}(\overline{\Omega})$ (and that we still denote by the same symbol).
Moreover, the map from $C^{0,\alpha}(\partial\Omega)$ to $C^{0,\alpha}(\overline{\Omega})$ that takes 
$\mu$ to $v_\Omega^+[S_-^t[\mu]]$ is linear and continuous.
\item[(v)] If $\mu\in C^{0,\alpha}(\partial\Omega)$, then $v_\Omega^-[S_-^t[\mu]]$ has an unique extension to $\overline{\Omega^-}$ that belongs to $C^{0,\alpha}_{b}(\overline{\Omega^-})$ (and that we still denote by the same symbol). 
Moreover, the map from $C^{0,\alpha}(\partial\Omega)$ to $C^{0,\alpha}_b(\overline{\Omega^-})$ that takes 
$\mu$ to $v_\Omega^-[S_-^t[\mu]]$ is linear and continuous.
\item[(vi)] If $\mu\in C^{0,\alpha}(\partial\Omega)$, then
\begin{equation}\label{thm:vpm2}
v_\Omega^+[S_-^t[\mu]](x)=v_\Omega^-[S_-^t[\mu]](x)\qquad\forall x\in\partial\Omega\,.
\end{equation}
\end{enumerate}
\end{theorem}
{\bf Proof.} Statements (i), (ii) follow  by Theorem \ref{thm:dlintesl} (i), by the  equalities of (\ref{thm:dlintesl1}), by Theorem \ref{thm:wpsi+-} on the double layer potential and   by Theorem \ref{thm:idwp}   on ${\mathcal{G}}_{d,+}[\cdot]$.
   Statement (iii) follows by 
 the  equalities of (\ref{thm:dlintesl1}), by the jump relations for the double layer potential on $\partial\Omega$ (cf.~Thm.~\ref{thm:wpsi+-}) and by equality ${\mathcal{G}}_{d,+}[\mu]_{|\partial\Omega}=\mu$.

Statements (iv), (v) follow  by Theorem \ref{thm:dlintesl} (ii), by the  equalities of (\ref{thm:dlintesl2}), by Theorem \ref{thm:wpsi+-} on the double layer potential and   by Theorem  \ref{thm:edwp} on ${\mathcal{G}}_{d,-}[\cdot]$.  

Statement (vi) follows by the  equalities of (\ref{thm:dlintesl2}), by the jump relations for the double layer potential on $\partial\Omega$ (cf.~Thm.~\ref{thm:wpsi+-}) and by equality ${\mathcal{G}}_{d,-}[\mu]_{|\partial\Omega}=\mu$.\hfill  $\Box$ 

\vspace{\baselineskip}

We note that according to the Definition \ref{defn:Vpm} of $V_\Omega$, we have 
\begin{eqnarray}\label{eq:Vpm}
V_\Omega[S_+^t[\mu]](x)&\equiv& v_\Omega^+[S_+^t[\mu]](x)=v_\Omega^-[S_+^t[\mu]](x)\quad\forall x\in\partial\Omega\,,
\\ \nonumber
V_\Omega[S_-^t[\mu]](x)&\equiv& v_\Omega^+[S_-^t[\mu]](x)=v_\Omega^-[S_-^t[\mu]](x)\quad\forall x\in\partial\Omega\,,
\end{eqnarray}
for all $\mu\in C^{0,\alpha}(\partial\Omega)$. 

By  Proposition \ref{prop:sldinfty} on the single layer potential, by the  classical regularity Theorems \ref{slay}, \ref{thm:slayh} on the single layer potential, by the continuity Lemma \ref{lem:cov-1apm} and by  Lemma \ref{lem:v-1at1},   we deduce the validity of the following immediate corollary of  Theorem \ref{thm:vpm}.
\begin{corollary}\label{cor:vpm}
 Let $\alpha\in ]0,1[$. Let $\Omega$ be a bounded open  subset of ${\mathbb{R}}^{n}$ of class $C^{1,\alpha}$. 
 \begin{enumerate}
\item[(i)] The map from $V^{-1,\alpha,+}(\partial\Omega)$ to $C^{0,\alpha}(\overline{\Omega})$ that takes 
$\tau$ to $v_\Omega^+[\tau]$ is linear and continuous.
\item[(ii)] If $n\geq 3$, then the map from $V^{-1,\alpha,+}(\partial\Omega)$ to $C^{0,\alpha}_b(\overline{\Omega^-})$ that takes 
$\tau$ to $v_\Omega^-[\tau]$ is linear and continuous.
If $n=2$, then the map from $V^{-1,\alpha,+}(\partial\Omega)$ to $C^{0,\alpha} (\overline{{\mathbb{B}}_n(0,r)}\setminus\Omega^-)$ that takes 
$\tau$ to $v_\Omega^-[\tau]_{|\overline{{\mathbb{B}}_n(0,r)}\setminus\Omega^-}$ is linear and continuous for all $r\in]0,+\infty[$ such that $\overline{\Omega^-}\subseteq {\mathbb{B}}_n(0,r)$ and the map from 
\[
V^{-1,\alpha,+}(\partial\Omega)_0\equiv
\{\tau\in 
V^{-1,\alpha,+}(\partial\Omega):\,<\tau,1>=0\}
\]
 to $C^{0,\alpha}_b(\overline{\Omega^-})$ that takes 
$\tau$ to $v_\Omega^-[\tau]$ is linear and continuous.
\item[(iii)] The 
 map  $V_\Omega[\cdot]$   from $V^{-1,\alpha,+}(\partial\Omega)$ to $C^{0,\alpha}(\partial\Omega)$ is linear and continuous.
\item[(iv)] The map from $V^{-1,\alpha,-}(\partial\Omega)$ to $C^{0,\alpha}(\overline{\Omega})$ that takes 
$\tau$ to $v_\Omega^+[\tau]$ is linear and continuous.
\item[(v)]  If $n\geq 3$, then the map from $V^{-1,\alpha,-}(\partial\Omega)$ to $C^{0,\alpha}_b(\overline{\Omega^-})$ that takes 
$\tau$ to $v_\Omega^-[\tau]$ is linear and continuous.
If $n=2$, then the map from $V^{-1,\alpha,-}(\partial\Omega)$ to $C^{0,\alpha} (\overline{{\mathbb{B}}_n(0,r)}\setminus\Omega^-)$ that takes 
$\tau$ to $v_\Omega^-[\tau]_{|\overline{{\mathbb{B}}_n(0,r)}\setminus\Omega^-}$ is linear and continuous for all $r\in]0,+\infty[$ such that $\overline{\Omega^-}\subseteq {\mathbb{B}}_n(0,r)$ and the map from 
\[
V^{-1,\alpha,-}(\partial\Omega)_0\equiv
\{\tau\in 
V^{-1,\alpha,-}(\partial\Omega):\,<\tau,1>=0\}
\]
 to $C^{0,\alpha}_b(\overline{\Omega^-})$ that takes 
$\tau$ to $v_\Omega^-[\tau]$ is linear and continuous.
\item[(vi)] The 
 map  $V_\Omega[\cdot]$   from $V^{-1,\alpha,-}(\partial\Omega)$ to $C^{0,\alpha}(\partial\Omega)$ is linear and continuous.
\end{enumerate}
\end{corollary}
 Finally, we show that the boundary operator $V_\Omega$ satisfies a symmetry property on the distributions of the form $S_\pm^t[\mu]$ when $\mu\in C^{0,\alpha}(\partial\Omega)$.
\begin{proposition}\label{prop:syvspm}
 Let $\alpha\in ]0,1[$. Let $\Omega$ be a bounded open  subset of ${\mathbb{R}}^{n}$ of class $C^{1,\alpha}$. 
Let $\mu\in C^{0,\alpha}(\partial\Omega)$. Then the following symmetry conditions are satisfied
 \begin{eqnarray}\label{prop:syvspm1}
<S_+^t[\mu],V_\Omega[\psi]>&=&\int_{\partial\Omega}V_\Omega[S_+^t[\mu]]\psi\,d\sigma\,,
\\ \nonumber
<S_-^t[\mu],V_\Omega[\psi]>&=&\int_{\partial\Omega}V_\Omega[S_-^t[\mu]]\psi\,d\sigma\,,
\qquad\forall\psi\in C^{0,\alpha}(\partial\Omega)\,.
\end{eqnarray}
\end{proposition}
{\bf Proof.} {\color{black} Let $\psi\in C^{0,\alpha}(\partial\Omega)$. We first prove the first equality in (\ref{prop:syvspm1}) in case $\mu\in  C^{1,\alpha}(\partial\Omega)$. By the definition of $S_+$ and by the classical second Green Identity (cf.~\textit{e.g.}, \cite[Thm.~4.3]{DaLaMu21}), we have
\begin{eqnarray}\label{prop:syvspm2}
\lefteqn{
<S_+^t[\mu],v>=\int_{\partial\Omega}\mu S_+[v]\,d\sigma
=\int_{\partial\Omega} {\mathcal{G}}_{d,+}[\mu] 
\frac{\partial}{\partial\nu_\Omega}{\mathcal{G}}_{d,+}[v]\,d\sigma
}
\\ \nonumber
&& 
=\int_{\partial\Omega} \frac{\partial}{\partial\nu_\Omega}{\mathcal{G}}_{d,+}[\mu] {\mathcal{G}}_{d,+}[v]\,d\sigma
=\int_{\partial\Omega} \frac{\partial}{\partial\nu_\Omega}{\mathcal{G}}_{d,+}[\mu] v\,d\sigma
\qquad\forall v\in C^{1,\alpha}(\partial\Omega)\,,
\end{eqnarray}
\textit{i.e.}, $S_+^t[\mu]$ is the distribution that is associated to $\frac{\partial}{\partial\nu_\Omega}{\mathcal{G}}_{d,+}[\mu] $. Then the Fubini Theorem implies that
\begin{eqnarray*}
\lefteqn{
<S_+^t[\mu],V_\Omega[\psi]>=\int_{\partial\Omega} \frac{\partial}{\partial\nu_\Omega}{\mathcal{G}}_{d,+}[\mu](x) \int_{\partial\Omega}S_n(x-y)\psi(y)\,d\sigma_y\,d\sigma_x
}
\\ \nonumber
&&\qquad
=\int_{\partial\Omega}\psi(y)\int_{\partial\Omega}S_n(x-y)\frac{\partial}{\partial\nu_\Omega}{\mathcal{G}}_{d,+}[\mu](x) \,d\sigma_x\,d\sigma_y
\\ \nonumber
&&\qquad
=<V_\Omega\left[\frac{\partial}{\partial\nu_\Omega}{\mathcal{G}}_{d,+}[\mu]\right],\psi>
=<V_\Omega\left[S_+^t[\mu]\right],\psi>
\end{eqnarray*}
and accordingly the first equality of (\ref{prop:syvspm1}) holds true. By precisely the same argument one can prove the second equality of (\ref{prop:syvspm1}) in case  $\mu\in  C^{1,\alpha}(\partial\Omega)$ by exploiting the classical second Green Identity in exterior domains (cf.~\textit{e.g.}, \cite[Cor.~4.8]{DaLaMu21}).

We now consider case  $\mu\in  C^{0,\alpha}(\partial\Omega)$. Let $\beta\in]0,\alpha[$. By the approximation Lemma \ref{lem:0a1aap} of the Appendix, there exists a sequence $\{\mu_j\}_{j\in{\mathbb{N}}}$ in $C^{1,\alpha}(\partial\Omega)$ that is bounded in the $C^{0,\alpha}(\partial\Omega)$-norm and that converges to $\mu$ in the $C^{0,\beta}(\partial\Omega)$-norm.

By the continuity of $S_+^t$ from $C^{0,\beta}(\partial\Omega)$ to $V^{-1,\beta,+}(\partial\Omega)$, that is continuously embedded into $\left(C^{1,\beta}(\partial\Omega)\right)'$, by the membership of $V_\Omega[\psi]$ in $C^{1,\beta}(\partial\Omega)$ that follows by Theorem \ref{slay}, by the membership of $\mu_j$ in $C^{1,\alpha}(\partial\Omega)$, by the validity of the first equality of (\ref{prop:syvspm1}) for the functions of $C^{1,\alpha}(\partial\Omega)$ and by the continuity of $V_\Omega[S_+^t[\cdot]]$ from $C^{0,\beta}(\partial\Omega)$ to
$C^{0,\beta}(\partial\Omega)$ that follows by Theorem \ref{thm:vpm}, 
we have
\begin{eqnarray*}
\lefteqn{
<S_+^t[\mu],V_\Omega[\psi]>=\lim_{j\to\infty}<S_+^t[\mu_j],V_\Omega[\psi]>
}
\\ \nonumber
&&\qquad
\lim_{j\to\infty}<V_\Omega[S_+^t[\mu_j]], \psi >=<V_\Omega[S_+^t[\mu]], \psi >
\end{eqnarray*}
and accordingly the first equality of (\ref{prop:syvspm1}) holds true for $\mu\in  C^{0,\alpha}(\partial\Omega)$. By precisely the same argument one can prove the second equality of (\ref{prop:syvspm1}) in case  $\mu\in  C^{0,\alpha}(\partial\Omega)$.
}\hfill  $\Box$ 

\vspace{\baselineskip}

\section{A distributional form of the jump formula for the normal derivative of the harmonic single layer potential and  a class of distributions on the boundary}
\label{sec:ndslv}

We first introduce  the following  classical lemma of Plemelj~\cite{Pl11}. For the sake of completeness we include a proof as in  Khavinson, Putinar and Shapiro~\cite[Lem.~2]{KhPuSh07} in a $L^p$ setting and in  Mitrea and Taylor \cite[(7.41)]{MitTa99} in a Sobolev space setting.
\begin{lemma}[Plemelj’s symmetrization principle]\label{lem:comvw}
Let $\alpha\in ]0,1[$. Let $\Omega$ be a bounded open  subset of ${\mathbb{R}}^{n}$ of class $C^{1,\alpha}$. If 
$\eta\in C^{0,\alpha}(\partial\Omega)$, then
\[
V_\Omega[W_\Omega^t[\eta]]=W_\Omega[V_\Omega[\eta]]\,.
\]
\end{lemma}
{\bf Proof.} Let $u\equiv v_\Omega^+[\eta]$. Then the third Green Identity implies that
\[
w_\Omega^-[u_{|\partial\Omega}](x)=v_\Omega^-[\frac{\partial u}{\partial\nu}](x)\qquad\forall x\in\Omega^-\,,
\]
(cf.~\textit{e.g.}, \cite[Thm.~Cor.~4.6]{DaLaMu21}). 
Then the jump formulas for the harmonic single and double layer potential imply that
\[
-\frac{1}{2}V_\Omega[\eta](x)+W_\Omega[V_\Omega[\eta]](x)=V_\Omega[-\frac{1}{2}\eta+W_\Omega^t[\eta]](x)\qquad\forall x\in\partial\Omega\,,
\]
and accordingly the formula of the statement holds true.
\hfill  $\Box$ 

\vspace{\baselineskip}

Next we prove that the validity of the  no-jump  condition  (\ref{defn:Vpm0}) and of a symmetry property for the distributional  harmonic  single layer potential on the boundary imply the validity of the jump formulas for the distributional  harmonic  single layer potential.
\begin{proposition}\label{prop:ndslv0a}
Let $\alpha\in ]0,1[$. Let $\Omega$ be a bounded open  subset of ${\mathbb{R}}^{n}$ of class $C^{1,\alpha}$. 
Let $\tau\in (C^{1,\alpha}(\partial\Omega))'$ satisfy the following assumptions
\begin{eqnarray}\label{prop:ndslv0a0}
v_\Omega^+[\tau]\in C^{0,\alpha}(\overline{\Omega}),\ \ 
v_\Omega^-[\tau]\in C^{0,\alpha}_{{\mathrm{loc}}}(\overline{\Omega^-}),\ \ 
v_\Omega^+[\tau]_{|\partial\Omega}=v_\Omega^-[\tau]_{|\partial\Omega}
\end{eqnarray}
and the symmetry condition
\begin{equation}\label{prop:ndslv0a1}
<\tau,V_\Omega[\psi]>=\int_{\partial\Omega}V_\Omega[\tau]\psi\,d\sigma\qquad\forall\psi\in C^{0,\alpha}(\partial\Omega) 
\end{equation}
(cf.~(\ref{defn:Vpm1}) for the definition of $V_\Omega$).
Then
\begin{equation}\label{prop:ndslv0a2}
\frac{\partial}{\partial\nu}v_\Omega^+[\tau]=-\frac{1}{2}\tau+W^t_\Omega[\tau]\,.
\end{equation}
If $n\geq 3$, then we also have
\begin{equation}\label{prop:ndslv0a3}
\frac{\partial}{\partial\nu_{\Omega^-}}v_\Omega^-[\tau]=-\frac{1}{2}\tau-W^t_\Omega[\tau]\,.
\end{equation}
If we also assume that $<\tau,1>=0$, then (\ref{prop:ndslv0a3}) holds also in case 
$n=2$.
\end{proposition}
{\bf Proof.} Let $\psi\in C^{1,\alpha}(\partial\Omega)$.  Let ${\mathcal{G}}_{d,+}$ has be as in Theorem \ref{thm:idwp}.
By a known representation theorem for harmonic functions there exists $(\eta,c)\in C^{0,\alpha}(\partial\Omega)_0\times{\mathbb{R}}$ such that
\begin{equation}\label{prop:ndslv0a4}
{\mathcal{G}}_{d,+}[\psi]=v_\Omega^+[\eta]+c\,, 
\end{equation}
(cf.~\textit{e.g.}, \cite[Thm.~6.48]{DaLaMu21}). Then the
 definition of distributional  normal derivative, the symmetry condition (\ref{prop:ndslv0a1}), the jump formula for the normal derivative of the single layer potential of Theorem \ref{slay} (ii), the classical  Plemelj’s symmetrization principle of Lemma \ref{lem:comvw} and Lemma \ref{lem:Sphi0a} imply that
 \begin{eqnarray*}
\lefteqn{
<\frac{\partial}{\partial\nu}v_\Omega^+[\tau],\psi> 
=\int_{\partial\Omega}V_\Omega[\tau]\frac{\partial}{\partial\nu}{\mathcal{G}}_{d,+}[\psi]\,d\sigma
}
 \\ \nonumber
&&\qquad
=\int_{\partial\Omega}V_\Omega[\tau]\frac{\partial}{\partial\nu}\left(v_\Omega^+[\eta]+c\right)\,d\sigma
\\ \nonumber
&&\qquad
=\int_{\partial\Omega}V_\Omega[\tau]\left(-\frac{1}{2}\eta+W_\Omega^t[\eta]\right)\,d\sigma
\\ \nonumber
&&\qquad
=-\frac{1}{2}\int_{\partial\Omega}V_\Omega[\tau]\eta\,d\sigma
+\int_{\partial\Omega}V_\Omega[\tau]W_\Omega^t[\eta]\,d\sigma
\\ \nonumber
&&\qquad
=-\frac{1}{2}<\tau,V_\Omega[\eta]>
+<V_\Omega[\tau],W_\Omega^t[\eta]>
\\ \nonumber
&&\qquad
=-\frac{1}{2}<\tau,V_\Omega[\eta]>
+<\tau,V_\Omega[W_\Omega^t[\eta]]>
\\ \nonumber
&&\qquad
=-\frac{1}{2}<\tau,V_\Omega[\eta]>
+<\tau,W_\Omega[V_\Omega[\eta]]>
\\ \nonumber
&&\qquad
=-\frac{1}{2}<\tau,V_\Omega[\eta]>
+<W_\Omega^t[\tau], V_\Omega[\eta]>
\\ \nonumber
&&\qquad
=<-\frac{1}{2}\tau+W_\Omega^t[\tau],V_\Omega[\eta]+c>
=<-\frac{1}{2}\tau+W_\Omega^t[\tau],\psi> 
\end{eqnarray*}
 and thus the proof of the jump formula (\ref{prop:ndslv0a2}) is complete. We now turn to prove the jump formula 
 (\ref{prop:ndslv0a3}). Let ${\mathcal{G}}_{d,-}$ has be as in Theorem \ref{thm:edwp}.
By a known representation theorem for harmonic functions that are harmonic at infinity, if $n\geq 3$, there exists $\eta\in C^{0,\alpha}(\partial\Omega)$ such that
\begin{equation}\label{prop:ndslv0a5}
{\mathcal{G}}_{d,-}[\psi]=v_\Omega^-[\eta]\qquad\text{in}\ \overline{\Omega^-}
\end{equation}
and if $n=2$ 
there exists $(\eta,c)\in C^{0,\alpha}(\partial\Omega)_0\times{\mathbb{R}}$ such that
\begin{equation}\label{prop:ndslv0a5a}
{\mathcal{G}}_{d,-}[\psi]=v_\Omega^-[\eta]+c\qquad\text{in}\ \overline{\Omega^-}
\end{equation}
and $v_\Omega^-[\eta]$ is harmonic at infinity (cf.~\textit{e.g.}, \cite[Thm.~6.48]{DaLaMu21}).  If $n\geq 3$, Proposition \ref{prop:sldinfty} implies that $v_\Omega^-[\tau]$ is harmonic at infinity. If $n=2$ our assumption
  $<\tau,1>=0$ and Proposition \ref{prop:sldinfty} imply that $v_\Omega^-[\tau]$ is harmonic at infinity. Then the
 definition of distributional  normal derivative, the symmetry condition (\ref{prop:ndslv0a1}), the classical jump formula for the normal derivative of the single layer potential of Theorem \ref{slay} (ii)  the classical Plemelj’s symmetrization principle of Lemma \ref{lem:comvw} and Lemma \ref{lem:Sphi0a} imply that
\begin{eqnarray*}
\lefteqn{
<\frac{\partial}{\partial\nu_{\Omega^-}}v_\Omega^-[\tau],\psi>=
\int_{\partial\Omega}V_\Omega[\tau]\frac{\partial}{\partial\nu_{\Omega^-}}{\mathcal{G}}_{d,-}[\psi]\,d\sigma
}
\\ \nonumber
&&\qquad
=\int_{\partial\Omega}V_\Omega[\tau]\frac{\partial}{\partial\nu_{\Omega^-}}\left(v_\Omega^-[\eta]+c\delta_{2,n}\right)\,d\sigma
 \\ \nonumber
&&\qquad
=\int_{\partial\Omega}V_\Omega[\tau]\left(-\frac{1}{2}\eta-W_\Omega^t[\eta]\right)\,d\sigma
\\ \nonumber
&&\qquad
=-\frac{1}{2}\int_{\partial\Omega}V_\Omega[\tau]\eta\,d\sigma
-\int_{\partial\Omega}V_\Omega[\tau]W_\Omega^t[\eta]\,d\sigma
\\ \nonumber
&&\qquad
=-\frac{1}{2}<\tau,V_\Omega[\eta]>
-<V_\Omega[\tau],W_\Omega^t[\eta]>
\\ \nonumber
&&\qquad
=-\frac{1}{2}<\tau,V_\Omega[\eta]>
-<\tau,V_\Omega[W_\Omega^t[\eta]]>
\\ \nonumber
&&\qquad
=-\frac{1}{2}<\tau,V_\Omega[\eta]>
-<\tau,W_\Omega[V_\Omega[\eta]]>
\\ \nonumber
&&\qquad
=-\frac{1}{2}<\tau,V_\Omega[\eta]>
-<W_\Omega^t[\tau], V_\Omega[\eta]>
\\ \nonumber
&&\qquad
=<-\frac{1}{2}\tau-W_\Omega^t[\tau],\psi-c\delta_{2,n}>
\\ \nonumber
&&\qquad
=<-\frac{1}{2}\tau-W_\Omega^t[\tau],\psi> -<-\frac{1}{2}\tau-W_\Omega^t[\tau],c\delta_{2,n}>
\\ \nonumber
&&\qquad
=<-\frac{1}{2}\tau-W_\Omega^t[\tau],\psi>+<\tau,\frac{1}{2}c\delta_{2,n}+W_\Omega[c]\delta_{2,n}>
\\ \nonumber
&&\qquad
=<-\frac{1}{2}\tau-W_\Omega^t[\tau],\psi>+<\tau,1>c\delta_{2,n}=<-\frac{1}{2}\tau-W_\Omega^t[\tau],\psi>
\end{eqnarray*}
 and thus the proof of the jump formula (\ref{prop:ndslv0a4}) is complete.
\hfill  $\Box$ 

\vspace{\baselineskip}

Proposition \ref{prop:ndslv0a} justifies the introduction the following class of distributions on $\partial\Omega$.
\begin{definition}\label{defn:njsdoa}
 Let $\alpha\in ]0,1[$. Let $\Omega$ be a bounded open  subset of ${\mathbb{R}}^{n}$ of class $C^{1,\alpha}$. Let 
 \begin{equation}\label{defn:njsdoa1}
 V^{-1,\alpha}(\partial\Omega)=\biggl\{
\tau\in (C^{1,\alpha}(\partial\Omega))':\tau \ \text{satisfies\ conditions}\ (\ref{prop:ndslv0a0}), (\ref{prop:ndslv0a1})
\biggr\}\,.
\end{equation}
\end{definition}
\begin{proposition}\label{prop:v-1apin}
Let $\alpha\in ]0,1[$. Let $\Omega$ be a bounded open  subset of ${\mathbb{R}}^{n}$ of class $C^{1,\alpha}$. Then   $C^{0,\alpha}(\partial\Omega)\subseteq V^{-1,\alpha}(\partial\Omega)$ and 
$V^{-1,\alpha,\pm}(\partial\Omega)\subseteq V^{-1,\alpha}(\partial\Omega)$.\end{proposition}
{\bf Proof.}  The first inclusion follows by the classical Theorem \ref{slay} (ii) and the Fubini Theorem. The second inclusions follow 
 by Theorem \ref{thm:vpm}, by Proposition \ref{prop:syvspm}
and by the inclusion of $C^{0,\alpha}(\partial\Omega)$ in $V^{-1,\alpha}(\partial\Omega)$.\hfill  $\Box$ 

\vspace{\baselineskip}

\section{A distributional form of the Plemelj’s symmetrization principle}
\label{sec:dcomvw}

We now exploit the distributional form of the  Green identities of Theorem \ref{thm:3rdGreen0a} and the distributional jump formulas to
show that the Plemelj's symmetrization principle has the following extension for densities in the class
$V^{-1,\alpha}(\partial\Omega)$.
\begin{lemma}[Plemelj’s symmetrization principle]\label{lem:dcomvw}
Let $\alpha\in ]0,1[$. Let $\Omega$ be a bounded open  subset of ${\mathbb{R}}^{n}$ of class $C^{1,\alpha}$. If 
$\tau\in V^{-1,\alpha}(\partial\Omega)$, then 
$W_\Omega^t[\tau]$ belongs to $V^{-1,\alpha}(\partial\Omega)$. Moreover, 
\begin{equation}\label{lem:dcomvw0}
V_\Omega[W_\Omega^t[\tau]]_{|\partial\Omega}=W_\Omega[V_\Omega[\tau]]\quad\text{on}\ \partial\Omega\,.
\end{equation}  
 \end{lemma}
{\bf Proof.} We split our proof into three parts. In the first part we show that  $v^{-}_\Omega[W_\Omega^t[\tau]]$ admits a continuous extension to $\overline{\Omega^-}$ and that 
 $v^{-}_\Omega[W_\Omega^t[\tau]]_{|\partial\Omega}=W_\Omega[V_\Omega[\tau]]$ on $\partial\Omega$. In the second part we show a corresponding statement for $v^{+}_\Omega[W_\Omega^t[\tau]]_{|\partial\Omega}$ and in the  third part we deduce the validity of our statement.
 
 For the first part of the proof, we set $u\equiv v_\Omega^+[\tau]$. By the Definition \ref{defn:njsdoa} of 
$V^{-1,\alpha}(\partial\Omega)$, we know that $u\in C^{0,\alpha}(\overline{\Omega})$. Then the third Green Identity in the distributional  form of Theorem \ref{thm:3rdGreen0a} (i) implies that
\begin{equation}\label{lem:dcomvw1}
w_\Omega^-[u_{|\partial\Omega}](x)=v_\Omega^-[\frac{\partial u}{\partial\nu}](x)\qquad\forall x\in\Omega^-\,.
\end{equation}
By the classical jump formula for the double layer potential, we know that $w_\Omega^-[u_{|\partial\Omega}]$ has a (unique) extension to $\overline{\Omega^-}$, which belongs to $C^{0,\alpha}_b(\overline{\Omega^-})$ and that
\begin{equation}\label{lem:dcomvw2}
w_\Omega^-[u_{|\partial\Omega}](x)=-\frac{1}{2} u(x)+W_\Omega[u_{|\partial\Omega}](x)\qquad\forall x\in\partial\Omega\,,
\end{equation}
(cf.~Theorem \ref{thm:wpsi+-}). Then equality (\ref{lem:dcomvw1}) implies that $v_\Omega^-[\frac{\partial u}{\partial\nu}]$ admits a unique continuous extension to $\overline{\Omega^-}$, which is of class $C^{0,\alpha}_b(\overline{\Omega^-})$. By the jump formula for the single layer potential with moments in $V^{-1,\alpha}(\partial\Omega)$ of Proposition \ref{prop:ndslv0a}, we have 
\[
\frac{\partial}{\partial\nu}u=-\frac{1}{2}\tau+W^t_\Omega[\tau] \qquad\text{on}\ \partial\Omega\,.
\] 
Then equalities (\ref{lem:dcomvw1}) and (\ref{lem:dcomvw2}) imply that
\[
-\frac{1}{2}V_\Omega[\tau](x)+W_\Omega[V_\Omega[\tau]](x)=v_\Omega^-[-\frac{1}{2}\tau +W_\Omega^t[\tau]](x)\qquad\forall x\in\partial\Omega\,,
\]
and accordingly that
\begin{equation}\label{lem:dcomvw3}
v^{-}_\Omega[W_\Omega^t[\tau]]_{|\partial\Omega}=W_\Omega[V_\Omega[\tau]] \quad\text{on} \ \partial\Omega\,.\end{equation}
For the second part of the proof, we set
\[
\tau_0\equiv \tau-\frac{<\tau,1>}{<1,1>}1\,,\qquad\tau_1\equiv \frac{<\tau,1>}{<1,1>}1
\]
and we note that
\[
\tau= \tau_0+\tau_1\,,\qquad<\tau_0,1>=0\,.
\]
Then the function $u_-\equiv v_\Omega^-[\tau_0]$ is harmonic at infinity and has limit at infinity equal to $0$ (cf.~Proposition \ref{prop:sldinfty}). By the Definition \ref{defn:njsdoa} of 
$V^{-1,\alpha}(\partial\Omega)$ and by the membership of $v_\Omega^-[1]$ in $C^{0,\alpha}_{{\mathrm{loc}}}(\overline{\Omega^-})$ (cf.~Theorem \ref{slay}), we know that $u_-\in C^{0,\alpha}_{{\mathrm{loc}}}(\overline{\Omega^-})$.
 Then the third Green Identity in the distributional  form of Theorem \ref{thm:3rdGreen0a} (ii) implies that
\begin{equation}\label{lem:dcomvw4}
w_\Omega^+[u_{-|\partial\Omega}](x)=-v_\Omega^+[\frac{\partial u_-}{\partial\nu_{\Omega^-}}](x)\qquad\forall x\in\Omega^-\,.
\end{equation}
By the classical jump formula for the double layer potential, we know that $w_\Omega^+[u_{|\partial\Omega}]$ has a (unique) extension to $\overline{\Omega}$, which belongs to $C^{0,\alpha}(\overline{\Omega})$ and that
\begin{equation}\label{lem:dcomvw5}
w_\Omega^+[u_{|\partial\Omega}](x)=\frac{1}{2} u(x)+W_\Omega[u_{|\partial\Omega}](x)\qquad\forall x\in\partial\Omega\,,
\end{equation}
(cf.~Theorem \ref{thm:wpsi+-}).  Then equality (\ref{lem:dcomvw4}) implies that $v_\Omega^+[\frac{\partial u}{\partial\nu_{\Omega^-}}]$ admits a unique continuous extension to $\overline{\Omega}$, which is of class $C^{0,\alpha}(\overline{\Omega})$. By the jump formula for the single layer potential with moments in $V^{-1,\alpha}(\partial\Omega)$ of Proposition \ref{prop:ndslv0a}, we have 
\[
\frac{\partial}{\partial\nu_{\Omega^-}}u=-\frac{1}{2}\tau_0-W^t_\Omega[\tau_0] \qquad\text{on}\ \partial\Omega\,.
\] 
Then equalities (\ref{lem:dcomvw4}) and (\ref{lem:dcomvw5}) imply that
\[
\frac{1}{2}V_\Omega[\tau_0](x)+W_\Omega[V_\Omega[\tau_0]](x)=-v_\Omega^+[-\frac{1}{2}\tau_0 -W_\Omega^t[\tau_0]](x)\qquad\forall x\in\partial\Omega\,,
\]
and accordingly that
\begin{equation}\label{lem:dcomvw6}
v^{+}_\Omega[W_\Omega^t[\tau_0]]_{|\partial\Omega}=W_\Omega[V_\Omega[\tau_0]] \quad\text{on} \ \partial\Omega\,.
\end{equation}
On the other hand, $\tau_1$ is a constant and thus we know that $W_\Omega^t[\tau_1]\in C^{0,\alpha}(\partial\Omega)$
(cf.~\textit{e.g.}, \cite[Thm.~6.8]{DaLaMu21}). Then Theorem \ref{slay} ensures that $v^{+}_\Omega[W_\Omega^t[\tau_1]]$ admits a unique continuous extension to $\overline{\Omega}$, which is of class $C^{0,\alpha}(\overline{\Omega})$ and  the classical Plemelj symmetrization principle implies that
\begin{equation}\label{lem:dcomvw7}
v^{+}_\Omega[W_\Omega^t[\tau_1]]_{|\partial\Omega}=W_\Omega[V_\Omega[\tau_1]] \quad\text{on} \ \partial\Omega\,.
\end{equation}
Hence, the sum $v^{+}_\Omega[W_\Omega^t[\tau]]=v^{+}_\Omega[W_\Omega^t[\tau_0]]+v^{+}_\Omega[W_\Omega^t[\tau_1]]$ admits a unique continuous extension to $\overline{\Omega}$, which is of class $C^{0,\alpha}(\overline{\Omega})$  and equalities (\ref{lem:dcomvw6}) and (\ref{lem:dcomvw7}) imply that
\begin{equation}\label{lem:dcomvw8}
v^{+}_\Omega[W_\Omega^t[\tau]]_{|\partial\Omega}=W_\Omega[V_\Omega[\tau]] \quad\text{on} \ \partial\Omega\,.
\end{equation}
For the third part of the proof, we note that 
  equalities (\ref{lem:dcomvw3}) and (\ref{lem:dcomvw8}) imply that
\[
V_\Omega[W_\Omega^t[\tau]]\equiv v^{+}_\Omega[W_\Omega^t[\tau]]_{|\partial\Omega}=v^{-}_\Omega[W_\Omega^t[\tau]]_{|\partial\Omega}=W_\Omega[V_\Omega[\tau]] \quad\text{on} \ \partial\Omega
\]
(see also the Definition \ref{defn:Vpm} of $V_\Omega[\cdot]$).  Hence we have proved that $W_\Omega^t[\tau]$ satisfies the conditions in (\ref{prop:ndslv0a0}) and equality (\ref{lem:dcomvw0}). Finally, we note that the classical Plemelj symmetrization Principle of Lemma \ref{lem:comvw}, the validity of condition condition (\ref{prop:ndslv0a1}) for $\tau$, and equality  (\ref{lem:dcomvw0})  imply that
\begin{eqnarray*}
\lefteqn{
<W_\Omega^t[\tau],V_\Omega[\psi]>
=<\tau,W_\Omega[V_\Omega[\psi]]>
=<\tau,V_\Omega[W_\Omega^t[\psi]]> 
}
\\ \nonumber
&&\qquad\qquad\qquad\quad
=<V_\Omega[\tau], W_\Omega^t[\psi]>
=<W_\Omega[V_\Omega[\tau]],\psi>
\\ \nonumber
&&\qquad\qquad\qquad\quad
=<V_\Omega[W_\Omega^t[\tau]],\psi>\qquad\forall\psi\in C^{0,\alpha}(\partial\Omega)
\end{eqnarray*}
and that accordingly $W_\Omega^t[\tau]$ satisfies  condition (\ref{prop:ndslv0a1}).\hfill  $\Box$ 

\vspace{\baselineskip}

\section{The null space of $\pm\frac{1}{2}I+W_\Omega^t$ and injectivity properties of the single layer in $V^{-1,\alpha}(\partial\Omega)$}
\label{sec:kpmi+wtv}

The null spaces of the restriction of $\pm\frac{1}{2}I+W_\Omega^t$ to $C^{0,\alpha}(\partial\Omega)$ are known (cf.~\textit{e.g.}, \cite[\S 6.5, \S 6.6]{DaLaMu21}). Here we show that the null space of  restriction of $\pm\frac{1}{2}I+W_\Omega^t$ to $V^{-1,\alpha}(\partial\Omega)$ coincides with the   null space of the restriction of $\pm\frac{1}{2}I+W_\Omega^t$ to $C^{0,\alpha}(\partial\Omega)$. To do so, we first introduce  the following elementary classical  lemma. For the convenience of the reader, we include a proof.
\begin{lemma}\label{lem:nzeker-12+wt}
 Let $\alpha\in ]0,1[$. Let $\Omega$ be a bounded open  subset of ${\mathbb{R}}^{2}$ of class $C^{1,\alpha}$. Then there exists
$\tilde{\mu}_1\in C^{0,\alpha}(\partial\Omega)$ such that
 $-\frac{1}{2}\tilde{\mu}_1+W_\Omega^t[\tilde{\mu}_1]=0$ and $<\tilde{\mu}_1,1>\neq 0$.\end{lemma}
{\bf Proof.} If $<\tau,1>=0$ for all $\tau \in C^{0,\alpha}(\partial\Omega)$ such that $-\frac{1}{2}\tau+W_\Omega^t[\tau]=0$, then $1$ would be orthogonal to the kernel of $-\frac{1}{2}I+W_\Omega^t[\cdot]$ in $C^{0,\alpha}(\partial\Omega)$ in the duality pairing in 
\begin{equation}\label{lem:nzeker-12+wt2}
 C^{0,\alpha}(\partial\Omega)\times C^{1,\alpha}(\partial\Omega)\,.
 \end{equation}
Since $W_\Omega[\cdot]$ is compact in $C^{1,\alpha}(\partial\Omega)$ and $W_\Omega^t[\cdot]$ is compact in $C^{0,\alpha}(\partial\Omega)$ (cf.~\textit{e.g.}, \cite[Thm.~6.8]{DaLaMu21}), the Fredholm Alternative  in the above duality pairing in the form of Wendland would imply the existence of a solution $\tau$ in $C^{1,\alpha}(\partial\Omega)$ for the equation
\[
-\frac{1}{2}\tau+W_\Omega[\tau]=1
\]
(cf.~\textit{e.g.}, \cite[Thm.~5.8]{DaLaMu21}). Then the function 
$w_\Omega^-[\tau]$ is harmonic in $\Omega^-$ and harmonic at infinity just as the constant $1$ and accordingly  the uniqueness Theorem for the exterior Dirichlet problem for functions that are harmonic in $\Omega^-$ and harmonic at infinity implies that 
\[
w_\Omega^-[\tau]=1\qquad\text{in}\ \Omega^-\,,
\]
a contradiction. Indeed, the limit of $w_\Omega^-[\tau]$ at infinity  equals $0$.
Hence, there exists $\tilde{\mu}_1\in C^{0,\alpha}(\partial\Omega)$ as in the statement.\hfill  $\Box$ 

\vspace{\baselineskip}

We are now ready to prove the following statement.
\begin{proposition}\label{prop:kpmi+wtv}
 Let   $\alpha\in ]0,1[$. Let $\Omega$ be a bounded open  subset of ${\mathbb{R}}^{n}$ of class $C^{1,\alpha}$. Then
 \begin{eqnarray}\label{prop:kpmi+wtv1}
\lefteqn{
\{\mu\in V^{-1,\alpha}(\partial\Omega):\,
 \pm\frac{1}{2}\mu+W_\Omega^t[S_n\mu]=0\}
}
\\ \nonumber
&&\qquad\qquad\qquad\qquad
=\{\mu\in C^{0,\alpha}(\partial\Omega):\,
\pm\frac{1}{2}\mu+W_\Omega^t[S_n\mu]=0\}\,.
\end{eqnarray}
\end{proposition}
{\bf Proof.} Since $C^{0,\alpha}(\partial\Omega)\subseteq V^{-1,\alpha}(\partial\Omega)$, it suffices to show that the set in the left hand side is contained in the set in the right hand side. Let $\mu\in V^{-1,\alpha}(\partial\Omega)$.

We first consider the case in which $-\frac{1}{2}\mu+W_\Omega^t[\mu]=0$. By  the distributional jump formula for the normal derivative of the harmonic single layer potential of Proposition \ref{prop:ndslv0a}, we have  $\frac{\partial}{\partial\nu}v_\Omega^+[\mu]=0$. Then the uniqueness Theorem \ref{thm:ivneu.uni} for the nonvariational interior Neumann problem implies that $v_\Omega^+[\mu]$ is constant on the connected components of $\Omega$. Since the membership of $\mu$ in $V^{-1,\alpha}(\partial\Omega)$ implies that $v_\Omega^+[\mu]$ has a continuous extension to the closure of $\Omega$, we have
\[
V_{\Omega}[\mu]\in <\chi_{\partial\Omega_j}:\,j\in \{1,\dots,\kappa^+\}>
=  \left\{\tau\in C^{1,\alpha}(\partial\Omega):\, -\frac{1}{2}\tau+W_\Omega[\tau]=0\right\}\,,
\]
(cf.~\textit{e.g.}, \cite[Thm.~6.24]{DaLaMu21}, Theorem 
\ref{thm:Ker-I+W50a} of the Appendix). If $n\geq 3$,  a known isomorphism Theorem implies that there exists a unique $\tilde{\mu}_-\in C^{0,\alpha}(\partial\Omega)$ such that
$v_\Omega^+[\mu]=v_\Omega^+[\tilde{\mu}_-]$ (cf.~\textit{e.g.}, \cite[Prop.~6.19]{DaLaMu21}). Hence, 
\[
v_\Omega^+[\mu-\tilde{\mu}_-]=0\qquad\text{on}\ \overline{\Omega}\,.
\]
Since $v_\Omega^-[\mu-\tilde{\mu}_-]=v_\Omega^-[\mu]-v_\Omega^-[\tilde{\mu}_-]$ is harmonic at infinity and equals 
$v_\Omega^+[\mu-\tilde{\mu}_-]=0$ on $\partial\Omega$, the uniqueness Theorem for the exterior Dirichlet problem for functions that are harmonic in $\Omega^-$ and harmonic at infinity implies that 
\[
v_\Omega^-[\mu-\tilde{\mu}_-]=0\qquad\text{on}\ \overline{\Omega^-}\,.
\]
Then the distributional jump formula for the normal derivative of the harmonic single layer potential of Proposition \ref{prop:ndslv0a}   implies that
\[
0=-\frac{\partial}{\partial\nu_{\Omega^-}}v_\Omega^-[\mu-\tilde{\mu}_-]-\frac{\partial}{\partial\nu}v_\Omega^+[\mu-\tilde{\mu}_-]=\mu-\tilde{\mu}_-\,.
\]
Hence, $\mu=\tilde{\mu}_-\in C^{0,\alpha}(\partial\Omega)$.

If $n=2$, then Lemma \ref{lem:nzeker-12+wt} implies the existence of 
$\tilde{\mu}_1\in C^{0,\alpha}(\partial\Omega)$ such that
 $-\frac{1}{2}\tilde{\mu}_1+W_\Omega^t[\tilde{\mu}_1]=0$ and $<\tilde{\mu}_1,1>\neq 0$. Next we set 
\[
\mu_0\equiv \mu-\frac{<\mu,1>}{<\tilde{\mu}_1,1>}\tilde{\mu}_1\,,\qquad\mu_1\equiv \frac{<\mu,1>}{<\tilde{\mu}_1,1>}\tilde{\mu}_1
\]
and we note that
\[
\mu= \mu_0+\mu_1\,,\qquad<\mu_0,1>=0\,.
\]
By definition of $\mu_1$, we have
\[
 -\frac{1}{2}\mu_0+W_\Omega^t[\mu_0]
=
-\frac{1}{2}\mu+W_\Omega^t[\mu]=0\,,
\]
(cf.~Lemma \ref{lem:Sphi0a}). By a known isomorphism theorem,   there exists a unique $(\tilde{\mu},c)$ in $C^{0,\alpha}(\partial\Omega)_0\times{\mathbb{R}}$ such that
$V_\Omega[\mu_0]=V_\Omega[\tilde{\mu}]+c$ (cf.~\textit{e.g.}, \cite[Prop.~6.47]{DaLaMu21}).     Then we have  
 \[
 V_\Omega[\mu_0-\tilde{\mu}]= c \qquad\text{on}\ \partial\Omega 
 \quad\text{and\ accordingly}\quad
 v_\Omega^+[\mu_0-\tilde{\mu}]= c\qquad\text{on}\ \overline{\Omega}\,.
  \]
Since $<\mu_0-\tilde{\mu},1>=0$, $v_\Omega^-[\mu_0-\tilde{\mu}]$ is harmonic at infinity and equals 
$v_\Omega^+[\mu_0-\tilde{\mu}]=c$ on $\partial\Omega$ (cf.~Proposition \ref{prop:sldinfty}), the uniqueness Theorem for the exterior Dirichlet problem for functions that are harmonic in $\Omega^-$ and harmonic at infinity implies that 
\[
v_\Omega^-[\mu_0-\tilde{\mu}]=c\qquad\text{on}\ \overline{\Omega^-}\,.
\]
Since $<\mu_0-\tilde{\mu},1>=0$,  the jump formula  in distributional form for the normal derivative of the harmonic single layer potential implies that
\[
0=-\frac{\partial}{\partial\nu_{\Omega^-}}v_\Omega^-[\mu_0-\tilde{\mu}]-\frac{\partial}{\partial\nu}v_\Omega^+[\mu_0-\tilde{\mu}]=\mu-\tilde{\mu}\,.
\]
Hence, $\mu_0=\tilde{\mu}\in C^{0,\alpha}(\partial\Omega)$. 
Since $\mu_1\in  C^{0,\alpha}(\partial\Omega)$, we conclude that 
$\mu=\mu_0+\mu_1\in C^{0,\alpha}(\partial\Omega)$.

Next we consider the case in which
 $\frac{1}{2}\mu+W_\Omega^t[\mu]=0$. By  Lemma \ref{lem:Sphi0a}, we have  
$<\mu,1>=0$ and accordingly, $v_\Omega^-[\mu]$ is harmonic at infinity and has limit at infinity equal to zero
(cf.~Proposition \ref{prop:sldinfty}).  Then the jump formula in distributional form for the normal derivative of the harmonic single layer potential implies that $\frac{\partial}{\partial\nu_{\Omega^-}}v_\Omega^-[\mu]=0$. 
  Then the uniqueness Theorem \ref{thm:evneu.uni} for the nonvariational   exterior Neumann problem implies that $v_\Omega^-[\mu]$ is constant on the connected components of $\Omega^-$ if $n\geq 2$ and that  $v_\Omega^-[\mu]$ equals zero on the unbounded connected component of $\Omega^-$ if $n\geq 3$. Since the limit of $v_\Omega^-[\mu]$ at infinity is zero, then $v_\Omega^-[\mu]$ equals zero on the unbounded connected component of $\Omega^-$ also for $n=2$. Since the membership of $\mu$ in $V^{-1,\alpha}(\partial\Omega)$ implies that $v_\Omega^-[\mu]$ has a continuous extension to the closure of $\Omega$, we have
  \[
  V_{\Omega}[\mu]\in <\chi_{\partial(\Omega^-)_j}:\,j\in \{1,\dots,\kappa^-\}>
= \left\{\tau\in C^{1,\alpha}(\partial\Omega):\, \frac{1}{2}\tau+W_\Omega[\tau]=0\right\}\,,
  \]
where we understand that the space  in the right hand side equals  $\{0\}$ if $\kappa^-=0$ (cf.~\textit{e.g.}, \cite[Thm.~6.14]{DaLaMu21}, Theorem 
\ref{thm:KerI+W3} of the Appendix). Then a  known isomorphism theorem implies that there exists a unique $\tilde{\mu}_+\in C^{0,\alpha}(\partial\Omega)$ such that
$\frac{1}{2}\tilde{\mu}_++W_\Omega^t[\tilde{\mu}_+]=0$ and
$V_\Omega[\mu]=V_\Omega[\tilde{\mu}_+]$ (cf.~\textit{e.g.}, \cite[Prop.~6.13]{DaLaMu21}). Hence, 
\[
v_\Omega^+[\mu-\tilde{\mu}_+]=0\qquad\text{on}\ \overline{\Omega}\,.
\]
By Lemma \ref{lem:Sphi0a}, we have $<\tilde{\mu}_+,1>=0$ and accordingly  $<\mu-\tilde{\mu}_+>=0$ and $v_\Omega^-[\mu-\tilde{\mu}_+]$ is harmonic at infinity. 
Since $v_\Omega^-[\mu-\tilde{\mu}_+]$ is harmonic in $\Omega^-$, harmonic at infinity and equals 
$v_\Omega^+[\mu-\tilde{\mu}_+]=0$ on $\partial\Omega$, the uniqueness Theorem for the exterior Dirichlet problem for functions that are harmonic in $\Omega^-$ and harmonic at infinity implies that \[
v_\Omega^-[\mu-\tilde{\mu}_-]=0\qquad\text{on}\ \overline{\Omega^-}\,.
\]
Then the jump formula  in distributional form  for the normal derivative of the harmonic single layer potential implies that
\[
0=-\frac{\partial}{\partial\nu_{\Omega^-}}v_\Omega^-[\mu-\tilde{\mu}_+]-\frac{\partial}{\partial\nu}v_\Omega^+[\mu-\tilde{\mu}_+]=\mu-\tilde{\mu}_+\,.
\]
Hence, $\mu=\tilde{\mu}_+\in C^{0,\alpha}(\partial\Omega)$.\hfill  $\Box$ 

\vspace{\baselineskip}

By the jump formula for the distributional normal derivative of the harmonic single layer potential, we can prove the following injectivity theorem for $v_\Omega^\pm[\cdot]$ and $V_\Omega[\cdot]$ on  $V^{-1,\alpha}(\partial\Omega)$.
\begin{theorem}\label{thm:injslv-1a}
Let   $\alpha\in ]0,1[$. Let $\Omega$ be a bounded open  subset of ${\mathbb{R}}^{n}$ of class $C^{1,\alpha}$. Let $\tau\in V^{-1,\alpha}(\partial\Omega)$, $c\in {\mathbb{R}}$.
\begin{enumerate}
\item[(i)] If $n\geq 3$ and $v_\Omega^+[\tau]=0$ in $\Omega$, then $\tau=0$.
\item[(ii)] If $n\geq 3$ and $v_\Omega^-[\tau]=c$ in $\Omega^-$, then $(\tau,c)=(0,0)$.
\item[(iii)] If $n\geq 3$ and $V_\Omega[\tau]=0$ in $\partial\Omega$, then $\tau=0$.
\item[(iv)] If $n\geq 2$, $<\tau,1>=0$ and $v_\Omega^+[\tau]=c$ in $\Omega$, then $(\tau,c)=(0,0)$.
\item[(v)] If $n\geq 2$, $<\tau,1>=0$ and $v_\Omega^-[\tau]=c$ in $\Omega^-$, then $(\tau,c)=(0,0)$.
\item[(vi)] If $n\geq 2$, $<\tau,1>=0$ and $V_\Omega[\tau]=c$ in $\partial\Omega$, then $(\tau,c)=(0,0)$.
\end{enumerate}
\end{theorem}
{\bf Proof.} (i) By the jump formula for the distributional normal derivative of the single layer potential, we have $-\frac{1}{2}\tau+W_\Omega^t[\tau]=0$. Then Proposition \ref{prop:kpmi+wtv} implies that $\tau\in C^{0,\alpha}(\partial\Omega)$ and the statement follows by \cite[Prop.~6.19]{DaLaMu21}. Statement (iii) is an immediate consequence of statement (i). 

(ii) Since the single layer has limit equal to zero  at infinity, we must have $c=0$. By the jump formula for the normal derivative of the single layer potential, we have $-\frac{1}{2}\tau-W_\Omega^t[\tau]=0$. Then Proposition \ref{prop:kpmi+wtv} implies that $\tau\in C^{0,\alpha}(\partial\Omega)$ and the statement follows by \cite[Prop.~6.13]{DaLaMu21}. 

(iv) By the jump formula for the normal derivative of the single layer potential, we have $-\frac{1}{2}\tau+W_\Omega^t[\tau]=0$. Then Proposition \ref{prop:kpmi+wtv} implies that $\tau\in C^{0,\alpha}(\partial\Omega)$ and the statement follows by \cite[Thm.~6.48 (i)]{DaLaMu21}. Statement (vi) is an immediate consequence of statement (iv).

(v) Since the single layer has limit equal to zero  at infinity for $n\geq 3$, we must have $c=0$ for $n\geq 3$.  By the jump formula for the normal derivative of the single layer potential, we have $-\frac{1}{2}\tau-W_\Omega^t[\tau]=0$. Then Proposition \ref{prop:kpmi+wtv} implies that $\tau\in C^{0,\alpha}(\partial\Omega)$ and the statement follows by \cite[Thm.~6.48 (ii), (iii)]{DaLaMu21}. \hfill  $\Box$ 

\vspace{\baselineskip}

 Then by combining the representation Theorem \ref{thm:slreth} for harmonic functions in terms of single layer potentials with the continuity Corollary \ref{cor:vpm}, with the inclusions of Proposition \ref{prop:v-1apin},
  with the injectivity Theorem \ref{thm:injslv-1a} and with the Open Mapping Theorem, we can prove the following isomorphism theorem, that extends   corresponding classical result in Schauder spaces (cf. \cite[Thms.~6.46--6.48]{DaLaMu21}). 
\begin{theorem}\label{thm:isoslv-1a}
Let   $\alpha\in ]0,1[$. Let $\Omega$ be a bounded open  subset of ${\mathbb{R}}^{n}$ of class $C^{1,\alpha}$. 
\begin{enumerate}
\item[(i)] If $n\geq 3$, then the map from $V^{-1,\alpha,-}(\partial\Omega)$ to $C^{0,\alpha}_h(\overline{\Omega})$, which  takes $\tau$ to
 $v_\Omega^+[\tau]$ is a linear homeomorphism.
\item[(ii)] If $n\geq 3$, then the map from $V^{-1,\alpha,+}(\partial\Omega)$ to $C^{0,\alpha}_{bh}(\overline{\Omega^-})$,
 which  takes $\tau$ to
 $v_\Omega^-[\tau]$ is a linear homeomorphism.
 \item[(iii)] If $n\geq 3$, then the map from $V^{-1,\alpha,\pm}(\partial\Omega)$ to $C^{0,\alpha} (\partial\Omega)$, which  takes $\tau$ to
 $V_\Omega[\tau]$ is a linear homeomorphism.
 \item[(iv)] If $n\geq 2$, then the map from $V^{-1,\alpha,-}(\partial\Omega)_0\times {\mathbb{R}}$ to $C^{0,\alpha}_h(\overline{\Omega})$, which  takes $(\tau,c)$ to
 $v_\Omega^+[\tau]+c$ is a linear homeomorphism.
 \item[(v)] If $n= 2$, then the map from $V^{-1,\alpha,+}(\partial\Omega)_0\times {\mathbb{R}}$ to $C^{0,\alpha}_{hb}(\overline{\Omega^-})$, which  takes $(\tau,c)$ to
 $v_\Omega^-[\tau]+c$ is a linear homeomorphism.
  \item[(vi)] If $n\geq 3$, then the map from $V^{-1,\alpha,-}(\partial\Omega)_0\times{\mathbb{R}}$ to $C^{0,\alpha} (\partial\Omega)$, which  takes $\tau$ to
 $V_\Omega[\tau]+c$ is a linear homeomorphism.
  \item[(vii)] If $n=2$, then the map from $V^{-1,\alpha,\pm}(\partial\Omega)_0\times{\mathbb{R}}$ to $C^{0,\alpha} (\partial\Omega)$, which  takes $\tau$ to
 $V_\Omega[\tau]+c$ is a linear homeomorphism.
  \end{enumerate}
   \end{theorem}
 {\bf Proof.}  By combining the representation Theorem \ref{thm:slreth} with the inclusions of $V^{-1,\alpha,\pm}(\partial\Omega)$ into $ V^{-1,\alpha}(\partial\Omega)$ of Proposition \ref{prop:v-1apin} and the injectivity Theorem \ref{thm:injslv-1a}, the maps of (i), (ii), (iv), (v) are linear isomorphisms and the maps of (iii),  (vi) 
  {\color{black} and (vii)	}	
   are injections. By the existence Theorems \ref{thm:existenceO0} and \ref{thm:existenceO-0} of the Appendix for the interior and exterior Dirichlet problem, each function of 
 $C^{0,\alpha} (\partial\Omega)$ is a restriction of a function $C^{0,\alpha}_{h}(\overline{\Omega})$ and of $C^{0,\alpha}_{hb}(\overline{\Omega^-})$. Then the surjectivity of the maps in statements (i), (ii) imply the surjectivity of the maps in (iii). Similarly, the surjectivity of the maps in statements
   (iv) and (v) imply the surjectivity of the maps in (vi) and (vii). Then
 all the maps of (i)--(vii) are linear isomorphisms and 
   the continuity Corollary \ref{cor:vpm}   and the Open Mapping Theorem imply that 
all the maps of (i)--(vii) are homeomorphisms.\hfill  $\Box$ 

\vspace{\baselineskip}

We also point out the validity of the following corollary  of Theorem \ref{thm:isoslv-1a}.
\begin{corollary}\label{corol:isoslvc-1a}
 Let   $\alpha\in ]0,1[$. Let $\Omega$ be a bounded open  subset of ${\mathbb{R}}^{n}$ of class $C^{1,\alpha}$.
 \begin{enumerate}
\item[(i)] If $n\geq 2$, then the map from $V^{-1,\alpha,-}(\partial\Omega)$ to $C^{0,\alpha}_h(\overline{\Omega})$, which  takes $\tau$ to
 $v_\Omega^+\left[\tau-\frac{<\tau,1>}{<1,1>}1\right]
+
\frac{<\tau,1>}{<1,1>}1$ is a linear homeomorphism.
 \item[(ii)] If $n= 2$, then the map from $V^{-1,\alpha,+}(\partial\Omega)$ to $C^{0,\alpha}_{hb}(\overline{\Omega^-})$, which  takes $\tau$ to
 $v_\Omega^-\left[\tau-\frac{<\tau,1>}{<1,1>}1\right]
+
\frac{<\tau,1>}{<1,1>}1$ is a linear homeomorphism.
  \item[(iii)] If $n\geq 3$, then the map from $V^{-1,\alpha,-}(\partial\Omega)$ to $C^{0,\alpha} (\partial\Omega)$, which  takes $\tau$ to
 $V_\Omega\left[\tau-\frac{<\tau,1>}{<1,1>}1\right]
+
\frac{<\tau,1>}{<1,1>}1$ is a linear homeomorphism.
 \item[(iv)] If $n= 2$, then the map from $V^{-1,\alpha,\pm}(\partial\Omega)$ to $C^{0,\alpha} (\partial\Omega)$, which  takes $\tau$ to
 $V_\Omega\left[\tau-\frac{<\tau,1>}{<1,1>}1\right]
+
\frac{<\tau,1>}{<1,1>}1$ is a linear homeomorphism.
 \end{enumerate}
\end{corollary}
{\bf Proof.} It suffices to note that the map from $V^{-1,\alpha,\pm}(\partial\Omega)$ to   $V^{-1,\alpha,\pm}(\partial\Omega)_0\times {\mathbb{R}}$ that takes $\tau$ to 
\[
\left(\tau-\frac{<\tau,1>}{<1,1>}1,
\frac{<\tau,1>}{<1,1>}1\right)
\]
is a linear homeomorphism and to apply Theorem \ref{thm:isoslv-1a} (iv), (v), (vi), (vii).\hfill  $\Box$ 

\vspace{\baselineskip}

\section{A characterization of the space $V^{-1,\alpha}(\partial\Omega)$}
\label{sec:v-1achar}
We are now able to prove the following characterization theorem for the space  $V^{-1,\alpha}(\partial\Omega)$.
\begin{theorem}\label{thm:v-1achar}
 Let  $\alpha\in ]0,1[$. Let $\Omega$ be a bounded open  subset of ${\mathbb{R}}^{n}$ of class $C^{1,\alpha}$. 
Then
 \begin{equation}\label{thm:v-1achar1}
V^{-1,\alpha}(\partial\Omega)=V^{-1,\alpha,+}(\partial\Omega)=V^{-1,\alpha,-}(\partial\Omega)\,.
\end{equation}
Moreover, the norms of $V^{-1,\alpha,+}(\partial\Omega)$ and of $V^{-1,\alpha,-}(\partial\Omega)$ generate the same topology on $V^{-1,\alpha}(\partial\Omega)$.
\end{theorem}
{\bf Proof.} By Proposition \ref{prop:v-1apin} it suffices to show that 
$V^{-1,\alpha}(\partial\Omega)$ is contained in $V^{-1,\alpha,\pm}(\partial\Omega)$  and we consider case  $n\geq 3$ and case $n=2$ separately. Let $\tau\in V^{-1,\alpha}(\partial\Omega)$.  

If $n\geq 3$, then the Definition \ref{defn:njsdoa} of $V^{-1,\alpha}(\partial\Omega)$ implies  that $V_\Omega[\tau]$ belongs to $C^{0,\alpha}(\partial\Omega)$ and accordingly 
the representation Theorem \ref{thm:isoslv-1a} (iii) implies the existence of $\tau_\pm\in V^{-1,\alpha,\pm}(\partial\Omega)$ such that
\[
V_\Omega[\tau]=V_\Omega[\tau_\pm]\qquad\text{on}\ \partial\Omega\,.
\]
Then the injectivity Theorem \ref{thm:injslv-1a} (iii) implies that $\tau=\tau_\pm$ and accordingly
$\tau\in V^{-1,\alpha,\pm}(\partial\Omega)$.

If $n=2$, then then the Definition \ref{defn:njsdoa} of $V^{-1,\alpha}(\partial\Omega)$ implies  that $V_\Omega[\tau]$ belongs to $C^{0,\alpha}(\partial\Omega)$. Since $1\in V^{-1,\alpha}(\partial\Omega)$ and 
$V_\Omega[1]\in C^{0,\alpha}(\partial\Omega)$, we have
\[
V_\Omega\left[\tau-\frac{<\tau,1>}{<1,1>}1\right] +
\frac{<\tau,1>}{<1,1>}1
\in  C^{0,\alpha}(\partial\Omega)\,.
\]
 Then the representation Corollary \ref{corol:isoslvc-1a} (iv) implies that there exists
 $\tau_\pm$ in the space $ V^{-1,\alpha,\pm}(\partial\Omega)$ such that
 \[
 V_\Omega\left[\tau-\frac{<\tau,1>}{<1,1>}1\right] +
\frac{<\tau,1>}{<1,1>}1
=V_\Omega\left[\tau_\pm-\frac{<\tau_\pm,1>}{<1,1>}1\right] +
\frac{<\tau_\pm,1>}{<1,1>}1\,.
 \]
 Then the injectivity Theorem \ref{thm:injslv-1a} (vi) implies that 
 \[
 \tau-\frac{<\tau,1>}{<1,1>}1=\tau_\pm-\frac{<\tau_\pm,1>}{<1,1>}1\,,
 \qquad
 \frac{<\tau,1>}{<1,1>}1=\frac{<\tau_\pm,1>}{<1,1>}1
 \] and accordingly $\tau=\tau_\pm$. Hence, 
$\tau\in V^{-1,\alpha,\pm}(\partial\Omega)$.

If $n\geq 3$, then Theorem \ref{thm:isoslv-1a} (iii) implies that 
the norms of $V^{-1,\alpha,+}(\partial\Omega)$ and $V^{-1,\alpha,-}(\partial\Omega)$ are equivalent
on $V^{-1,\alpha,-}(\partial\Omega)=V^{-1,\alpha,+}(\partial\Omega)$.

If $n=2$, then Theorem \ref{thm:isoslv-1a} (vii) implies that the norms of $V^{-1,\alpha,+}(\partial\Omega)_0\times{\mathbb{R}}$ and $V^{-1,\alpha,-}(\partial\Omega)_0\times{\mathbb{R}}$ are equivalent. Since
the map from $V^{-1,\alpha,\pm}(\partial\Omega)_0\times{\mathbb{R}}$ to
$V^{-1,\alpha,\pm}(\partial\Omega)$ that takes $(\tau,c)$ to $\tau+c$ is a linear homeomorphism, we conclude that the norms of $V^{-1,\alpha,+}(\partial\Omega)$ and $V^{-1,\alpha,-}(\partial\Omega)$ are equivalent
on $V^{-1,\alpha,-}(\partial\Omega)=V^{-1,\alpha,+}(\partial\Omega)$.\hfill  $\Box$ 

\vspace{\baselineskip}

By exploiting the isomorphisms of Theorem \ref{thm:isoslv-1a}, Corollary \ref{corol:isoslvc-1a} and the characterization of Theorem \ref{thm:v-1achar}, we can introduce a norm on $V^{-1,\alpha}(\partial\Omega)$ that is equivalent to the norms of $V^{-1,\alpha,\pm}(\partial\Omega)$. We do so by means of the following theorem.
\begin{theorem}\label{thm:v-1an}
 Let $\alpha\in]0,1[$.  Let  $\Omega$ be a  bounded open subset of ${\mathbb{R}}^{n}$ of class $C^{1,\alpha}$. 
 \begin{enumerate}
\item[(i)] If  $n\geq 3$, then the map $\|\cdot\|_{V^{-1,\alpha}(\partial\Omega)}$ from $V^{-1,\alpha}(\partial\Omega)$ to $[0,+\infty[$ that is defined by
\[
\|\tau\|_{V^{-1,\alpha}(\partial\Omega)}\equiv
\|V_\Omega[\tau]\|_{C^{0,\alpha}(\partial\Omega)}\qquad\forall\tau\in V^{-1,\alpha}(\partial\Omega)\,,
\]
is a norm on $V^{-1,\alpha}(\partial\Omega)$ that is equivalent to the norms 
$\|\cdot\|_{V^{-1,\alpha,\pm}(\partial\Omega)}$ on $V^{-1,\alpha}(\partial\Omega)$. In particular, $\left(V^{-1,\alpha}(\partial\Omega),\|\cdot\|_{V^{-1,\alpha}(\partial\Omega)}\right)$ is a Banach space and $V_\Omega[\cdot]$ is a linear isometry from $\left(V^{-1,\alpha}(\partial\Omega),\|\cdot\|_{V^{-1,\alpha}(\partial\Omega)}\right)$ onto $C^{0,\alpha}(\partial\Omega)$.
\item[(ii)] If  $n\geq 2$, then the map $\|\cdot\|_{V^{-1,\alpha}(\partial\Omega)}'$ from $V^{-1,\alpha}(\partial\Omega)$ to $[0,+\infty[$ that is defined by
\begin{eqnarray*}
\lefteqn{
\|\tau\|_{V^{-1,\alpha}(\partial\Omega)}'
}
\\ \nonumber
&&\quad\equiv
\left\|
V_\Omega\left[\tau-\frac{<\tau,1>}{<1,1>}1\right]
+
\frac{<\tau,1>}{<1,1>}1
\right\|_{C^{0,\alpha}(\partial\Omega)}\quad\forall\tau\in V^{-1,\alpha}(\partial\Omega)\,,
\end{eqnarray*}
is a norm on $V^{-1,\alpha}(\partial\Omega)$ that is equivalent to the norms 
$\|\cdot\|_{V^{-1,\alpha,\pm}(\partial\Omega)}$ on $V^{-1,\alpha}(\partial\Omega)$. In particular, $\left(V^{-1,\alpha}(\partial\Omega),\|\cdot\|_{V^{-1,\alpha}(\partial\Omega)}'\right)$ is a Banach space and the operator
 $J^{-1,\alpha}$ 
from  $\left(V^{-1,\alpha}(\partial\Omega),\|\cdot\|_{V^{-1,\alpha}(\partial\Omega)}'\right)$ onto $C^{0,\alpha}(\partial\Omega)$ that is defined by 
\begin{equation}\label{thm:v-1an1}
J^{-1,\alpha}[\tau]\equiv V_\Omega\left[\tau-\frac{<\tau,1>}{<1,1>}1\right]
+
\frac{<\tau,1>}{<1,1>}1 
\qquad\forall \tau\in V^{-1,\alpha}(\partial\Omega)
\end{equation}
is a linear isometry  .
\end{enumerate}
\end{theorem}
{\bf Proof.} (i) By Theorem \ref{thm:isoslv-1a} (iii), $V_\Omega[\cdot]$ is a linear homeomorphism from the space $V^{-1,\alpha,\pm}(\partial\Omega)$ onto 
$C^{0,\alpha}(\partial\Omega)$. By the definition of norm $\|\cdot\|_{V^{-1,\alpha}(\partial\Omega)}$, $V_\Omega[\cdot]$ is a linear isometry from $\left(V^{-1,\alpha}(\partial\Omega),\|\cdot\|_{V^{-1,\alpha}(\partial\Omega)}\right)$ onto $C^{0,\alpha}(\partial\Omega)$. Hence, statement (i) holds true. 

(ii) By  Theorem \ref{thm:v-1achar}, we have $V^{-1,\alpha,-}(\partial\Omega)=V^{-1,\alpha,+}(\partial\Omega)$ both algebraically and topologically. Then
Corollary \ref{corol:isoslvc-1a} (iii), (iv) implies that the map $J^{-1,\alpha}$ from $V^{-1,\alpha,\pm}(\partial\Omega)$ onto $C^{0,\alpha} (\partial\Omega)$  is a linear homeomorphism. By the definition of norm $\|\cdot\|_{V^{-1,\alpha}(\partial\Omega)}'$, 
 the map $J^{-1,\alpha}$ is a linear isometry from $\left(V^{-1,\alpha}(\partial\Omega),\|\cdot\|_{V^{-1,\alpha}(\partial\Omega)}\right)$ onto $C^{0,\alpha}(\partial\Omega)$. Hence, statement (ii) holds true. \hfill  $\Box$ 

\vspace{\baselineskip}

\begin{remark}{\em
 Let $\alpha\in]0,1[$.  Let  $\Omega$ be a  bounded open subset of ${\mathbb{R}}^{n}$ of class $C^{1,\alpha}$. 
Since $V^{-1,\alpha}(\partial\Omega)$ coincides both algebraically and topologically with $V^{-1,\alpha,\pm}(\partial\Omega)$, Corollary \ref{cor:vpm} (iii) implies that $V_\Omega[\cdot]$ is linear and continuous from $V^{-1,\alpha}(\partial\Omega)$ to $C^{0,\alpha}(\partial\Omega)$. 
}\end{remark}

{\color{black}
\begin{theorem}\label{thm:contidnuv-1a}
 Let $\alpha\in ]0,1[$. Let $\Omega$ be a bounded open  subset of ${\mathbb{R}}^{n}$ of class $C^{1,\alpha}$. Then the  following statements hold.
 \begin{enumerate}
\item[(i)] The map $\frac{\partial}{\partial\nu}$ from the closed subspace $C^{0,\alpha}_h(\overline{\Omega}) $ of $ C^{0,\alpha}(\overline{\Omega})$ to $V^{-1,\alpha}(\partial\Omega)$ is linear and continuous. Moreover, formula (\ref{eq:trispo1}) holds true. 
\item[(ii)] The   map $\frac{\partial u}{\partial\nu_{\Omega^-}}$ from the closed subspace
$
C^{0,\alpha}_{bh}(\overline{\Omega^-}) $ of $ C^{0,\alpha}_b(\overline{\Omega^-})$ to $V^{-1,\alpha}(\partial\Omega)$ is linear and continuous. Moreover, formula (\ref{eq:trespo2}) holds true.
\end{enumerate}
\end{theorem}
{\bf Proof.} (i) By Theorem \ref{thm:contidnu}, we know that $\frac{\partial}{\partial\nu}$ is linear and continuous from $C^{0,\alpha}_h(\overline{\Omega}) $ to $\left(C^{1,\alpha}(\partial\Omega)\right)'$
and that formula (\ref{eq:trispo1}) holds true. Then formula (\ref{eq:trispo1}), the continuity of the restriction from $C^{0,\alpha}_h(\overline{\Omega}) $ to $C^{0,\alpha}(\partial\Omega)$ and the definition of the norm in $V^{-1,\alpha, +}(\partial\Omega)$, $\frac{\partial}{\partial\nu}$ is linear and continuous from $C^{0,\alpha}_h(\overline{\Omega}) $ to $V^{-1,\alpha, +}(\partial\Omega)$,  that coincides both algebraically and topologically with $V^{-1,\alpha}(\partial\Omega)$ (cf. Theorem \ref{thm:v-1achar}). The proof of (ii) follows that lines of that of (i), by invoking Theorem \ref{thm:contednu} instead of Theorem  \ref{thm:contidnu}.\hfill  $\Box$ 

\vspace{\baselineskip}
}

The following lemma is useful in order to verify the continuity of a linear operator with values in $V^{-1,\alpha}(\partial\Omega)$.
\begin{lemma}\label{lem:cov-1a}
 Let $\alpha\in]0,1[$.  Let  $\Omega$ be a  bounded open subset of ${\mathbb{R}}^{n}$ of class $C^{1,\alpha}$.  Let $X$ be a normed space. Let $L$ be a linear map from $X$ to $V^{-1,\alpha}(\partial\Omega)$. Then the following statements hold.
 \begin{enumerate}
\item[(i)]  
 $L$ is continuous if and only if the map from $X$ to $C^{0,\alpha}(\partial\Omega)$ that takes $x$ to
 \[
 J^{-1,\alpha}[L[x]]=V_\Omega\left[L[x]-\frac{< L[x],1>}{<1,1>}1\right]
+
\frac{<L[x],1>}{<1,1>}1
\]
from $X$ to $C^{0,\alpha}(\partial\Omega)$ is 
 continuous.
 \item[(ii)]    $L$ is compact if and only if the map from $X$ to $C^{0,\alpha}(\partial\Omega)$ that takes $\tau$ to $J^{-1,\alpha}[L[\tau]]$ is compact.
\end{enumerate}
\end{lemma}
{\bf Proof.} By Theorem \ref{thm:v-1an} (ii), $J^{-1,\alpha}$ is a linear isometry. Thus $L$ is continuous if and only if $J^{-1,\alpha}\circ L$ is continuous and 
$L$ is compact if and only if $J^{-1,\alpha}\circ L$ is compact.\hfill  $\Box$ 

\vspace{\baselineskip}

Then we can readily prove the following embedding Theorem.
\begin{theorem}\label{thm:coim-1ab}
 Let $\alpha,\beta\in]0,1[$, $\beta<\alpha$.  Let  $\Omega$ be a  bounded open subset of ${\mathbb{R}}^{n}$ of class $C^{1,\alpha}$. Then $V^{-1,\alpha}(\partial\Omega)$ is compactly embedded into $V^{-1,\beta}(\partial\Omega)$.
 \end{theorem}
{\bf Proof.} By Theorem \ref{thm:v-1an} (ii), $J^{-1,\alpha}$ is a linear isometry from $V^{-1,\alpha}(\partial\Omega)$ to $C^{0,\alpha}(\partial\Omega)$, that is compactly imbedded into $C^{0,\beta}(\partial\Omega)$. Then $J^{-1,\alpha}$ is compact from $V^{-1,\alpha}(\partial\Omega)$ to $C^{0,\beta}(\partial\Omega)$. Since
\[
J^{-1,\beta}[\tau]=J^{-1,\alpha}[\tau]\qquad\forall \tau\in V^{-1,\alpha}(\partial\Omega)\,,
\]
$J^{-1,\beta}$ is compact from $V^{-1,\alpha}(\partial\Omega)$ to $C^{0,\beta}(\partial\Omega)$. Then 
by applying Lemma \ref{lem:cov-1a} (ii) with $X=V^{-1,\alpha}(\partial\Omega)$ and with $L$ equal to the embedding of $V^{-1,\alpha}(\partial\Omega)$ into $V^{-1,\beta}(\partial\Omega)$, 
we conclude that  $V^{-1,\alpha}(\partial\Omega)$ is compactly embedded into $V^{-1,\beta}(\partial\Omega)$.\hfill  $\Box$ 

\vspace{\baselineskip}

\section{The transpose of the double layer potential in the duality pairing
 $\left(V^{-1,\alpha}(\partial\Omega),C^{1,\alpha} (\partial\Omega)\right)$
}\label{sec:trdlwd}
Since $V^{-1,\alpha}(\partial\Omega)$ is a subspace of the dual $\left(C^{1,\alpha} (\partial\Omega)\right)'$, we can restrict the canonical duality pairing of $\left(\left(C^{1,\alpha} (\partial\Omega)\right)',C^{1,\alpha} (\partial\Omega)\right)$ to $\left(V^{-1,\alpha}(\partial\Omega),C^{1,\alpha} (\partial\Omega)\right)$ in the spirit of the Wendland duality approach of
\cite{We67}, \cite{We70} (cf.~Kress \cite{Kr14}). By the  Plemelj's symmetrization principle in distributional form of Lemma \ref{lem:dcomvw}, the operator  
$W_{\Omega}^t[\cdot]$ maps $V^{-1,\alpha}(\partial\Omega)$ to itself. 
We now turn to show that  $W_{\Omega}^t[\cdot]$  is compact in the space $\left(V^{-1,\alpha}(\partial\Omega),\|\cdot\|_{V^{-1,\alpha}(\partial\Omega)}'\right)$. In order to do so, we first prove the following statement.
\begin{theorem}\label{thm:wt-1ac}
 Let  $\alpha,\beta\in]0,1[$, $\beta\leq\alpha$.  Let  $\Omega$ be a  bounded open subset of ${\mathbb{R}}^{n}$ of class $C^{1,\alpha}$. 
 Let $W_{\Omega}^t[\cdot]$ be the transpose operator to $W_{\Omega}[\cdot]$ in $C^{1,\alpha} (\partial\Omega)$. 
 Then the operator $W_{\Omega}^t[\cdot]$ is linear and continuous from  $V^{-1,\beta}(\partial\Omega)$ to $V^{-1,\alpha}(\partial\Omega)$.
\end{theorem}
{\bf Proof.} We first show that $J^{-1,\beta}W_{\Omega}^t[\cdot]$ is linear and continuous from 
$V^{-1,\beta}(\partial\Omega)$ to $C^{0,\alpha}(\partial\Omega)$. To do so, we show that $J^{-1,\beta}W_{\Omega}^t[\cdot]$ equals a continuous operator from 
$V^{-1,\beta}(\partial\Omega)$ to $C^{0,\alpha}(\partial\Omega)$. Let
 $\tau\in V^{-1,\beta}(\partial\Omega)$. Then the Plemelj symmetrization principle of Lemma \ref{lem:dcomvw} implies that
\begin{eqnarray}\label{thm:wt-1ac1}
\lefteqn{
 J^{-1,\beta}W_{\Omega}^t[\tau]
}
\\ \nonumber
&& 
=
V_\Omega\left[W_{\Omega}^t[\tau]-\frac{<W_{\Omega}^t[\tau],1>}{<1,1>}1\right]
+
\frac{<W_{\Omega}^t[\tau],1>}{<1,1>}1
\\ \nonumber
&& 
=
V_\Omega\left[W_{\Omega}^t[\tau]\right]
-V_\Omega\left[\frac{<W_{\Omega}^t[\tau],1>}{<1,1>}1\right]
+
\frac{<W_{\Omega}^t[\tau],1>}{<1,1>}1
\\ \nonumber
&& 
=
W_\Omega\left[V_{\Omega}[\tau]\right]
-\frac{<\tau,W_{\Omega}[1]>}{<1,1>}V_\Omega\left[1\right]
+
\frac{<\tau,W_{\Omega}[1]>}{<1,1>}1
\\ \nonumber
&& 
=
W_\Omega\left[
V_\Omega\left[\tau-\frac{<\tau,1>}{<1,1>}1\right]
+
\frac{<\tau,1>}{<1,1>}1\right]
\\ \nonumber
&&\quad 
+W_\Omega\left[
V_\Omega\left[\frac{<\tau,1>}{<1,1>}1\right]\right]
-W_\Omega\left[\frac{<\tau,1>}{<1,1>}1\right]
\\ \nonumber
&&\quad
-\frac{<\tau,W_{\Omega}[1]>}{<1,1>}V_\Omega\left[1\right]
+
\frac{<\tau,W_{\Omega}[1]>}{<1,1>}1
\\ \nonumber
&& 
=
W_\Omega\left[
V_\Omega\left[\tau-\frac{<\tau,1>}{<1,1>}1\right]
+
\frac{<\tau,1>}{<1,1>}1\right]
\\ \nonumber
&&\quad 
+W_\Omega\left[
V_\Omega\left[ 1\right]\right]\frac{<\tau,1>}{<1,1>}
-\frac{1}{2}\frac{<\tau,1>}{<1,1>}
\\ \nonumber
&&\quad
-\frac{1}{2}\frac{<\tau,1>}{<1,1>}V_\Omega\left[1\right]
+
\frac{1}{2}\frac{<\tau,1>}{<1,1>} 
\\ \nonumber
&& 
=
W_\Omega\left[
J^{-1,\beta}[\tau]\right]
\\ \nonumber
&&\quad 
+\left(
W_\Omega\left[V_\Omega\left[ 1\right]\right]
-\frac{1}{2}V_\Omega\left[1\right] \right)\frac{<\tau,1>}{<1,1>} \,.
\end{eqnarray}
Since $J^{-1,\beta}$ is an isometry from $V^{-1,\beta}(\partial\Omega)$ onto $C^{0,\beta}(\partial\Omega)$ and $W_\Omega[\cdot]$ is linear and continuous from $C^{0,\beta}(\partial\Omega)$ to $C^{0,\alpha}(\partial\Omega)$ (cf.~Theorem \ref{thm:schauin} (ii)), we conclude that $W_\Omega\left[
J^{-1,\beta}[\cdot]\right]$ is linear and continuous from $V^{-1,\beta}(\partial\Omega)$ to $C^{0,\alpha}(\partial\Omega)$.
Since the factor
\[
\left(
W_\Omega\left[V_\Omega\left[1\right]\right]
-\frac{1}{2}V_\Omega\left[1\right] \right)\in C^{1,\alpha}(\partial\Omega)
\]
 and the map from $V^{-1,\beta}(\partial\Omega)$  to ${\mathbb{R}}$ that takes $\tau$ to $<\tau,1>$ is linear and continuous, we conclude that the map from $V^{-1,\beta}(\partial\Omega)$  to $C^{0,\alpha}(\partial\Omega)$ that takes $\tau$ to the right hand side of (\ref{thm:wt-1ac1}) is continuous.
 
 Next we prove that  $W_\Omega^t[V^{-1,\beta}(\partial\Omega)]$ is contained in $V^{-1,\alpha}(\partial\Omega)$. To do so, we note that
 \begin{equation}\label{thm:wt-1ac2}
J^{-1,\beta}[\mu]=J^{-1,\alpha}[\mu]\qquad\forall \mu\in V^{-1,\alpha}(\partial\Omega)\,,
\end{equation}
that $J^{-1,\beta}$ is a bijection from $V^{-1,\beta}(\partial\Omega)$ onto $C^{0,\beta}(\partial\Omega)$, that $C^{0,\alpha}(\partial\Omega)\subseteq C^{0,\beta}(\partial\Omega)$, that $J^{-1,\alpha}$ is a bijection from $V^{-1,\alpha}(\partial\Omega)$ onto $C^{0,\alpha}(\partial\Omega)$ and that accordingly the inclusion of $J^{-1,\beta}W_\Omega^t[V^{-1,\beta}(\partial\Omega)]$ in $C^{0,\alpha}(\partial\Omega)$ implies that
$W_\Omega^t[V^{-1,\beta}(\partial\Omega)]$ is contained in $V^{-1,\alpha}(\partial\Omega)$. Then equalities (\ref{thm:wt-1ac1}) and (\ref{thm:wt-1ac2}) and the above continuity of  $J^{-1,\beta}W_\Omega^t$ imply that $J^{-1,\alpha}W_\Omega^t=J^{-1,\beta}W_\Omega^t$ is linear and continuous from 
$V^{-1,\beta}(\partial\Omega)$ to $C^{0,\alpha}(\partial\Omega)$ and thus   Lemma \ref{lem:cov-1a} (i) implies the continuity of $W_\Omega^t$ from $V^{-1,\beta}(\partial\Omega)$  to $V^{-1,\alpha}(\partial\Omega)$   and thus the proof is complete.\hfill  $\Box$ 

\vspace{\baselineskip}

By the compact embedding of $V^{-1,\alpha}(\partial\Omega)$ into $V^{-1,\beta}(\partial\Omega)$ of Theorem \ref{thm:coim-1ab} and by Theorem \ref{thm:wt-1ac}, we conclude that the following statement holds true.
\begin{corollary}\label{corol:thm:wt-1acp}
 Let  $\alpha\in]0,1[$.  Let  $\Omega$ be a  bounded open subset of ${\mathbb{R}}^{n}$ of class $C^{1,\alpha}$. Then 
the operator $W_{\Omega}^t[\cdot]$ is compact from  $V^{-1,\alpha}(\partial\Omega)$ to itself.
\end{corollary}
 
\section{Existence of solutions for a nonvariational form of the Neumann Problem in $\Omega$ and $\Omega^-$}\label{sec:eineu}

We now prove the existence of solutions to the interior and exterior Neumann problems. We begin with Theorem \ref{thm:nvINexist},  where we consider the interior case. 
\begin{theorem}\label{thm:nvINexist}
Let $\alpha\in]0,1[$.  Let  $\Omega$ be a  bounded open subset of ${\mathbb{R}}^{n}$ of class $C^{1,\alpha}$. 
If 
{\color{black}$g\in V^{-1,\alpha}(\partial\Omega)$	} 
 satisfies the compatibility condition \eqref{lem:nvINcc2}, then the integral equation
\begin{equation}\label{thm:nvINexist1}
\left(-\frac{1}{2}I+W^t_\Omega\right)[\phi]=g
\end{equation}
has a solution $\phi\in V^{-1,\alpha}(\partial\Omega)$ and the distributional single layer potential
\[
v^+_\Omega[\phi]\in C^{0,\alpha}(\overline\Omega)
\]
is a solution of the interior Neumann boundary value problem (\ref{eq:nvneupb}). The set of all solutions of (\ref{eq:nvneupb}) in $C^0(\overline{\Omega})$ consists of the functions $u\in C^{0,\alpha}(\overline\Omega)$ such that the difference $u-v^+_\Omega[\phi]$ is constant on each connected component $\Omega_1$, \dots, $\Omega_{\kappa^+}$ of $\Omega$.
\end{theorem}
{\bf Proof.} If {\color{black}$g\in V^{-1,\alpha}(\partial\Omega)$	}  satisfies the compatibility condition (\ref{lem:nvINcc2}), then  it is orthogonal to the null space of the integral operator $-\frac{1}{2}I+W_\Omega$ in $C^{0,\alpha}(\partial\Omega)$ (cf. Theorem 
{\color{black}\ref{thm:Ker-I+W50a}		}
 of the Appendix). Since $W_\Omega[\cdot]$ is compact in $C^{0,\alpha}(\partial\Omega)$ and $W_\Omega^t[\cdot]$ is compact in $V^{-1,\alpha}(\partial\Omega)$ (cf.~Theorem \ref{thm:schauin} and Corollary \ref{corol:thm:wt-1acp}), the Fredholm Alternative Theorem 
of Wendland \cite{We67}, \cite{We70} (cf.~\textit{e.g.}, \cite[Thm.~5.8]{DaLaMu21})  implies that equation (\ref{thm:nvINexist1})  has a solution $\phi\in V^{-1,\alpha}(\partial\Omega)$. Then, Proposition \ref{prop:sldinfty} (i) implies that  $v^+_\Omega[\phi]$ is harmonic in $\Omega$ and the Definition \ref{defn:njsdoa} of $V^{-1,\alpha}(\partial\Omega)$ implies that 
 it is a function of $C^{0,\alpha}(\overline{\Omega})$. By the jump formula for the normal derivative of single layer potential of Proposition \ref{prop:ndslv0a}, it follows that $v^+_\Omega[\phi]$ is a solution of (\ref{eq:nvneupb}). The last statement of the theorem is a consequence of the uniqueness result of Theorem \ref{thm:ivneu.uni}.\hfill  $\Box$ 

\vspace{\baselineskip}

For the exterior Neumann problem we find convenient to describe separately the case of dimension $n\ge 3$ and the case of dimension $n=2$.

\begin{theorem}\label{thm:nvENexist3} Let $n\in {\mathbb{N}}$, $n\geq 3$, $\alpha\in]0,1[$.  Let  $\Omega$ be a  bounded open subset of ${\mathbb{R}}^{n}$ of class $C^{1,\alpha}$.  If  
{\color{black}$g\in V^{-1,\alpha}(\partial\Omega)$	} 
satisfies the compatibility condition (\ref{lem:nvENcc.eq2}), then the integral equation 
\[
\left(\frac{1}{2}I+W^t_\Omega\right)[\phi]=g
\]
has a solution $\phi\in V^{-1,\alpha}(\partial\Omega)$ and the distributional single layer
\[
v^-_\Omega[\phi]\in C^{0,\alpha}_{b}(\overline{\Omega^-})
\]
is a solution of the exterior Neumann boundary value problem (\ref{eq:nveneupb}). The set of all solutions of (\ref{eq:nveneupb}) in $C^{0}_{b}(\overline{\Omega^-})$
consists of the functions $u\in C^{0,\alpha}_{b}(\overline{\Omega^-})$ such that the difference $u-v^-_\Omega[\phi]$ is constant on each bounded connected component $(\Omega^-)_1$, \dots, $(\Omega^-)_{\kappa^-}$ of $\Omega^-$ and  is identically equal to $0$ on the unbounded connected component $(\Omega^-)_0$. 
\end{theorem}
{\bf Proof.}  The proof follows the lines to that of Theorem \ref{thm:nvINexist}. We just have to invoke  Theorem {\color{black}\ref{thm:KerI+W30a} 	}
instead of Theorem 
{\color{black}\ref{thm:Ker-I+W50a}		}
  of the Appendix and notice that $v^-_\Omega[\phi]$ is harmonic at infinity by Proposition \ref{prop:sldinfty} (ii). The last statement of the theorem follows by the uniqueness result of Theorem \ref{thm:evneu.uni}.\hfill  $\Box$ 

\vspace{\baselineskip}

Finally, for the two-dimensional case we have the following statement.

\begin{theorem}\label{thm:nvENexist2}
Let $n=2$, $\alpha\in]0,1[$.  Let  $\Omega$ be a  bounded open subset of ${\mathbb{R}}^{n}$ of class $C^{1,\alpha}$.   If  $g\in V^{-1,\alpha}(\partial\Omega)$ satisfies the compatibility conditions (\ref{lem:nvENcc.eq2}) and (\ref{lem:nvENcc.eq3}), then the integral equation 
\begin{equation}\label{thm:nvENexist2.eq1}
\left(\frac{1}{2}I+W^t_\Omega\right)[\phi]=g
\end{equation}
has a solution $\phi\in V^{-1,\alpha}(\partial\Omega)_0$ and the distributional single layer
\[
v^-_\Omega[\phi]\in C^{0,\alpha}_{b}(\overline{\Omega^-})
\]
is a solution of the exterior Neumann boundary value problem (\ref{eq:nveneupb}) (see notation (\ref{eq:dx0})). Moreover, the set of all solutions of (\ref{eq:nveneupb}) in $C^{0}_{b}(\overline{\Omega^-})$ consists of the functions $u\in C^{0,\alpha}_{b}(\overline{\Omega^-})$ such that the difference $u-v^-_\Omega[\phi]$ is constant on each (bounded and unbounded) connected component  $(\Omega^-)_0$, \dots, $(\Omega^-)_{\kappa^-}$ of $\Omega^-$. 
\end{theorem}
{\bf Proof.} The proof follows the lines of those of Theorems \ref{thm:nvINexist} and \ref{thm:nvENexist3}. We only observe that, if $g$ satisfies (\ref{lem:nvENcc.eq2}) and (\ref{lem:nvENcc.eq3}) and $\phi$ is a solution of (\ref{thm:nvENexist2.eq1}), then Lemma \ref{lem:Sphi0a} implies that 
\[
0=<g,1>=<\left(\frac{1}{2}I+W^t_\Omega\right)[\phi],1>=<\phi,1>\,.
\] 
Accordingly, the jump formula of (\ref{prop:ndslv0a3}) and Proposition  \ref{prop:sldinfty}  (iii) imply that $v^-_\Omega[\phi]$ satisfies the boundary condition of the Neumann problem (\ref{eq:nveneupb}) and that 
$v^-_\Omega[\phi]$ is harmonic in $\Omega^-$ and harmonic at infinity.\hfill  $\Box$ 

\vspace{\baselineskip}

 \appendix

\section{Appendix: Classical   solvability of the Dirichlet problem with data in  $C^{0,\alpha}(\partial\Omega)$}
  \label{sec:clsohad} 
  
 We now present some classical results on the 
 solvability of the Dirichlet problem for  H\"{o}lder continuous harmonic functions in the specific form that we need in the paper. We first note   that  the following classical known result holds true  (cf.~\textit{e.g.}, \cite[Thm.~6.14]{DaLaMu21}).
\begin{theorem}\label{thm:KerI+W3}
 Let $\alpha\in ]0,1[$. Let $\Omega$ be a bounded open  subset of ${\mathbb{R}}^{n}$ of class $C^{1,\alpha}$.   Then the following statements hold true.
 \begin{enumerate}
\item[(i)] If $\varkappa^->0$, then the null space 
\begin{equation}\label{thm:KerI+W31}
\left\{\mu\in C^{1,\alpha}(\partial\Omega):\,\left(\frac{1}{2}I+W_{\Omega}\right)[\mu]=0\right\}
\end{equation}
consists of the functions from $\partial\Omega$ to $\mathbb{R}$ that are constant on $\partial(\Omega^-)_j$ for all $j\in\{1,\dots,\varkappa^-\}$ and that are identically equal to $0$ on $\partial(\Omega^-)_0$.
\item[(ii)] If $\varkappa^-=0$, then
\begin{equation}\label{thm:KerI+W32}
\left\{\mu\in C^{1,\alpha}(\partial\Omega):\,\left(\frac{1}{2}I+W_{\Omega}\right)[\mu]=0\right\}=\{0\}\,.
\end{equation}
\end{enumerate}
\end{theorem}
We now turn to prove that the equalities (\ref{thm:KerI+W31}), (\ref{thm:KerI+W32}) hold also in case we replace the membership of $\mu$ in $C^{1,\alpha}(\partial\Omega)$ with that of $\mu$ in $C^{m,\alpha}(\partial\Omega)$ for any natural number $m$, including $m=0$ in the sense of the following statement.
\begin{theorem}\label{thm:KerI+W30a}
 Let $m\in {\mathbb{N}}$, $\alpha\in ]0,1[$. Let $\Omega$ be a bounded open  subset of ${\mathbb{R}}^{n}$ of class $C^{\max\{m,1\},\alpha}$.   Then the following statements hold true.
 \begin{enumerate}
\item[(i)] If $\varkappa^->0$, then the null space 
\begin{equation}\label{thm:KerI+W30a1}
\left\{\mu\in C^{m,\alpha}(\partial\Omega):\,\left(\frac{1}{2}I+W_{\Omega}\right)[\mu]=0\right\}
\end{equation}
consists of the functions from $\partial\Omega$ to $\mathbb{R}$ that are constant on $\partial(\Omega^-)_j$ for all $j\in\{1,\dots,\varkappa^-\}$ and that are identically equal to $0$ on $\partial(\Omega^-)_0$.
\item[(ii)] If $\varkappa^-=0$, then
\begin{equation}\label{thm:KerI+W30a2}
\left\{\mu\in C^{m,\alpha}(\partial\Omega):\,\left(\frac{1}{2}I+W_{\Omega}\right)[\mu]=0\right\}=\{0\}\,.
\end{equation}
\end{enumerate}
\end{theorem}
{\bf Proof.} If $m\geq 1$, then $C^{m,\alpha}(\partial\Omega)\subseteq
C^{1,\alpha}(\partial\Omega)$ and the statements (i), (ii)  are an immediate consequence of Theorem \ref{thm:KerI+W3}. We now consider case $m=0$. 
It suffices to prove that if $\mu\in C^{0,\alpha}(\partial\Omega)$ and
\begin{equation}\label{thm:KerI+W30a3}
\frac{1}{2}\mu+W_{\Omega}[\mu]=0
\end{equation}
then  $\mu\in C^{1,\alpha}(\partial\Omega)$. We do so by means of known regularization results for $W_{\Omega}$. We first note that equality (\ref{thm:KerI+W30a3}) and a simple inductive argument imply that
\[
(-1/2)^j\mu=W_{\Omega}^j[\mu]\qquad\forall j\in {\mathbb{N}}\setminus\{0\}\,.
\]
Now let 
\[
M\equiv\max\{j\in {\mathbb{N}}\setminus\{0\}:\,j\alpha < 1\}\,.
\] 
By applying $M$ times the classical result  \cite[Thm.~1 (i)]{La24}, we deduce that  $W_{\Omega}^{M}[\mu]$ belongs to $C^{0,M\alpha}(\partial\Omega)$.

If $(M+1)\alpha=1$,  then     
 \cite[Thm.~1 (ii)]{La24} implies  that $W_{\Omega}^{M+1}[\mu]$ belongs to 
$C^{0,\gamma}(\partial\Omega)$ for all $\gamma\in]1-\alpha,1[$. We now fix $\gamma\in]1-\alpha,1[$.
Then $\gamma+\alpha>1$ and \cite[Thm.~2 (i)]{La24} implies that
$W_{\Omega}^{M+2}[\mu]$ belongs to $C^{1,\gamma+\alpha-1}(\partial\Omega)$. Hence,  Thm.~9.2 (ii) of Dondi and the author
\cite{DoLa17} ensures that  $(-1/2)^{M+3}\mu=W_{\Omega}^{M+3}[\mu]$ belongs to $C^{1, \alpha}(\partial\Omega)$.

If instead $(M+1)\alpha>1$, then  
 \cite[Thm.~2 (ii)]{La24} implies  that $W_{\Omega}^{M+1}[\mu]$ belongs to $C^{1,M\alpha+\alpha-1}(\partial\Omega)$.  Hence, 
Thm.~9.2 (ii) of Dondi and the author
\cite{DoLa17} ensures that  $(-1/2)^{M+2}\mu=W_{\Omega}^{M+2}[\mu]$ belongs to $C^{1, \alpha}(\partial\Omega)$ and thus the proof is complete.
\hfill  $\Box$ 

\vspace{\baselineskip}

Next we prove the following decomposition statement (for the  classical statement in $L^2(\partial\Omega)$ {\color{black}with $\Omega$ of class $C^2$,	}
 we refer to Folland \cite[Cor.~3.39]{Fo95} and for a corresponding  statement in $C^{1,\alpha}(\partial\Omega)$ 
 {\color{black}with $\Omega$ of class $C^{1,\alpha}$,	}
 we refer to \cite[Cor.~6.17]{DaLaMu21}). {\color{black} The proof is a variant of the proof of
 Folland \cite[Cor.~3.39]{Fo95}.		 }
\begin{proposition}\label{prop.KerI+W50a}  Let $m\in {\mathbb{N}}$, $\alpha\in ]0,1[$. Let $\Omega$ be a bounded open  subset of ${\mathbb{R}}^{n}$ of class $C^{\max\{m,1\},\alpha}$.      Then 
 \begin{equation}\label{prop.KerI+W50a1}
C^{m,\alpha}(\partial\Omega)=\mathrm{Im}\left(\frac{1}{2}I+W_\Omega\right)_{|C^{m,\alpha}(\partial\Omega)}
\oplus\mathrm{Ker}\left(\frac{1}{2}I+W_\Omega\right)_{|C^{m,\alpha}(\partial\Omega)}\,, 
\end{equation}
where the sums are direct but not necessarily orthogonal.
\end{proposition}
{\bf Proof.} By Corollary 9.1 of reference \cite{DoLa17} of Dondi and the author,    $\frac{1}{2}I+W_\Omega$ is a compact perturbation of the identity in $C^{m,\alpha}(\partial\Omega)$. Then  $\frac{1}{2}I+W_\Omega$ is a Fredholm operator of index zero in $C^{m,\alpha}(\partial\Omega)$. Moreover,  $\mathrm{Ker}\left(\frac{1}{2}I+W_\Omega\right)$ is closed in $C^{m,\alpha}(\partial\Omega)$ and has dimension $\varkappa^-$ by  Theorem  \ref{thm:KerI+W30a}. Then the image of $\frac{1}{2}I+W_\Omega$  is a closed subspace of $C^{m,\alpha}(\partial\Omega)$ and has co-dimension $\varkappa^-$. So, it remains to show that $\mathrm{Im}\left(\frac{1}{2}I+W_\Omega\right)$  and $\mathrm{Ker}\left(\frac{1}{2}I+W_\Omega\right)$ have a trivial intersection. Let 
\[
\phi\in \mathrm{Im}\left(\frac{1}{2}I+W_\Omega\right)_{|C^{m,\alpha}(\partial\Omega)}\cap \mathrm{Ker}\left(\frac{1}{2}I+W_\Omega\right)_{|C^{m,\alpha}(\partial\Omega)}\,.
\]
By Theorem \ref{thm:KerI+W30a}, we have  $\phi\in C^{1,\alpha}(\partial\Omega)$. Then \cite[Cor.~6.17]{DaLaMu21} implies that there exist uniquely determined elements $\phi_1\in 
\{ \frac{1}{2}\psi+W^t_\Omega[\psi]:\,\psi\in C^{0,\alpha}(\partial\Omega)\}$ and $\phi_2\in C^{1,\alpha}(\partial\Omega)$ such that $\frac{1}{2}\phi_2+W^t_\Omega[\phi_2]=0$
 and $\phi=\phi_1+\phi_2$.  In particular, there exists $\psi_1\in C^{0,\alpha}(\partial\Omega)$ such that $\phi_1=\frac{1}{2}\psi_1+W^t_\Omega[\psi_1]$. But then 
 \[
 \int_{\partial\Omega}\phi_1\phi d\sigma=\int_{\partial\Omega}\frac{1}{2}\psi_1\phi+W^t_\Omega[\psi_1]\phi\,d\sigma=\int_{\partial\Omega}\psi_1\left(
 \frac{1}{2}\phi+W_\Omega[\phi]
 \right)\,d\sigma=0\,,
 \]
 by the membership of $\phi$ in $\mathrm{Ker}\left(\frac{1}{2}I+W_\Omega\right)_{|C^{m,\alpha}(\partial\Omega)}$. Moreover, 
   the membership of  $\phi$ in
  $\mathrm{Im}\left(\frac{1}{2}I+W_\Omega\right)_{|C^{m,\alpha}(\partial\Omega)}$ implies the existence of $\psi\in C^{m,\alpha}(\partial\Omega)$ such that $\phi= \frac{1}{2}\psi+W_\Omega[\psi]$ and accodingly
 \[
 \int_{\partial\Omega}\phi_2\phi d\sigma=\int_{\partial\Omega}\phi_2\frac{1}{2}\psi
 +\phi_2W_\Omega[\psi] d\sigma
 \int_{\partial\Omega}\left(\frac{1}{2}\phi_2+W_\Omega^t[\phi_2]\right)\psi\,d\sigma
 =0\,,
 \]
 by equality $\frac{1}{2}\phi_2+W^t_\Omega[\phi_2]=0$.  Hence, $\int_{\partial\Omega}(\phi_1+\phi_2) \phi\,d\sigma=0$ and $\phi=0$.\hfill  $\Box$ 

\vspace{\baselineskip}
  
 Then we have the following variant of \cite[Prop.~6.13]{DaLaMu21}. 
\begin{proposition}\label{prop:KerI+W2} 
 Let $m\in {\mathbb{N}}\setminus\{0\}$, $\alpha\in ]0,1[$. Let $\Omega$ be a bounded open  subset of ${\mathbb{R}}^{n}$ of class $C^{m,\alpha}$.   The map from $\mathrm{Ker}(\frac{1}{2}I+W^t_\Omega)_{|C^{m-1,\alpha}(\partial\Omega)}$ to $\mathrm{Ker}(\frac{1}{2}I+{W}_\Omega)_{|C^{m,\alpha}(\partial\Omega)}$ that takes $\mu$ to $V_\Omega[\mu]$ is an isomorphism.
\end{proposition}
{\bf Proof.} Since $C^{m-1,\alpha}(\partial\Omega)\subseteq C^{0,\alpha}(\partial\Omega)$, Proposition 6.13 of \cite{DaLaMu21} 
and Theorem \ref{thm:KerI+W30a} 
imply that $V_{\Omega}$  maps $\mathrm{Ker}(\frac{1}{2}I+W^t_\Omega)_{|C^{m-1,\alpha}(\partial\Omega)}$ into 
\[
\mathrm{Ker}(\frac{1}{2}I+{W}_\Omega)_{|C^{1,\alpha}(\partial\Omega)}
=\mathrm{Ker}(\frac{1}{2}I+{W}_\Omega)_{|C^{m,\alpha}(\partial\Omega)}\,.
\]
By Dondi and the author \cite[Cors.~9.1, 10.1]{DoLa17}, ${W}_\Omega$ is compact in $C^{m,\alpha}(\partial\Omega)$ and ${W}_\Omega^t$ is compact in $C^{m-1,\alpha}(\partial\Omega)$ and thus the Fredholm Alternative Theorem in the duality pairing of
\begin{equation}\label{prop:KerI+W21}
C^{m-1,\alpha}(\partial\Omega)\times C^{m,\alpha}(\partial\Omega)
\end{equation}
  in the form of Wendland \cite{We67}, \cite{We70} (cf. \cite[Thm.~5.8]{DaLaMu21})  implies that the null spaces $\mathrm{Ker}(\frac{1}{2}I+W^t_\Omega)_{|C^{m-1,\alpha}(\partial\Omega)}$ and  $\mathrm{Ker}(\frac{1}{2}I+{W}_\Omega)_{|C^{m,\alpha}(\partial\Omega)}$ have the same finite dimension. As a consequence, it suffices to show that the map that takes $\mu$ to $V_\Omega[\mu]$ is injective to prove that it is an isomorphism. By \cite[Prop.~6.13]{DaLaMu21}, such a map is injective on $\mathrm{Ker}(\frac{1}{2}I+W^t_\Omega)_{|C^{0,\alpha}(\partial\Omega)}$ and accordingly on its subset $\mathrm{Ker}(\frac{1}{2}I+W^t_\Omega)_{|C^{m-1,\alpha}(\partial\Omega)}$.\hfill  $\Box$ 

\vspace{\baselineskip}

Then by means of the decomposition of Proposition \ref{prop.KerI+W50a}, one can solve the interior Dirichlet problem for the Laplace equation with boundary data in $C^{m,\alpha}(\partial\Omega)$ by following the scheme of the proof of the case of boundary data in $C^{1,\alpha}(\partial\Omega)$ of \cite[Thm.~6.27]{DaLaMu21}. We do by means of the following statement.
\begin{theorem}\label{thm:existenceO0}
Let $m\in {\mathbb{N}}$, $\alpha\in ]0,1[$. Let $\Omega$ be a bounded open  subset of ${\mathbb{R}}^{n}$ of class $C^{\max\{m,1\},\alpha}$. If $g\in C^{m,\alpha}(\partial\Omega)$, then the boundary value problem
\begin{equation}\label{thm:existenceO00}
\begin{cases}
\Delta u=0&\text{in }\Omega\,,\\
u=g&\text{on }\partial\Omega
\end{cases}
\end{equation}
has one and only one solution $u\in C^{m,\alpha}(\overline\Omega)$. Moreover, there exist
 $\psi\in C^{m,\alpha}(\partial\Omega)$, $\mu\in C^{\max\{m,1\}-1,\alpha}(\partial\Omega)$ such that
 \[
 u=w_\Omega^+[\psi]+v_\Omega^+[\mu]\,.
 \] 
In particular, the solution $u$ equals the sum of a double layer potential of class $C^{m,\alpha}(\overline\Omega)$ and of a single layer potential of class $C^{\max\{m,1\},\alpha}(\overline\Omega)$. 
\end{theorem}
{\bf Proof.}  By the direct sum in (\ref{prop.KerI+W50a1}) there exist   $g_1\in \mathrm{Im}(\frac{1}{2}I+{W}_\Omega)_{|C^{m,\alpha}(\partial\Omega)}$ and  $g_2\in \mathrm{Ker}(\frac{1}{2}I+{W}_\Omega)_{|C^{m,\alpha}(\partial\Omega)}$ such that 
\[
g=g_1+g_2
\]
(here $g_2=0$ when $\varkappa^-=0$).
Then there exists $\psi\in C^{m,\alpha}(\partial\Omega)$ such that
\[
\frac{1}{2}\psi+{W}_\Omega[\psi]=g_1
\] 
and by the jump formula for the double layer potential (cf.~\textit{e.g.}, \cite[Thm.~6.6]{DaLaMu21}), it follows that 
\begin{equation}\label{thm:existenceO01}
w_\Omega^+[\psi]_{|\partial\Omega}=g_1\,.
\end{equation}
Moreover, by the membership of $g_2$ in 
\[
\mathrm{Ker}(\frac{1}{2}I+{W}_\Omega)_{|C^{m,\alpha}(\partial\Omega)}
=\mathrm{Ker}(\frac{1}{2}I+{W}_\Omega)_{|C^{\max\{m,1\},\alpha}(\partial\Omega)}\,,
\]
there exists a function $\mu\in\mathrm{Ker}(\frac{1}{2}I+{W}_\Omega^t)_{|C^{\max\{m,1\}-1,\alpha}(\partial\Omega)}$ such that 
\begin{equation}\label{thm:existenceO02}
V_\Omega[\mu]=v_\Omega[\mu]_{|\partial\Omega}=g_2\,,
\end{equation}
(cf.~Theorem \ref{thm:KerI+W30a} and  Proposition \ref{prop:KerI+W2}). Since the single and double  layer potentials are harmonic in $\Omega$ (cf.~\textit{e.g.}, \cite[Thm.~4.22 (i), Prop.~ 4.28 (i)]{DaLaMu21}),   equalities (\ref{thm:existenceO01}) and (\ref{thm:existenceO02}),    and   the regularity Theorems for the single and double harmonic layer potentials (see Theorems \ref{slay}, \ref{thm:wpsi+-}), imply that    $w_\Omega^+[\psi]\in C^{m,\alpha}(\overline\Omega)$,  
$v^+_\Omega[\mu]\in C^{\max\{m,1\},\alpha}(\overline\Omega)$, and that
\[
u\equiv w_\Omega^+[\psi]+v^+_\Omega[\mu]
\]
is a solution of the Dirichlet problem (\ref{thm:existenceO00}) and belongs to $C^{m,\alpha}(\overline{\Omega})$ (cf.~\textit{e.g.}, Proposition \ref{prop:sldinfty}, \cite[Prop.~4.28]{DaLaMu21}). The uniqueness of $u$ is a consequence of the Maximum Principle.\hfill  $\Box$ 

\vspace{\baselineskip}

Next we note that  the following classical known result holds true  (cf.~\textit{e.g.}, \cite[Thm.~6.24]{DaLaMu21}).
\begin{theorem}\label{thm:Ker-I+W5}
Let $\alpha\in ]0,1[$. Let $\Omega$ be a bounded open  subset of ${\mathbb{R}}^{n}$ of class $C^{1,\alpha}$.   
Then the null space 
\begin{equation}\label{thm:Ker-I+W51}
\left\{\mu\in C^{1,\alpha}(\partial\Omega):\,\left(-\frac{1}{2}I+W_{\Omega}\right)[\mu]=0\right\}
\end{equation}
consists of the functions from $\partial\Omega$ to $\mathbb{R}$ which are constant on $\partial\Omega_j$ for all $j\in\{1,\dots,\varkappa^+\}$.
\end{theorem}
By exploiting the very same regularization argument of the proof of Theorem \ref{thm:KerI+W30a}, one can prove
  that the equality (\ref{thm:Ker-I+W51})  holds also in case we replace the membership of $\mu$ in $C^{1,\alpha}(\partial\Omega)$ with that of $\mu$ in $C^{m,\alpha}(\partial\Omega)$ for any natural number $m$, including $m=0$ in the sense of the following statement.\begin{theorem}\label{thm:Ker-I+W50a}
Let $m\in {\mathbb{N}}$, $\alpha\in ]0,1[$. Let $\Omega$ be a bounded open  subset of ${\mathbb{R}}^{n}$ of class $C^{\max\{m,1\},\alpha}$.      
Then the null space 
\begin{equation}\label{thm:Ker-I+W50a1}
\left\{\mu\in C^{m,\alpha}(\partial\Omega):\,\left(-\frac{1}{2}I+W_{\Omega}\right)[\mu]=0\right\}
\end{equation}
consists of the functions from $\partial\Omega$ to $\mathbb{R}$ which are constant on $\partial\Omega_j$ for all $j\in\{1,\dots,\varkappa^+\}$.
\end{theorem}
Next we prove the following decomposition statement 
The proof can be obtained by a straightforward modification of that of  Proposition \ref{prop.KerI+W50a} and using Theorems  \ref{thm:Ker-I+W50a} and   \cite[Cor.~6.26]{DaLaMu21}  instead of Theorems \ref{thm:KerI+W30a} and \cite[Cor.~6.17]{DaLaMu21} 
(for the classical statement in $L^2(\partial\Omega)$ {\color{black}with $\Omega$ of class $C^2$,	} we refer to Folland \cite[Corollary 3.39]{Fo95} and for a corresponding  statement in $C^{1,\alpha}(\partial\Omega)$  {\color{black}with $\Omega$ of class $C^{1,\alpha}$,	} we refer to \cite[Cor.~6.26]{DaLaMu21}).
\begin{proposition}\label{prop.Ker-I+W80a} 
Let $m\in {\mathbb{N}}$, $\alpha\in ]0,1[$. Let $\Omega$ be a bounded open  subset of ${\mathbb{R}}^{n}$ of class $C^{\max\{m,1\},\alpha}$.    Then 
 \begin{equation}\label{prop.Ker-I+W80a1}
C^{m,\alpha}(\partial\Omega)=\mathrm{Im}\left(-\frac{1}{2}I+W_\Omega\right)_{|C^{m,\alpha}(\partial\Omega)}\oplus\mathrm{Ker}\left(-\frac{1}{2}I+W_\Omega\right)_{|C^{m,\alpha}(\partial\Omega)}\,, 
\end{equation}
where the sums are direct but not necessarily orthogonal. \end{proposition}
Then we have the following variant of \cite[Props.~6.19, 6.23]{DaLaMu21}.
\begin{proposition}\label{prop:Ker-I+W2}
Let $m\in {\mathbb{N}}\setminus\{0\}$, $\alpha\in ]0,1[$. Let $\Omega$ be a bounded open  subset of ${\mathbb{R}}^{n}$ of class $C^{m,\alpha}$.   Then the following statements hold.
\begin{enumerate}
\item[(i)]  If $n\geq 3$, then the map from $\mathrm{Ker}(-\frac{1}{2}I+W^t_\Omega)_{|C^{m-1,\alpha}(\partial\Omega)}$ to $\mathrm{Ker}(-\frac{1}{2}I+{W}_\Omega)_{|C^{m,\alpha}(\partial\Omega)}$ that takes $\mu$ to $V_\Omega[\mu]$ is an isomorphism.
\item[(ii)] The null space $\mathrm{Ker}(-\frac{1}{2}I+{W}_\Omega)_{|C^{m,\alpha}(\partial\Omega)}$ is the topological direct sum of 
\[
V_\Omega\left[\left(\mathrm{Ker}\left(-\frac{1}{2}I+W^t_\Omega\right)_{|C^{m-1,\alpha}(\partial\Omega)}\right)_0\right]
\]
 and of the space of constant functions from $\partial\Omega$ to $\mathbb{R}$ (see notation (\ref{eq:x0})).
 \end{enumerate}
\end{proposition}\label{prop:bvps.Ker-I+W2}
{\bf Proof.} (i) Since $C^{m-1,\alpha}(\partial\Omega)\subseteq C^{0,\alpha}(\partial\Omega)$, Proposition 6.19 of \cite{DaLaMu21} 
and Theorem \ref{thm:Ker-I+W50a} 
imply that $V_{\Omega}$  maps $\mathrm{Ker}(-\frac{1}{2}I+W^t_\Omega)_{|C^{m-1,\alpha}(\partial\Omega)}$ into 
\[
\mathrm{Ker}(-\frac{1}{2}I+{W}_\Omega)_{|C^{1,\alpha}(\partial\Omega)}
=\mathrm{Ker}(-\frac{1}{2}I+{W}_\Omega)_{|C^{m,\alpha}(\partial\Omega)}\,.
\]
By Dondi and the author \cite[Cors.~9.1, 10.1]{DoLa17}, ${W}_\Omega$ is compact in $C^{m,\alpha}(\partial\Omega)$ and ${W}_\Omega^t$ is compact in $C^{m-1,\alpha}(\partial\Omega)$ and thus the Fredholm Alternative Theorem in the duality pairing of (\ref{prop:KerI+W21}) in the form of Wendland \cite{We67}, \cite{We70} (cf. \cite[Thm.~5.8]{DaLaMu21})  implies that the null spaces $\mathrm{Ker}(-\frac{1}{2}I+W^t_\Omega)_{|C^{m-1,\alpha}(\partial\Omega)}$ and  $\mathrm{Ker}(-\frac{1}{2}I+{W}_\Omega)_{|C^{m,\alpha}(\partial\Omega)}$ have the same finite dimension. As a consequence, it suffices to show that the map that takes $\mu$ to $V_\Omega[\mu]$ is injective to prove that it is an isomorphism. By \cite[Prop.~6.19]{DaLaMu21}, such a map is injective on $\mathrm{Ker}(\frac{1}{2}I+W^t_\Omega)_{|C^{0,\alpha}(\partial\Omega)}$ and accordingly on its subset $\mathrm{Ker}(\frac{1}{2}I+W^t_\Omega)_{|C^{m-1,\alpha}(\partial\Omega)}$.

(ii) The closed subspace $\left(\mathrm{Ker}\left(-\frac{1}{2}I+W^t_\Omega\right)_{|C^{m-1,\alpha}(\partial\Omega)}\right)_0$ of  
\[
\mathrm{Ker}\left(-\frac{1}{2}I+W^t_\Omega\right)_{|C^{m-1,\alpha}(\partial\Omega)}\]
 is defined by the vanishing of the linear functional that takes $\phi$ to $\int_{\partial\Omega}\phi\, d\sigma$, therefore it has co-dimension smaller than or equal to one. By Dondi and the author \cite[Thms.~9.2 (ii), 10.1(ii)]{DoLa17}, ${W}_\Omega$ is compact in $C^{m,\alpha}(\partial\Omega)$ and ${W}_\Omega^t$ is compact in $C^{m-1,\alpha}(\partial\Omega)$ and thus the Fredholm Alternative Theorem in the duality pairing of (\ref{prop:KerI+W21}) in the form of Wendland \cite{We67}, \cite{We70} (cf. \cite[Thm.~5.8]{DaLaMu21})  implies that the null spaces 
$\mathrm{Ker}\left(-\frac{1}{2}I+W^t_\Omega\right)_{|C^{m-1,\alpha}(\partial\Omega)}$ and $\mathrm{Ker}\left(-\frac{1}{2}I+W_\Omega\right)_{|C^{m,\alpha}(\partial\Omega)}$ have the same finite dimension. By Proposition 6.20 of  \cite{DaLaMu21}, the map $V_\Omega$ is injective from $\left(\mathrm{Ker}\left(-\frac{1}{2}I+W^t_\Omega\right)_{|C^{0,\alpha}(\partial\Omega)}\right)_0$ to 
$\mathrm{Ker}\left(-\frac{1}{2}I+W_\Omega\right)_{|C^{1,\alpha}(\partial\Omega)}$  and accordingly on its subset $\left(\mathrm{Ker}\left(-\frac{1}{2}I+W^t_\Omega\right)_{|C^{m-1,\alpha}(\partial\Omega)}\right)_0$. Hence, the subspace  $\left(\mathrm{Ker}\left(-\frac{1}{2}I+W^t_\Omega\right)_{|C^{m-1,\alpha}(\partial\Omega)}\right)_0$ of  $\mathrm{Ker}\left(-\frac{1}{2}I+W^t_\Omega\right)_{|C^{m-1,\alpha}(\partial\Omega)}$ and the subspace $V_\Omega\left[\left(\mathrm{Ker}\left(-\frac{1}{2}I+W^t_\Omega\right)_{|C^{m-1,\alpha}(\partial\Omega)}\right)_0\right]$ of $\mathrm{Ker}(-\frac{1}{2}I+{W}_\Omega)_{|C^{m,\alpha}(\partial\Omega)}$ have the same (finite) dimension and co-dimension.  In particular, 
\[
V_\Omega\left[\left(\mathrm{Ker}\left(-\frac{1}{2}I+W^t_\Omega\right)_{|C^{m-1,\alpha}(\partial\Omega)}\right)_0\right]
\]
 has co-dimension less than or equal to one in $\mathrm{Ker}(-\frac{1}{2}I+{W}_\Omega)_{|C^{m,\alpha}(\partial\Omega)}$. 
By Theorem \ref{thm:Ker-I+W50a}, the constant functions on $\partial\Omega$ are contained in $\mathrm{Ker}(-\frac{1}{2}I+{W}_\Omega)_{|C^{m,\alpha}(\partial\Omega)}$. Instead,  Lemma 6.22 of \cite{DaLaMu21}  implies that $V_\Omega\left[\left(\mathrm{Ker}\left(-\frac{1}{2}I+W^t_\Omega\right)_{|C^{m-1,\alpha}(\partial\Omega)}\right)_0\right]$ does not contain non-trivial constant functions. Since the space of constant functions on $\partial\Omega$ has dimension one, we conclude that $V_\Omega\left[\left(\mathrm{Ker}\left(-\frac{1}{2}I+W^t_\Omega\right)_{|C^{m-1,\alpha}(\partial\Omega)}\right)_0\right]$ has co-dimension exactly equal to one and that $\mathrm{Ker}(-\frac{1}{2}I+{W}_\Omega)_{|C^{m,\alpha}(\partial\Omega)}$ is the algebraic direct sum of $V_\Omega\left[\left(\mathrm{Ker}\left(-\frac{1}{2}I+W^t_\Omega\right)_{|C^{m-1,\alpha}(\partial\Omega)}\right)_0\right]$ and of the space of constant functions. Finally, since all the spaces involved have finite dimension, they are closed and the algebraic direct sum is topological (cf.~\textit{e.g.}, \cite[Thms.~2.3, 2.6]{DaLaMu21}).\hfill  $\Box$ 

\vspace{\baselineskip}

 By means of the decomposition of Proposition \ref{prop.Ker-I+W80a}, one could solve the exterior Dirichlet problem for the Laplace equation with boundary data in $C^{m,\alpha}(\partial\Omega)$ by following the scheme of the proof of the case of boundary data in $C^{1,\alpha}(\partial\Omega)$ of \cite[Thm.~6.33]{DaLaMu21}. We do so by means of the following statement.
 \begin{theorem}\label{thm:existenceO-0} 
 Let $m\in {\mathbb{N}}$, $\alpha\in ]0,1[$. Let $\Omega$ be a bounded open  subset of ${\mathbb{R}}^{n}$ of class $C^{\max\{m,1\},\alpha}$.      
Let $g\in C^{m,\alpha}(\partial\Omega)$. The boundary value problem
\begin{equation}\label{thm:existenceO-00}
\begin{cases}
\Delta u=0&\text{in }\Omega^-\,,\\
u=g&\text{on }\partial\Omega\,,\\
u\text{ is harmonic at infinity} 
\end{cases}
\end{equation}
has one and only one solution $u\in C^{m,\alpha}_{b}(\overline{\Omega^-})$. Moreover, if $n\geq 3$, then there exist
 $\psi\in C^{m,\alpha}(\partial\Omega)$, $\mu\in C^{\max\{m,1\}-1,\alpha}(\partial\Omega)$ such that
 \[
 u=w_\Omega^-[\psi]+v_\Omega^-[\mu]\,,
  \] 
and if $n
{\color{black}=	}	
2$, then  there exist
 $\psi\in C^{m,\alpha}(\partial\Omega)$, $\mu\in C^{\max\{m,1\}-1,\alpha}(\partial\Omega)_0$, $c\in {\mathbb{R}}$ such that
 \[
 u=w_\Omega^-[\psi]+v_\Omega^-[\mu]+c\,.
 \] 
In particular, the solution $u$ equals the sum of a double layer potential of class $C^{m,\alpha}_{b}(\overline{\Omega^-})$ and of a single layer potential of class $C^{\max\{m,1\},\alpha}_{b}(\overline{\Omega^-})$ for $n\geq 3$ 
and  the sum of a double layer potential of class $C^{m,\alpha}_{b}(\overline{\Omega^-})$, of a single layer potential of class $C^{\max\{m,1\},\alpha}_{b}(\overline{\Omega^-})$  and of a constant for $n=2$. 
\end{theorem}
 {\bf Proof.} By the direct sum of Proposition \ref{prop.Ker-I+W80a} there exist   $g_1\in \mathrm{Im}(-\frac{1}{2}I+{W}_\Omega)_{|C^{m,\alpha}(\partial\Omega)}$ and  $g_2\in \mathrm{Ker}(-\frac{1}{2}I+{W}_\Omega)_{|C^{m,\alpha}(\partial\Omega)}$ such that 
\[
g=g_1+g_2\,.
\]
We first consider the Dirichlet problem with boundary datum $g_1$. Since $g_1$ belongs to $\mathrm{Im}(-\frac{1}{2}I+{W}_\Omega)_{|C^{m,\alpha}(\partial\Omega)}$,  there exists $\psi\in C^{m,\alpha}(\partial\Omega)$ such that
\[
-\frac{1}{2}\psi+{W}_\Omega[\psi]=g_1
\] 
and the jump formula for the double layer potential implies that 
\[
w_\Omega^-[\psi]_{|\partial\Omega}=g_1\,,
\]
{\color{black} (cf.~\textit{e.g.}, Theorem \ref{thm:wpsi+-} (iii)).	}
Moreover, the function $w_\Omega^-[\psi]$ is harmonic and harmonic at infinity and belongs to  $C^{m,\alpha}_{b}(\overline{\Omega^-})$ by classical results on the regularity of double layer potentials  (see Theorem \ref{thm:wpsi+-} (ii),  \cite[Prop.~4.28]{DaLaMu21}).
 
We now turn to the problem with datum $g_2$ and we consider separately the cases of dimension $n\ge 3$ and $n=2$. For $n\geq 3$, the membership of $g_2$ in 
\[
\mathrm{Ker}(-\frac{1}{2}I+{W}_\Omega)_{|C^{m,\alpha}(\partial\Omega)}
=\mathrm{Ker}(-\frac{1}{2}I+{W}_\Omega)_{|C^{\max\{m,1\},\alpha}(\partial\Omega)}\,,
\]
implies  that
 there exists a function $\mu$ in $\mathrm{Ker}(-\frac{1}{2}I+{W}_\Omega^t)_{|C^{\max\{m,1\}-1,\alpha}(\partial\Omega)}$ such that 
 \[
V_\Omega[\mu]=v_\Omega[\mu]_{|\partial\Omega}=g_2\,.
\]
(cf.~Theorem \ref{thm:Ker-I+W50a} and Proposition \ref{prop:Ker-I+W2} (i)). Since  the single layer $v^-_\Omega[\mu]$ is harmonic, harmonic at infinity, and solves the Dirichlet problem with boundary datum $g_2$ and $v^-_\Omega[\mu]\in C^{\max\{m,1\},\alpha}_{b}(\overline{\Omega^-})$, the function   
\[
u\equiv w_\Omega^-[\psi]+v^-_\Omega[\mu]
\]
is a solution of the exterior Dirichlet problem (\ref{thm:existenceO-00}) and belongs to $C^{m,\alpha}_{b}(\overline{\Omega^-})$
(cf.~\textit{e.g.}, Proposition \ref{prop:sldinfty}, Theorem \ref{thm:slayh}).

Finally, we consider the  case when $n=2$. The membership of $g_2$ 
  in 
\[
\mathrm{Ker}(-\frac{1}{2}I+{W}_\Omega)_{|C^{m,\alpha}(\partial\Omega)}
=\mathrm{Ker}(-\frac{1}{2}I+{W}_\Omega)_{|C^{\max\{m,1\},\alpha}(\partial\Omega)}\,,
\]
implies  that
 there exist  a function $\mu$ in ${\color{black}\biggl(	}
 \mathrm{Ker}(-\frac{1}{2}I+{W}_\Omega^t)_{|C^{\max\{m,1\}-1,\alpha}(\partial\Omega)}
 {\color{black}\biggr)_0	}$  and a constant function $\rho$ on $\partial\Omega$ such that   
\[
V_\Omega[\mu]+\rho=v_\Omega[\mu]_{|\partial\Omega}+\rho=g_2 
\]
(cf.~Theorem \ref{thm:Ker-I+W50a} and Proposition \ref{prop:Ker-I+W2} (ii)). Next we note that the function $v^-_\Omega[\mu]$ is harmonic in $\Omega^-$, harmonic at infinity, and equals $g_2-\rho$ on $\partial\Omega$ (we note that $\int_{\partial\Omega} \mu\,d\sigma=0$).   Moreover, a classical regularity  result for the single layer potential implies that $v^-_\Omega[\mu]$ belongs to $C^{\max\{m,1\},\alpha}_{b}(\overline{\Omega^-})$ (cf.~\textit{e.g.}, Proposition \ref{prop:sldinfty}, Theorem \ref{thm:slayh}). 
   Then, if we still denote by $\rho$ the constant extension of $\rho$ to $\overline{\Omega^-}$, we have that
\[
u\equiv w_\Omega^-[\psi]+v^-_\Omega[\mu]+\rho
\]
is a solution of  the exterior Dirichlet problem (\ref{thm:existenceO-00}) and belongs to $C^{m,\alpha}_{b}(\overline{\Omega^-})$.

Both for $n\geq 3$ and for $n=2$, the uniqueness of the solution $u$ is a consequence of the uniqueness theorem for the solutions of an exterior Dirichlet problem for harmonic functions that are harmonic at infinity (cf.~\textit{e.g.}, \cite[Prop.~6.1]{DaLaMu21}).\hfill  $\Box$ 

\vspace{\baselineskip}

{\color{black}

Finally, we prove a technical approximation lemma in H\"{o}lder spaces that we have exploited in the paper.
\begin{lemma}\label{lem:0a1aap}
 Let  $\alpha\in ]0,1]$. Let $\Omega$ be a bounded open  subset of ${\mathbb{R}}^{n}$ of class $C^{1,\alpha}$. If $f\in C^{0,\alpha}(\partial\Omega)$, then there exists a sequence $\{f_j\}_{j\in {\mathbb{N}}}$ in $C^{1,\alpha}(\partial\Omega)$   that converges to $f$ in the $ C^{0,\beta}(\partial\Omega)$-norm for all $\beta\in]0,\alpha[$ and that is bounded in the $ C^{0,\alpha}(\partial\Omega)$-norm.
\end{lemma}
{\bf Proof.} Let $r\in ]0,+\infty[$ be such that $\overline{\Omega}\subseteq {\mathbb{B}}_n(0,r-1)$.  Since $f\in C^{0,\alpha}(\partial\Omega)$ and $\Omega$ is of class $C^{1,\alpha}$, there exists $\tilde{f}\in C^{0,\alpha} (\overline{{\mathbb{B}}_n(0,r)})$ such that the support of $\tilde{f}$ is contained in ${\mathbb{B}}_n(0,r-1)$ and $f=\tilde{f}_{|\partial\Omega}$ (cf.~\textit{e.g.}, \cite[Thm.~2.85]{DaLaMu21}). Then we take the convolution with a standard family of mollifiers $\{\eta_\epsilon\}_{\epsilon\in]0,+\infty[}$  and we set
\[
\tilde{f}_j(x)\equiv\int_{{\mathbb{R}}^n}\tilde{f}(x-y)\eta_{2^{-j}}(y)\,dy\qquad\forall x\in {\mathbb{R}}^n
\]
for all $j\in {\mathbb{N}}$ (cf.~\textit{e.g.}, \cite[Thm.~2.72]{DaLaMu21}). By known properties of the convolution, the function $\tilde{f}_j$ is of class $C^\infty$ in ${\mathbb{R}}^n$ and we have  
\[
\lim_{j\to\infty}\tilde{f}_j=\tilde{f}\qquad\text{uniformly\ in}\  \overline{{\mathbb{B}}_n(0,r)} \,.
 \]
Then the Young inequality for the convolution and the 
elementary inequality
\begin{eqnarray*}
\lefteqn{
|\tilde{f}_j(x_1)-\tilde{f}_j(x_2)| =
\biggl|\int_{{\mathbb{R}}^n}\tilde{f} (x_1-y)\eta_{2^{-j}}(y)\,dy
-
\int_{{\mathbb{R}}^n}\tilde{f} (x_2-y)\eta_{2^{-j}}(y)\,dy\biggr|
}
\\ \nonumber
&&\qquad
\leq 
|\tilde{f} :\, \overline{{\mathbb{B}}_n(0,r)}|_\alpha |x_1-x_2|^\alpha\int_{{\mathbb{R}}^n}\eta_{2^{-l}}(y)\,dy=
|\tilde{f} :\, \overline{{\mathbb{B}}_n(0,r)}|_\alpha |x_1-x_2|^\alpha
\end{eqnarray*}
for all $x_1$, $x_2\in \overline{{\mathbb{B}}_n(0,r)}$,  $j\in {\mathbb{N}}$, imply that 
\[
\sup_{j\in {\mathbb{N}}}\|\tilde{f}_{j}\|_{C^{0,\alpha}(\overline{{\mathbb{B}}_n(0,r)})}<+\infty
 \,,
\]
(cf.~\textit{e.g.}, Folland \cite[8.7]{Fo95}). Here  $|\tilde{f} :\, \overline{{\mathbb{B}}_n(0,r)}|_\alpha$ denotes the $\alpha$-H\"{o}lder constant of $\tilde{f}$ in $\overline{{\mathbb{B}}_n(0,r)}$. Thus we can take $f_j\equiv \tilde{f}_{j|\partial\Omega}$   and we have $f_j\in C^{1,\alpha}(\partial\Omega)$ 
for all $j\in {\mathbb{N}}$, $\sup_{j\in {\mathbb{N}}}\|f_j\|_{C^{0,\alpha}(\partial\Omega)}$ is finite and  the compactness of the embedding of 
$C^{0,\alpha}(\partial\Omega)$ into $C^{0,\beta}(\partial\Omega)$ implies that 
$\{f_j\}_{j\in {\mathbb{N}}}$ converges to $f$ in $C^{0,\beta}(\partial\Omega)$ for all $\beta\in]0,\alpha[$.\hfill  $\Box$ 

\vspace{\baselineskip}
}

  \noindent
{\bf Statements and Declarations}\\

 \noindent
{\bf Competing interests:} This paper does not have any  conflict of interest or competing interest.

 \noindent
{\bf Acknowledgement.} The author  is indebted to Prof.~Cristoph Schwab for pointing out references \cite{AzKe95}, \cite[Chapt.~II, \S 6]{LiMa68}, \cite{RoSe69} during a short visit at the University of Padova {\color{black}
and to Prof.~Paolo Musolino for a number of comments on the paper.  }

The author  acknowledges  the support of the Research  Project GNAMPA-INdAM   $\text{CUP}\_$E53C22001930001 `Operatori differenziali e integrali in geometria spettrale' and   of the Project funded by the European Union – Next Generation EU under the National Recovery and Resilience Plan (NRRP), Mission 4 Component 2 Investment 1.1 - Call for tender PRIN 2022 No. 104 of February, 2 2022 of Italian Ministry of University and Research; Project code: 2022SENJZ3 (subject area: PE - Physical Sciences and Engineering) ``Perturbation problems and asymptotics for elliptic differential equations: variational and potential theoretic methods''.


\begin{thebibliography}{11}


\bibitem{Arte15}
W. Arendt , A. F. M. ter Elst. {\em The Dirichlet-to-Neumann Operator on Exterior
Domains}. Potential Anal.,  {\bf 43} (2015),   313--340.

\bibitem{AzKe95}
A. K.~Aziz and R.B.~Kellogg, {\em On homeomorphisms for an elliptic equation in domains with corners},   Differential Integral Equations, {\bf  8} (1995),  333--352.


\bibitem{BoGr13}
U.~Bottazzini and J.~Gray, {\em Hidden harmony – geometric fantasies. The rise of complex function theory}. New York, NY: Springer 2013.


\bibitem{BrDaLu23}
R.~Bramati,  M.~Dalla Riva and  B.~Luczak. {\em Continuous harmonic functions on a ball that are not in $H^s$ for $s>1/2$}. Preprint 2023. https://arxiv.org/abs/2203.04744

 
\bibitem{Co88}
M.~Costabel, {\em  Boundary integral operators on Lipschitz domains: elementary results}, 
SIAM J. Math. Anal. {\bf 19} (1988),  613--626.

 



\bibitem{Da13}
M.~Dalla Riva, {\em A family of fundamental solutions of elliptic partial differential operators with real constant coefficients}. Integral Equations Operator Theory, {\bf 76} (2013),  1--23.

 

\bibitem{DaLaMu21}
M.~Dalla Riva, M.~Lanza de Cristoforis, and P.~Musolino. {\em Singularly Perturbed Boundary Value Problems. A Functional Analytic Approach}.  Springer, Cham,  2021.

\bibitem{DaMi15}
M.~Dalla Riva and G.~Mishuris.  {\em Existence results for a nonlinear transmission problem}. 
J. Math. Anal. Appl. 430 (2015), no. 2, 718–741.

\bibitem{DaMoMu13} 
 M.~Dalla Riva, J.~Morais, and P.~Musolino, {\em A family of fundamental solutions of elliptic partial differential operators with quaternion constant coefficients}. Math.~Methods Appl.~Sci.,  {\bf 36} (2013),  1569--1582.
 
 \bibitem{DoLa17}
F.~Dondi and M.~Lanza de Cristoforis, {\em  Regularizing properties
 of the double layer potential
of  second order elliptic differential operators},   Mem. Differ. Equ.
Math. Phys. 71 (2017), 69--110.

\bibitem{Fo95}
G.~B. Folland.
  {\em Introduction to partial differential equations}.
 Princeton University Press, Princeton, NJ, second edition, 1995.


 

\bibitem{Ha1906}
J.~Hadamard, {\em Sur le principe de Dirichlet}, Bull. Soc. Math. France, 34, 135--138, 1906.

 
\bibitem{KhPuSh07}
D.~Khavinson, M.~Putinar and H.S.~Shapiro,  {\em  Poincar\'{e}'s Variational Problem in Potential Theory}
Arch. Rational Mech. Anal. 185 (2007) 143--184.

\bibitem{Kr14}
R.~Kress. {\em Linear integral equations}, volume~82 of {\em Applied
  Mathematical Sciences}.
 Springer, New York, third edition, 2014.


   
  \bibitem{La24}
 M.~Lanza de Cristoforis, {\em
A survey on the boundary behavior of the double layer potential in Schauder spaces in the frame of an abstract approach}, to appear in Exact and Approximate Solutions for Mathematical Models in Science and Engineering,
C. Constanda, P. Harris, B. Bodmann (eds.), Springer, Cham, 2024.	 

\bibitem{LaMu18}
M.~Lanza de Cristoforis and  P.~Musolino, {\em Two-parameter anisotropic homogenization
for a Dirichlet problem for the Poisson equation
in an unbounded periodically perforated domain. 
A functional analytic approach},    Math.~Nachr. 291 (2018), no. 8-9, 1310--1341.

\bibitem{LiMa68}
J. L. Lions and E. Magenes. {\em Problemes Aux Limites Non-Homogenes et Applications}, Vol. 1,  Dunod,  Paris, 1968. 

  
\bibitem{MaSh98}  
  V. Maz’ya and T. Shaposhnikova, {\em Jacques Hadamard, a universal mathematician}. Providence, RI: American Mathematical Society;   London Mathematical Society, 1998.

\bibitem{McL00}
W.~McLean, {\em Strongly elliptic systems and boundary integral equations},  Cambridge University Press, Cambridge, 2000.



\bibitem{Mi11}
S.~Mikhailov, {\em Traces, extensions and co-normal derivatives for elliptic systems on Lipschitz domains}, J.~Math.~Anal.~Appl., {\bf 378} (2011), 324--342. 

\bibitem{Mi65}
C.~Miranda, {\em Sulle propriet\`{a} di regolarit\`{a} di certe 
trasformazioni integrali}, Atti Accad. Naz.   
Lincei Mem. Cl. Sci. Fis. Mat. Natur. Sez I, {\bf 7} (1965), 303--336. 

 

 
\bibitem{MitMitMit22}
D.~Mitrea, I.~Mitrea and M.~Mitrea, {\em Geometric harmonic analysis I -- a sharp divergence theorem with nontangential pointwise traces.} Developments in Mathematics, 72. Springer, Cham, 2022.   

\bibitem{MitTa99}
M.~Mitrea and M.~Taylor, {\em Boundary Layer Methods for Lipschitz Domains in Riemannian Manifolds}
Journal of Functional Analysis, {\bf  163}, 181--251 (1999).
 

\bibitem{Ne12}
J.~Ne\v{c}as. {\em
Direct methods in the theory of elliptic equations}.
Translated from the 1967 French original by Gerard Tronel and Alois Kufner. Springer Monogr. Math.
Springer, Heidelberg, 2012.  

\bibitem{NePl73}
J.C. Nedelec and J. Planchard. {\em Une m\'{e}thode variationnelle d'\'{e}l\'{e}ments
finis pour la r\'{e}solution num\'{e}rique d'un
probl\`{e}me ext\'{e}rieur dans $R^3$}. 
Revue fran\c{c}aise d'automatique informatique recherche op\'{e}rationnelle. Analyse num\'{e}rique, tome 7, n$^o$ R3 (1973), p. 105--129.


\bibitem{Pl11}
J.~Plemelj. {\em PotentialtheoretischeUntersuchungen},Teubner,Leipzig,1911.

\bibitem{Pr1871}
F.E.~Prym, {\em Zur Integration der Differentialgleichung $\frac{\partial^2 u}{\partial x^2}+\frac{\partial^2 u}{\partial y^2}=0$}, J. Reine Angew. Math., 
73, 340--364, 1871.

\bibitem{RoSe69}
Ja A. R\u{o}itberg and Z. G S\v{e}ftel',  {\em A theorem on homeomorphisms for elliptic systems and its
applications}, Mat. USSR Sbornik, {\bf 7} (1969), 439--465.

\bibitem{Sc31}
J.~Schauder, {\em Potentialtheoretische {U}ntersuchungen}.
  Math. Z., 33(1):602--640, 1931.


 \bibitem{We67}
W.~Wendland.
{\em  Die {F}redholmsche {A}lternative f\"{u}r {O}peratoren, die
  bez\"{u}glich eines bilinearen {F}unktionals adjungiert sind}.
  Math. Z., 101:61--64, 1967.

\bibitem{We70}
W.~{Wendland}.
{\em Bemerkungen \"uber die Fredholmschen S\"atze.}
 Methoden Verfahren Math. Phys. 3, B.I.-Hochschulskripten 722/722a,
  141-176 (1970)., 1970.

\bibitem{Wi93}
M.~Wiegner, {\em Schauder estimates for boundary layer potentials}, Math. Methods Appl.~Sci., {\bf 16}   (1993),   877--894.
 


\end{thebibliography}
\end{document}